\documentclass[12pt]{book}

\usepackage{a4,fancyhdr}
\renewcommand{\chaptermark}[1]%
                 {\markboth{#1}{}}
\renewcommand{\sectionmark}[1]%
                 {\markright{\thesection\ #1}}
\usepackage{makeidx}     
\usepackage{graphicx}    
\usepackage{multicol}    

\usepackage{dsfont}
\usepackage{eufrak}
\DeclareMathAlphabet{\varmathds}{U}{dsss}{m}{n}

\begin{document}

\pagenumbering{roman}
\setcounter{page}{1}
\pagenumbering{arabic}

\pagenumbering{roman}
\setcounter{page}{1}
\pagenumbering{arabic}

\count1=1
\centerline{\bf Structure of Tate-Shafarevich groups of }

\centerline{\bf  elliptic curves over global function fields}

\vskip0.2in
\centerline{\bf M.L. Brown}
\vskip0.2in
\centerline{ Institut Fourier, B.P. 74, 38402 Saint Martin d'H\`eres, France}

\vskip0.4in \centerline{\bf Shortened Title} 

\centerline{Tate-Shafarevich Groups}

\vskip0.4in
\centerline{\bf Abstract}

The structure of the Tate-Shafarevich groups of a class of elliptic curves 
over global function fields is determined. These are known to be finite abelian 
groups from the monograph [1] and hence they are direct sums of finite cyclic 
groups where the orders of these cyclic components are invariants of the
Tate-Shafarevich group. This decomposition of the Tate-Shafarevich groups into direct sums of 
finite cyclic groups depends on the behaviour of Drinfeld-Heegner points on these
elliptic curves. These are points  analogous to Heegner points on elliptic curves over the rational numbers.

\vskip0.2in
\vskip0.4in
\centerline{\bf AMS 2012 Mathematics Subject Classification Numbers}

\centerline{11G05, 11G09, 11G20, 11G40, 14G10, 14G17, 14G25, 14H52}
\vskip0.4in

\centerline{\bf Keywords: Elliptic curves, Tate-Shafarevich groups, Function fields}
\vskip0.4in
\centerline{\bf Table of Contents}
\vskip0.2in

\noindent  Chapter 1. Preliminaries

1.1.  Introduction

1.2. Global fields of positive characteristic

1.3. Orders in imaginary quadratic fields

1.4. Ring class fields

1.5. Elliptic curves over global fields of positive characteristic

1.6. The Drinfeld modular curve $X_0^{\rm Drin}(I)$

1.7. Analogue for $F$ of the Shimura-Taniyama-Weil conjecture

1.8. Drinfeld-Heegner points

1.9. Groups and cohomology 

1.10. Torsion on elliptic curves $E/F$

1.11. Igusa's theorem

1.12 Consequences of Igusa's theorem

\noindent Chapter 2. Local duality, Cassels pairings, Tate-Shafarevich groups

2.1. Local duality of elliptic curves

2.2. Selmer groups and  Tate-Shafarevich groups

2.3. The Cassels pairing

\noindent Chapter 3. The cohomology classes $\gamma_n(c),\delta_n(c)$

3.1. The set $\cal P$ of prime numbers

3.2. Frobenius elements and the set $\Lambda(n)$ of divisors

3.3. A refined Hasse principle for finite group schemes

3.4. Drinfeld-Heegner points and the cohomology classes $\gamma_n(c)$, $\delta_n(c)$

\noindent Chapter 4. Structure of the Tate-Shafarevich group and the Selmer group

4.1. Statement of the main theorems

4.2. Cochains for the the cohomology classes
$\gamma_n(c), \delta_n(c)$

4.3. Points $P_c$  defined over local fields

4.4. The map $\chi_z$  

4.5. Localizations of the classes $\gamma_n(c)$ and $\delta_n(c)$

4.6. The Cassels pairing with a class $\delta_n(c)$

\noindent Chapter 5. Construction of cohomology classes and proofs of the
main theorems

5.1. $M_r$ is finite for some $r$

5.2. A class $\gamma_n(c)$ in the Selmer group

5.3. Proof of  Theorem 4.1.15

5.4. Proof of Theorems 4.1.9 and 4.1.13  

5.5. Proof of Theorems 1.1.1 and 4.1.14

5.6. Generators of  Tate-Shafarevich groups

\vskip0.4in
\noindent{\bf  Chapter 1. Preliminaries}

\vskip0.2in

\noindent {\bf 1.\the\count1.  Introduction}
\advance\count1 by 1

\vskip0.2in
Let $F$ be a global field of positive characteristic $p>0$. Let $E/F$ be an elliptic curve with an origin, that is to say a 1-dimensional abelian variety.

In [1]  it is shown that for  a class of these elliptic curves $E/F$, the Tate-Shafarevich group ${\textstyle\coprod\!\!\!\coprod}(E/F)$  is finite and for prime numbers $l$ belonging to a set $\cal S$ of prime numbers  given by 
arithmetic   conditions then the $l$-primary
component ${\textstyle\coprod\!\!\!\coprod}(E/F)_{l^\infty}$ has order which is  
explicitly bounded.

In this paper, we determine the structure of 
the finite abelian group
${\textstyle\coprod\!\!\!\coprod}(E/F)_{l^\infty}$ for the same class of elliptic 
curves and for all 
$l$ in the same set of prime numbers $\cal S$ (in the notation of Theorem 4.1.9 below $\cal S$ is the set $\cal P$ with the exclusion of the prime divisors of the order of the 
Picard group Pic$(A)$). 
We also determine the structure of the Selmer groups of the elliptic curves in question.

Let $E/F$ be an elliptic curve and let $K$ be an imaginary quadratic extension 
of $F$ with respect to the place $\infty$ of $F$ (see \S1.2). Let Spec $A$ be the non-singular affine curve with function field $F$ and whose point at infinity is 
$\infty$ (see \S1.2). Assume that $E,K,\infty$ satisfy the hypotheses (a), (b), (c) of (4.1.1) in chapter 4. This provides  an infinite  set of prime numbers $\cal P$ of positive Dirichlet density and 
defined by arithmetic conditions (see \S3.1). Indeed, $\cal P$ contains all except
finitely many prime numbers $l\in {\mathds Z}$ of the form $2^sn+1$ where 
$s\geq 1$ and $n$ is odd such that $q$ is a $2^s$th power non-residue modulo $l$  where $q$ is the order of the 
exact finite field of constants of $F$.

Fix a prime number 
$l\in \cal P$. There are sets of divisors
$\Lambda^r(n)$, relative to $l$,  on $F$ for all integers $n\geq 1$ such that each divisor 
in $\Lambda^r(n)$ is a sum of $r$ distinct prime divisors and 
there is  a decreasing filtration on $\Lambda^r(1)$. 
$$\Lambda^r=\Lambda^r(1)\supseteq\Lambda^r(2)\supseteq\ldots$$
For any divisor
$c\in\Lambda^r(n)$ there is a corresponding Drinfeld-Heegner point 
$P_c$ of $E(K[c])$, the group of $K[c]$-rational points of $E$, where $K[c]$ is the ring class field of $K$ with conductor 
$c$ (see \S1.4 and (3.4.9)).

On $E(K[c])$  there is the decreasing $l$-adic filtration
$$E(K[c])\supseteq lE(K[c])\supseteq l^2E(K[c])\supseteq\ldots$$
\nopagebreak Define
\nopagebreak$$M_r=\min_{c\in \Lambda^r}(\max (n \in {\mathds N}\ \vert  P_c\in l^nE(K[c])))
{\rm \ \ for \ all \ integers \ }r\in {\mathds N}.$$

If the point $P_0\in E(K)$ has infinite order in the group of 
$K$-rational points $E(K)$ then it can be shown that $M_0,M_1,\ldots$ is a decreasing
sequence of non-negative integers (see Lemma 5.1.4). One of the main results
of this paper is the following.
\vskip0.2in
\noindent{\bf 1.1.1. Theorem.} {\sl Suppose that $P_0$ has infinite order in $E(K)$, the group of $K$-rational points of $E$. Let $l$ be 
a prime number in $\cal P$ and which is coprime to the order of the Picard group of the 
affine curve ${\rm Spec } \ A$. Let $\epsilon=\pm1$ be the sign in the 
functional equation of the $L$-function of $E/F$. Then the Tate-Shafarevich group ${\coprod\!\!\!\coprod}(E/F)$ of $E/F$ is finite and its $l$-primary component is given by
$${\coprod\!\!\!\!\coprod}(E/F)_{l^\infty}\cong \prod_{ (-1)^i=\epsilon\atop i\geq 0}
({\mathds Z}/l^{M_i-M_{i+1}}{\mathds Z})^2$$
where the product runs over integers $i\in {\mathds N}$ such that 
$(-1)^i=\epsilon$.}
\vskip0.2in
A similar statement holds for the Tate-Shafarevich group 
 of the elliptic curve $E\times_FK$ over $K$
(see Theorem 4.1.9) as well as the Selmer groups of these curves (corollary 4.1.19).  The main results of this paper are stated in \S4.1. 

It may be conjectured that for every global field $F$ of characteristic $>0$ there are
infinitely many non-isomorphic elliptic curves $E/F$ and infinitely many imaginary quadratic 
field extensions $K/F$ such that $E,K,\infty$ satisfy the hypotheses of this Theorem
1.1.1 and those of \S4.1.
If this conjecture holds, then the above theorem and those of \S4.1 give infinitely many non-isomorphic 
elliptic curves over a given global field 
of positive characteristic 
whose $l$-primary components of the Tate-Shafarevich group are structurally known
for infinitely many prime numbers $l$ satisfying arithmetic conditions.

The method  of this paper is related to  that of Kolyvagin's   determination 
of the structure of Tate-Shafarevich groups of a class of elliptic curves over the rational numbers (see [6] and [8]).

The proofs of the main theorems of this paper stated in \S4.1 and Theorem 1.1.1 above require many preliminary results which 
are explained in chapters 1-5.

Chapter 2 contain basics  on Tate local duality, Selmer groups, and the Cassels pairing on Tate-Shafarevich groups.
In section 3.1, the set $\cal P$ of prime numbers is defined by arithmetic conditions.
In section \S3.2, the sets $\Lambda^r(n)$ 
of divisors on the global field  $F$ are defined. Sections \S\S3.3,3.4 constructs
the cohomology classes $\gamma_n(c),\delta_n(c)$ in the cohomology of 
the elliptic curve $E/F$. 

In section \S4.1., the main results of this paper are stated which are then proved in \S\S5.3-5.5 after some further properties 
of $\gamma_n(c),\delta_n(c)$ are proved in Chapters 4 and 5. We show in particular that the cohomology classes
$\delta_n(c)$ define characters, via the Cassels pairing on ${\textstyle\coprod\!\!\!\coprod}(E/F)_{l^\infty}$ which determine the structure of this group. The method of proof of the main results in \S4.1. is by the  construction of many independent elements of the Tate-Shafarevich 
group ${\textstyle\coprod\!\!\!\coprod}(E/F)_{l^\infty}$.
 Finally, \S5.6 contains complements to the main results.

While care has been taken to minimize the number of hypotheses required for the main theorems of this paper, these hypotheses  are still numerous (see for example \S3.1, Definition 3.1.3).
The assiduous reader will have an abundance of interesting problems in their elimination.

\vskip0.4in 
\noindent {\bf 1.\the\count1.  Global fields of positive characteristic}
\advance\count1 by 1

\vskip0.2in 
The notation of this paper is mainly that of the monograph [1] and is detailed in the rest of this chapter 1. A few  differences arise, notably the sets of divisors $\Lambda(n)$,
which are  required for the more refined results of this paper. 

Let

\vskip0.2in
$k$ be a finite field of characteristic $p$ with $q=p^m$ elements;

${\overline k}$ be an algebraic closure of $k$;

$C/k$ be a smooth projective irreducible curve over $k$;

$F$ be the function field of $C$; these hypotheses imply that the finite field $k$ is 

\qquad the exact field of constants of the global field $F$;

$\Sigma_L$, for any global field $L$, be the set of all places of the field $L$;

$\infty\in \Sigma_F$ be a closed point of $C/k$;

$\kappa(z)$ be the residue field at a place $z\in \Sigma_F$ of $F$;

$F_v$ be the completion of $F$ at the place $v\in \Sigma_F$;

$F^{\rm sep}$ be the separable closure of $F$;

 $A$ be the coordinate ring $\Gamma(C\setminus\{\infty\}, {\cal O}_C)$ of the affine curve $C\setminus\{\infty\}$;

 Div$_+(A)$ be the semi-group of effective $k$-rational divisors on Spec $A$; that 

\qquad is to say Div$_+(A)$ is the semi-group of effective $k$-rational divisors on $C/k$ 

\qquad which are coprime to the place $\infty$; an element of 
Div$_+(A)$ may be written 

\qquad as a finite linear combination $\sum_i n_i z_i$
where $n_i\in \mathds N$ and $z_i$ are prime 

\qquad divisors  on Spec $A$ for all $i$;

Supp$(c)$, for  a divisor $c\in $ Div$_+(A)$, be the support of the 
divisor $c$ which  is 

\qquad the set of prime divisors with non-zero coefficient in $c$;

Pic$(A)$ be the Picard group of $A$, the group of projective  $A$-modules
of 

\qquad rank $1$;

$K$ be a separable imaginary quadratic extension field of $F$ with respect to $\infty$, 

\qquad that is to say $K$ is a quadratic extension field of $F$ in which the place $\infty$

\qquad  remains  inert;

$B$ be the integral closure of $A$ in $K$;

$\tau$ be the non-trivial element of the galois group Gal$(K/F)$.

\vskip0.4in\noindent {\bf 1.\the\count1.  Orders in imaginary quadratic field extensions}\advance\count1 by 1
\vskip0.2in 
Let $K/F$ be the imaginary quadratic field extension with respect to $\infty$ of \S1.2.

An {\it order} $O$ in $K$ with respect to $A$ is an $A$-subalgebra of $B$ whose fraction field is equal to $K$. 

There is a bijection between orders $O_c$ of $K$ with respect to $A$ and effective $k$-rational divisors $c$ in Div$_+(A)$ and it
is given by
$$c\mapsto A+BI(c)$$
where $I(c)$ is the ideal of $A$ cutting out the divisor $c$. The divisor $c$ is the {\it conductor} of the order $O_c$.

[For more details on orders in imaginary quadratic extensions, see [1, Chapter 2, \S2.2]

\vskip0.4in\noindent {\bf 1.\the\count1.  Ring class fields}\advance\count1 by 1
\vskip0.2in 
Let 

\vskip0.2in

$O_c$ be the order of $K$  with respect to $A$ and  with conductor $c$ 

\qquad where $c\in$ Div$_+(A)$ (see \S1.3);

   $A_v$, for each place $v$ of $A$, be the localisation of $A$ at $v$;

${\widehat O}_{c,v}$ be the completion of the semi-local ring $O_c\otimes_AA_v$;

$G_c=K_\infty^*\prod_v {\widehat O}_{c,v}^*$ be the subgroup of the 
id\`ele group of the global field $K$ 

\qquad whose components are the units of ${\widehat O}_{c,v}$ for all places $v\not=\infty$ of $F$ 

\qquad and $K_\infty^*$ for the 
place $v=\infty$ and where in the product $v$ runs over all 

\qquad places  of $F$;

$K[c]$, for any divisor $c\in$ Div$_+(A)$, be the ring class field with conductor $c$ 

\qquad with respect to $\infty$;  this is the finite abelian extension field of $K$ defined 

\qquad by the subgroup $G_c$ of the id\`ele group of $K$
 via the reciprocity map;

$G(c/c')$ be the Galois group of the field extension $K[c]/K[c']$ for divisors $c\geq c'$ 

\qquad of Div$_+(A)$.

\vskip0.2in
We have these properties of the decomposition of primes in ring class fields (for the proofs, see [1, Chapter 2, \S2.3.13]):
\vskip0.2in
(a) The primes ramified in $K[c]/K$ are precisely the primes in the support 

\qquad of $c$.

(b) The extension $K[c]/K$ is split completely at the place of $K$ lying 

\qquad above
$\infty$.

(c) If $z\not\in $ Supp$(c)$ then for any positive integer $n$, the galois extension 

\qquad $K[c+nz]/K[c]$ is totally ramified at all places of $K[c]$ above $z$.
\vskip0.2in [See [1, Chapter 2, \S2.3] for 
more details on ring class fields.]

\vskip0.4in 
\noindent {\bf 1.\the\count1.  Elliptic curves over global fields of positive characteristic} \advance\count1 by 1
\vskip0.2in 
 Let

\vskip0.2in

$C/k$ be a smooth projective irredicible curve over $k$ (as in \S1.2);

$X/k$ be an elliptic surface over $C$, that is to say $X/k$ is a smooth

\qquad projective irreducible surface, equipped with a morphism $f:X\to C$ 

\qquad which has a section, such that all fibres of $f$, except a finite number, are 

\qquad elliptic curves;

$E/F$ be the generic fibre of $f:X\to C$, which is an elliptic curve $E$ over $F$ 

\qquad equipped with an origin where $F$ is the function field of $C$.

\vskip0.2in

The conductor of an elliptic curve $E/F$ is an effective $k$-rational divisor on $F$
supported only at the places of bad reduction of $E$ and whose multiplicities are defined in terms of the Galois representation of Gal$(F^{\rm sep}/F)$ given by $E$.

 [See [1, \S\S B.11.1-B.11.4 for the definition of the conductor of $E/F$.]

\vskip0.4in\noindent {\bf 1.\the\count1.  The Drinfeld modular curve 
$X_0^{\rm Drin}(I)$ }\advance\count1 by 1
\vskip0.2in 
Let 
\vskip0.2in

$I$ be a non-zero 
ideal of $A$.

  $X_0^{\rm Drin}(I)$ be the curve which is the coarse moduli scheme of Drinfeld modules 

\qquad of rank $2$ for $A$ equipped with an $I$-cyclic structure (see
[1, Definition 

\qquad 2.4.2, p. 23]); this curve is compactified by a finite 
number of cusps which 

\qquad correspond to ``degenerate" Drinfeld modules.
\vskip0.2in
This curve $X_0^{\rm Drin}(I)$ is an analogue for the global field $F$
of the classical modular curve $X_0(N)$ which is the coarse moduli scheme of elliptic curves equipped with a cyclic subgroup of order $N$,
where $N\in \mathds N$.

[For more details, see [1, \S2.4].]

\vskip0.4in\noindent {\bf 1.\the\count1.  Analogue for $F$ of the Shimura-Taniyama-Weil conjecture} \advance\count1 by 1

\vskip0.2in 
Let $E/F$ be an elliptic curve with split (Tate) multiplicative reduction at $\infty$. Let $I$ be the non-zero ideal of the ring $A$ which is the conductor, without the place at $\infty$, of the elliptic curve $E/F$.

According to the work of Drinfeld on the Langlands conjecture, there is a finite 
surjective morphism of curves over $F$
$$X_0^{\rm Drin}(I)\to E.$$
This result is an analogue for the global field $F$ of  the Shimura-Taniyama-Weil conjecture proved by Wiles for semi-stable elliptic curves over 
the rational numbers.

[For more details see [1, \S4.7] and [1, Appendix B].]

\vskip0.4in\noindent {\bf 1.\the\count1. Drinfeld-Heegner points} 
\advance\count1 by 1

\vskip0.2in 

Let $K$ be an imaginary quadratic extension field of $F$ with respect to
$\infty$ (as in \S1.2) and let $I$ be a non-zero ideal of $A$.

Let $D$ be a rank $2$ Drinfeld module for $A$ with complex multiplication by an  order $\cal O$ of the field $K$ with respect to $A$, that is to say 
$\cal O$ is a subring of $K$ which is integral over $A$. Let $Z$ be an 
$I$-cyclic subgroup of $D$. Then the pair $(D,Z)$ represents a non-cuspidal point
of the modular curve $X_0^{\rm Drin}(I)$. Such points $(D,Z)$ exist if the 
prime divisors in the support of $I$ split completely in the field 
extension $K/F$.

If the quotient Drinfeld module $D/Z$ has the same ring of endomorphisms $\cal O$ as $D$ then the point $(D,Z)$ on $X_0^{\rm Drin}(I)$ is called a {\it Drinfeld-Heegner point}.

If $f:X_0^{\rm Drin}(I)\to E$ is a finite morphism of curves where $E/F$ is an elliptic curve (see \S\S1.5,1.6,1.7), then the point $f(D,Z)$ of the elliptic curve $E$ is also called a Drinfeld-Heegner point.

The Drinfeld-Heegner points $(D,Z)$ and $f(D,Z)$ are rational over the ring class field
$K[c]$ where $c$ is the conductor of the order $\cal O$ of $K$ relative to $A$. 

[See [1, \S\S2.2,2.3] or \S3.4 below for more details.]

\vskip0.4in\noindent {\bf 1.\the\count1. Groups and cohomology}\advance\count1 by 1

\vskip0.2in 

 If $G$ is a discrete abelian group denote by 
\vskip0.2in
$G_m$  the kernel of multiplication by the integer $m \in \mathds N$ on $G$; 

$_mG$ the cokernel $G/mG$ of multiplication by the integer $m\in \mathds N$;

$\vert G\vert$  the order of the group $G$ which is either a positive
integer or $+\infty$;

ord$(g)$ the order of an element $g\in G$ which is the cardinality of the subgroup 

\qquad generated by $g$;

exp$(G)$ the exponent of $G$ which is maximum order of an element of $G$;

 ${\widehat {G }}$ the Pontrjagin dual of $G$ namely the topological group Hom$(G, {\mathds Q}/{\mathds Z})$;

\vskip0.2in

If $G$ is a finite abelian group then ${\widehat G}$ may be identified with the group 
of $1$-dimensional complex characters of $G$,  that is to say  the 
group of  homomorphisms ${\rm Hom}(G, {\mathds C}^*).$ A character of a 
finite abelian group is always assumed to be irreducible.

If $F'/F$ is a galois extension of a global field $F$  and if $z\in 
\Sigma_F$ is a prime divisor of $F$ unramified in $F'$ we denote by Frob$(z)$,  or Frob$_{F'/F}(z)$, the conjugacy class in Gal$(F'/F)$ of Frobenius substitutions associated to $z$.

If $L$ is a field, we shall write
$H^i(L,M)$ for the Galois cohomology group $H^i({\rm Gal}(L^{\rm sep}/L),M)$,
where $L^{\rm sep}$ is the separable closure of $L$. If $L'/L$ is a 
finite galois field extension, we write $H^i(L'/L,M)$
for $H^i({\rm Gal}(L'/L),M)$.

If $F$ is a global field and $z$ is a place of $F$ then the restriction, or localization, of a class $c\in H^i(L,M)$ is written $c_z\in H^i(F_z,M)$
where $F_z$ is the completion of $F$ at $z$.

\vskip0.4in\noindent {\bf 1.\the\count1. Torsion on elliptic curves $E/F$} 

\vskip0.2in  
Let
$E$ be an elliptic curve over a global field $F$ of positive characteristic $p>0$
as in \S1.5. Let $K$ be an imaginary quadratic extension field of $F$ with respect
to the place $\infty$ of $F$. As in \S1.4,  let $K[c]$ be the ring class field over $K$
with conductor $c\in $ Div$_+(A)$.
Put
$$K[A]=\bigcup_{c\in {\rm Div}_+(A)} K[c].$$
That is to say $K[A]$ is a field which is the union of all the ring class
fields $K[c]$ in some algebraic closure of $K$.

Let $S$ be a subset of ${\mathds N}^*$ such that if $n_1\in S$ and $n_0$ is any 
divisor of $n_1$ then $n_0\in S$. A {\it quasi-group}
 $\{ G_n\}_{n\in S}$ relative to $S$ is a family of abelian groups $G_n$ indexed by the elements
$n$ of $S$ such that $nG_n=0$ and if $n_0,n_1\in S$ are elements where
$n_0$ divides $n_1$ there is a group homomorphism $f_{n_1n_0}:G_{n_1}\to 
G_{n_0}$
satisfying the compatibility condition that if $n_0,n_1,n_2\in S$ and $n_0$ divides $n_1$ and 
$n_1$ divides $n_2$  then 
$f_{n_2n_0}= f_{n_1n_0}\circ f_{n_2n_1}$. [See [1, Chapter 7, \S7.1, p. 
330].]

\vskip0.2in\noindent {\bf 1.\the\count1.1 Proposition.} ([1, Proposition 7.3.8.])
{\sl   The quasi-group
$$\{ E(K[A])_n\}_{n\in \mathds N}$$
is trivial, that is to say the order of the torsion group $E(K[A])_n$ is bounded
 independently of $n$ and there is a finite set $\cal E$ of prime numbers
such that for all integers $n$ prime to all elements of $\cal E$ we have
$$E(K[A])_n\cong 0.\ \ \ {\sqcap \!\!\!\!\sqcup}$$}

\vskip0.2in\noindent {\bf 1.\the\count1.2. Proposition.} ([1, Proposition 7.14.2])
{\sl Let $\cal E$ be the finite set of prime numbers of Proposition 1.\the\count1.1. For any divisor $c$ of {\rm Div}$_+(A)$, the restriction homomorphism 
$$H^1(K,E_n)\to H^1(K[c],E_n)^{{\cal G}_c}$$
is an isomorphism for all integers $n$ prime to $\cal E$ where ${\cal G}_c={\rm Gal}(K[c]/K)$.}
\vskip0.2in
This follows from Proposition 1.\the\count1.1 and the Hochschild-Serre spectral sequence
$$H^i({\cal G}_c, H^j(K[c], E_n(K^{\rm sep})))\Rightarrow H^{i+j}
(K, E_n(K^{\rm sep}))$$(more details are given in  [1, Proposition 7.14.2]).${\sqcap \!\!\!\!\sqcup}$

\advance\count1 by 1
\vskip0.4in 
\noindent { \bf 1.\the\count1. Igusa's theorem}
\vskip0.2in 
This section is a summary of  the results of
Igusa for the Galois action on torsion
points of elliptic curves over global fields of positive characteristic.

\vskip0.2in
\noindent (1.\the\count1.1)
As in \S1.2, let $F$ be a global field of positive characteristic $p$ where
$k$ is the exact field of constants of $F$ and  let $E/F$ be an
elliptic curve. Let

\vskip0.2in

$G=$ Gal$(F^{\rm sep}/F)$, where $F^{\rm sep}$ is the separable closure of $F$;

$n$ be an integer prime to $p$;

$E_n$ be the finite 
$F$-group scheme of $n$-torsion points of $E$;

$E_\infty$ be the torsion subgroup of $E(F^{\rm sep})$ of order prime to $p$.
\vskip0.2in

\noindent The elliptic curve $E/F$ is said to be \label{isotrivial} {\it isotrivial} if there is a finite galois extension
field $F'$ of $F$ such that the curve $E\times_FF'$ is definable over a finite subfield of 
$F'$.

\vskip0.2in
\noindent (1.\the\count1.2)
The action of the galois group $G$ on $E_n$ provides  a homomophism
$$\rho_n:G\to {\rm Aut}(E_n)\cong {\rm GL}_2({\mathds Z}/n{\mathds Z}).$$
The determinant
$${\rm det}:{\rm Aut}(E_n)\to ({\mathds Z}/n{\mathds Z})^*$$
induces a homomorphism
$$G\to ({\mathds Z}/n{\mathds Z})^*.$$
Let $H_n$ be the subgroup of $({\mathds Z}/n{\mathds Z})^*$ generated by the powers
of $q=\vert k\vert$ modulo $n$. Then $H_n$ is naturally isomorphic to the
Galois group of the field of $n$th roots of unity over $k$. Let $\Gamma_n$ be
the subgroup of GL$_2({\mathds Z}/n{\mathds Z})$ defined by the exact sequence of
finite groups, where det is the restriction to $\Gamma_n$ of the determinant homomorphism on GL$_2({\mathds Z}/n{\mathds Z})$,
$$\matrix{&&&&&{\rm det}&&&\cr
0&\to& {\rm SL}_2({\mathds Z}/n{\mathds Z})&\to &\Gamma_n&\to& H_n&\to
&0.\cr}\leqno{(1.\the\count1.3)}$$ 
\vskip0.2in
\noindent (1.\the\count1.4)
 Passing to the projective limit of the previous exact sequence over all integers $n$
prime to $p$ we obtain the exact sequence of profinite groups
$$0\quad \longrightarrow\quad
 {\rm SL}_2({\widehat{\mathds Z}}^{(p)})\quad \longrightarrow\quad {\widehat\Gamma}\quad \longrightarrow\quad{\widehat H}\quad \longrightarrow\quad 0$$
where $\widehat H$ is the subgroup of ${\widehat{\mathds Z}}^{(p)*}$ topologically
generated by $q$,  where
$${\widehat{\mathds Z}}^{(p)} =\prod_{l\not=p}{\mathds Z}_l$$
is the profinite prime-to-$p$ completion of $\mathds Z$, and
$\widehat{\Gamma}$ is a closed subgroup of ${\rm GL}_2({\widehat{\mathds
Z}}^{(p)})$. 

\vskip0.2in  
\noindent (1.\the\count1.5)
Passing to the projective limit of the exact sequence $(1.\the\count1.3)$ where $n$
runs over all powers of a prime number  $l$ where $l\not=p$, we obtain the
exact sequence
$$0\quad \longrightarrow\quad
 {\rm SL}_2({{\mathds Z}}_l)\quad \longrightarrow\quad {\widehat\Gamma}_l\quad \longrightarrow\quad{\widehat H}_l\quad \longrightarrow\quad 0.$$

\vskip0.2in
\noindent {\bf  1.\the\count1.6. Theorem.} (Igusa [I]). {\sl Suppose that $E/F$ is not
isotrivial. Then the profinite group ${\rm Gal}(F(E_\infty)/F)$ is an open
subgroup of $\widehat \Gamma$.\qquad ${\sqcap \!\!\!\!\sqcup}$}

\vskip0.2in This result has the following consequence.
\vskip0.2in
\noindent {\bf  1.\the\count1.7. Theorem.} {\sl  Suppose that $E/F$ is not
isotrivial. Then  for all prime numbers $l\not=p$ the profinite group ${\rm
Gal}(F(E_{l^\infty})/F)$ is an open subgroup of ${\widehat \Gamma}_l$ and is
equal to  ${\widehat \Gamma}_l$ for all but finitely many $l$.\qquad ${\sqcap \!\!\!\!\sqcup}$}
\vskip0.2in
\noindent {\it  1.\the\count1.8. Remarks.} (1) Suppose that the curve $E/F$ is isotrivial. Then it is
easy to show that 
the group ${\rm Gal}(F(E_\infty)/F)$  is an extension of a finite group by the abelian profinite
group ${\widehat{\mathds Z}}^{(p)}$.

\noindent (2)  Let $E$ be an elliptic curve defined over a number field $L$. The galois action on the torsion points of $E/L$ is known and depends 
principally  on whether or not $E$ has complex multiplication. 

[See  [S5]  and [S5, \S 4.5] for both the cases of complex multiplication and   without complex multiplication. See also [1, Chapter 7, Remarks 7.2.8] for more details.]

\advance\count1 by 1

\vskip0.4in
 
\noindent {\bf 1.\the\count1. Consequences of Igusa's theorem}
 \vskip0.2in  For a finite group $G$ and a  ${\mathds
Z}[G]$-module $M$, let   ${ H}^i(G,M)$ denote the standard cohomology
groups of $G$ acting on $M$  (see [1, \S5.6], see also  [9, Chap. 1] for the Tate modified cohomology groups).\vskip0.2in

 \noindent {\bf  1.\the\count1.1. Proposition.} {\sl  Let $E/F$ be an elliptic curve and 
let ${\mathds N}^{(p)}$ be the set of positive
integers coprime to $p$, where $p$ is the characteristic of $F$.
Write $G_n$ for the group ${\rm Gal}(F(E_n)/F)$.}

\noindent (i) {\sl Let $i=0$ or $1$. Then   }
$$\{ { H}^i(G_n,E_{n})\}_{n\in {\mathds N}^{(p)}}$$
{\sl is a trivial quasi-group. }

\noindent (ii) {\sl There is a finite set $\cal N$ of prime numbers including
$p$ such that for all prime numbers $l\not\in {\cal N}$ we have
   } $$ H^i(G_{l^n},E_{l^n})=0{\textsl{ \ \ for \ all }} \ n\geq 1{\textsl{
 \ and \ all
\ }} i\geq 0.$$  

[Part (i) may be restated as:  for $i=0$ or $1$ and for all $n$ coprime to $p$, the order of the  group
${ H}^i(G_n,E_{n}) $ is bounded
 independently of $n$ and there is a finite set  of prime numbers
such that for all integers $n$ coprime to this set of prime numbers we have
${ H}^i(G_n,E_{n}) \cong 0$. For the proof of this proposition, see [1, Chapter 7, Proposition 7.3.1] ${\sqcap \!\!\!\!\sqcup}$\vskip0.2in

\vskip0.2in\noindent (1.12.2) For each prime number $l$ different from 
$p$,  
select once and for all a basis of the Tate
module $T_l(E)$ over ${\mathds Z}_l$, the $l$-adic completion of $\mathds Z$; this fixes for the rest of this paper, for every  prime number $l\not=p$, an
isomorphism of ${\rm
Gal}(F(E_{l^\infty})/F)$ with a subgroup of ${\rm GL}_2({\mathds Z}_l)$ and the two groups may then be identified with each other.

\vskip0.2in
 \noindent {\bf  1.\the\count1.3. Proposition.} {\sl Let $E/F$ be an elliptic curve
which is not isotrivial. Let $L$ be a finite extension field of $F$ in which $k$
is algebraically closed. Then there is an infinite set $S$ of prime numbers  of
positive Dirichlet density   
 such that  for all   $l\in S$ we have

\noindent {\rm (a)} the fields $F(E_{l^\infty})$ and $L$ are linearly disjoint
over $F$;

\noindent {\rm (b)} $E(L)_{l^\infty}=0$;

\noindent {\rm (c)} $\pmatrix{1&0\cr 0&-1\cr}\in {\rm Gal}(F(E_{l^\infty})/F)$.}

\vskip0.2in
\noindent [For the proof, see [1, Proposition 7.3.10].] ${\sqcap \!\!\!\!\sqcup}$
\vskip0.4in

\vfil\eject

\noindent {\bf Chapter 2. Local duality, Cassels pairings, Tate-Shafarevich groups}
\vskip0.4in
\noindent {\bf 2.1. Local duality of elliptic curves}
\vskip0.2in
This section is a  brief summary of  
  Tate-Poitou local duality for elliptic curves over a
local field. 

[For more details on local duality of abelian varieties, see [9, Chapter I and Chapter III, \S7].]

\vskip0.2in\noindent (2.1.1) Let
\vskip0.2in
$L$ be a non-archimedian complete local field;

$L^{\rm sep}$ be the separable closure of $L$;

$E/L$ be  an elliptic curve over $L$;

$n\geq 1$ be an integer coprime to   the  characteristic of $L$;

$G$ be the Galois group Gal$(L^{\rm sep}/L)$.

\vskip0.2in\noindent (2.1.2)
Let $\mu_{n}$ be the multiplicative subgroup of $L^{{\rm sep}}$ of $n$th roots
of unity. Then $\mu_{n}$ is a finite $G$-module. Let $E_{n}$ denote the
$G$-module  of $n$-torsion points of $E(L^{\rm sep})$. We have   an abelian
group isomorphism $$E_{n}\cong {{\mathds Z}\over n{\mathds Z}}\oplus {{\mathds Z}\over
n{\mathds Z}}.$$
 Denote by $\{,\}$ the Weil pairing
$$\{,\}:E_{n}\times E_{n}\to \mu_{n}.$$ This is a perfect pairing of $G$-modules. In
particular, we have an isomorphism of $G$-modules $$E_{n}\cong {\rm Hom}_G( E_{n}, \mu_{n}).$$

\vskip0.2in\noindent (2.1.3) \label{WP1} The Weil pairing  induces a cup-product pairing in Galois cohomology
$$H^1(L,E_{n})\times H^1(L,E_{n})\to H^2(L,\mu_{n}).$$
This is a non-degenerate anti-symmetric pairing of abelian groups.
By local class field theory, we have a canonical isomorphism of groups, where Br$(L)$ is the
Brauer group of $L$,
$${\rm Br}(L)\cong {{\mathds Q}\over {\mathds Z}}.$$
This induces an isomorphism
$$ H^2(L,\mu_{n})={\rm Br}(L)_{n}\cong {{\mathds Z}\over n{\mathds Z}}.$$
\eject
\vskip0.2in\noindent {\bf 2.1.4. Theorem.} (Tate-Poitou local duality). {\sl 
The cup product pairing  
$$<,>_v:H^1(L,E_{n})\times H^1(L,E_{n})\to {{\mathds Z}\over n{\mathds Z}}\leqno{(2.1.5)}\label{aa}$$
 obtained from the Weil pairing  is  an
alternating and non-degenerate pairing of  ${{\mathds Z}\over n{\mathds Z}}$-modules.}
\vskip0.2in
 [For the proof, see [9, Chapter I, Cor.
2.3].]${\sqcap \!\!\!\!\sqcup}$

 \vskip0.2in\noindent {\bf 2.1.6. Theorem.} {\sl Assume that $n$ is prime
to the residue field characteristic of $L$.}

\noindent (i)  {\sl  The
subgroup ${}_n\!E(L)$ of $ H^1(L,E_{n})$ is isotropic for the alternating
pairing $<,>_v$.} 

\noindent (ii)  {\sl  The cup product pairing $<,>_v$ on $ H^1(L,E_{n})$ 
induces a non-degenerate pairing of abelian groups, where 
${}_nE(L)=E(L)/nE(L)$,
$$[,]_v:{}_n\!E(L)\times H^1(L,E)_{n}\to {{\mathds Z}\over n{\mathds Z}}.$$
}
\vskip0.2in

 [For the  proof, see [1, Chapter 7, Theorem 7.15.6, p.403] for part (i)
and [9, Chap. I, Cor. 3.4 and Remark 3.6, Chapter III, Theorem 7.8] for 
part (ii). Note that part (ii) holds without the hypothesis that $n$ be coprime
to the characteristic of $L$, see [9, Chap. III, \S7].] ${\sqcap \!\!\!\!\sqcup}$

\vskip0.2in\noindent(2.1.7) Suppose now that $K$ is a global field of positive characteristic, $\Sigma_K$ is the set of all places of $K$, and that 
the integer
$n$ is coprime to the characteristic of $K$. Let $E/K$ be an elliptic curve.

\vskip0.2in
\noindent{\bf 2.1.8. Proposition.}
{\sl Let $c$ and $c'$ be two elements of $H^1(K,E_n)$. Denote by $c_v$ and $c'_v$ the induced elements of $H^1(K_v,E_n)$ for all places $v\in \Sigma_K$ of $K$ where $K_v$ is the completion of $K$ at $v$. Then we have
$$\sum_{v\in \Sigma_K}<c_v, c_v'>_v=0.$$}

\vskip0.2in
\noindent{\it Proof.}
The sum of the local invariants of a global class in $H^2(K,{\mathds G}_m)$ is zero.${\sqcap \!\!\!\!\sqcup}$
\vfil\eject
\vskip0.4in
\noindent {\bf 2.2.  Selmer groups and   Tate-Shafarevich groups}
\vskip0.2in

\noindent (2.2.1) Let 
\vskip0.2in
$E/F$ be an elliptic curve as in \S1.5;

 $n\in {\mathds N}$ be an integer coprime to the characteristic
$p$ of the global field $F$.

\vskip0.2in

\noindent (2.2.2)    For a place   $v$  of the  field $F$, we write $F_v$ for
the completion of $F$ at the place $v$ (as in \S1.2). 
 The exact sequence of commutative group schemes over $F$, obtained 
 from the morphism of multiplication by $n$,
$$0\longrightarrow  E_n\longrightarrow  E\ {\buildrel n\over \longrightarrow }\ E\longrightarrow  0$$
gives rise to  a commutative diagram, where the maps
res$_v$ are the restriction homomorphisms at $v$ and the rows are exact sequences
of abelian groups,
$$\matrix{0&\to & _{n}E(F)&\to & H^1(F,E(F^{\rm sep})_{n})&\to & H^1(F,E(F^{\rm
sep}))_{n}&\to & 0\cr &&\downarrow{\rm res}_v&&\downarrow{\rm res}_v&&\downarrow{\rm
res}_v&&\cr 0&\to & _{n}E(F_v)&\to & H^1(F_v,E(F^{\rm sep}_v)_{n})&\to &
H^1(F_v,E(F^{\rm sep}_v))_{n}&\to & 0\cr}$$
 As in \S1.9, 
$_nE(F)$ denotes the cokernel  
$E(F)/nE(F)$ and $E(F^{\rm sep})_n$ denotes the
$n$-torsion subgroup of $E(F^{\rm sep})$. 

\vskip0.2in

\noindent (2.2.3) The {\it Tate-Shafarevich group}
${\textstyle\coprod\!\!\!\coprod}(E/F)$ of $E/F$ is defined
as
$${\textstyle\coprod\!\!\!\coprod}(E/F)={\rm ker}\{ H^1(F,E)\to
\prod_{v\in \Sigma_F}H^1(F_v,E)\}.$$ 
This group
${\textstyle\coprod\!\!\!\coprod}(E/F)$ is known to be a torsion 
of cofinite type (see [10]).

The {\it $n$-Selmer group} is
defined as, in terms of the commutative diagram of (2.2.2) and where 
${\rm res}_v$ denotes the middle vertical homomorphism of the diagram,   
$${\rm Sel}_{n}(E/F)=\bigcap_{v\in \Sigma_F}{\rm
res}_v^{-1}(_{n}E(F_v)).$$ Therefore ${\rm Sel}_{n}(E/F)$
is  a subgroup of $H^1(F,E(F^{\rm sep})_n)$ and is a finite abelian
group. We then have the exact sequence of
torsion abelian groups from the commutative diagram of (2.2.2),
where ${\textstyle\coprod\!\!\!\coprod}(E/F)_{n}$ is the $n$-torsion
subgroup of ${\textstyle\coprod\!\!\!\coprod}(E/F)$, $$0\to \ _{n}E(F)\to {\rm Sel}_{n}(E/F)\to
{\textstyle\coprod\!\!\!\coprod}(E/F)_{n}\to 0.$$ 

\vskip0.2in

\noindent (2.2.4) Let $F'/F$ be a finite separable galois field extension.
We write 
${\textstyle\coprod\!\!\!\coprod}(E/F')$ in place of 
${\textstyle\coprod\!\!\!\coprod}(E\times_FF'/F')$ for the Tate-Shafarevich group of $E\times_FF'$ over $F'$ obtained by ground field extension from $F$ to $F'$; similarly, for
the Selmer quasi-group, we write ${\rm Sel}_{n}(E/F')$ in place of
${\rm Sel}_{n}(E\times_FF'/F')$.

\vskip0.4in
\noindent {\bf  2.2.5. Proposition.} {\sl  Let $r$ be the degree of the 
finite separable field extension   $F'/F$.  
For any torsion abelian group $\cal A$ write $\cal A_{({\rm non}\  r )} $ for the 
torsion subgroup of $\cal A$ of order coprime to $r$. The  restriction
homomorphism provides  isomorphisms for all integers $n$ coprime to $r$}
$${\textstyle\coprod\!\!\!\coprod} (E/F)_{
({\rm non}\  r )}\ \ {\buildrel {\buildrel{\rm
res}\over\cong }\over  \longrightarrow} \ \
{\textstyle\coprod\!\!\!\coprod} (E/F')_{({\rm non}\  r )}{}^{{\rm Gal}(F'/F)}$$
$${\rm Sel}_n(E/F)  \ \ 
{\buildrel {\buildrel{\rm res}\over\cong }\over  \longrightarrow} 
\ \ 
{\rm Sel}_n(E/F'){}^{{\rm Gal}(F'/F)}.$$
\noindent {\it Proof.}
%
The
definition of the Tate-Shafarevich groups  provides  a commutative diagram 
with exact rows $$\matrix{ 0&\!\!\to\!\! &\!\!
{\textstyle\coprod\!\!\!\coprod}(E/F')^{{\rm
Gal}(F'/F)}\!\! &\!\!\to\!\! & \!\!H^1(F',E)^{{\rm Gal}(F'/F)}\!\!
&\!\!\to\!\! &\!\!\displaystyle{\bigg(\prod_{v\in
\Sigma_{F'}}H^1(F'_v,E)}\bigg)^{{\rm Gal}(F'/F)}\cr &&\uparrow&&\uparrow&&\uparrow\cr
 0&\!\!\to\!\! &\!\!{\textstyle\coprod\!\!\!\coprod}(E/F)\!\!
&\!\!\to\!\! & \!\!H^1(F,E)\!\!&\!\!\to\!\!
&\!\!\displaystyle{\prod_{v\in \Sigma_{F}}H^1(F_v,E)}\cr}$$
For any place $v$ of $F$, the $F_v$-algebra $F_v\otimes_FF'$ is \'etale and is
the product of the completions of $F'$ at the places lying over $v$. The
inflation restriction sequence provides isomorphisms for any integer $s$ coprime to the order of Gal$(F'/F)$
 \begin{eqnarray*}
H^1(F,E)_{s}&\cong H^1(F',E)^{{\rm
Gal}(F'/F)}_{s}\phantom{{{\rm Gal}(F'/F)}
{\rm \  \ all\ places \ } v{\rm \ of \ } F}\\
H^1(F_v,E)_{s}&\cong H^1(F_v\otimes_FF',E)_{s}^{{\rm Gal}(F'/F)}
{\rm \ for \ all\ places \ } v{\rm \ of \ } F.\\
\end{eqnarray*}
The isomorphism of the
proposition for the Tate-Shafarevich groups now follows by a diagram chase.
The corresponding isomorphism for the $n$-Selmer groups  follows from the 
commutative diagram with exact rows
$$\matrix{0&\to& \ _{n}E(F)&\to& {\rm Sel}_{n}(E/F)&\to&
{\textstyle\coprod\!\!\!\coprod}(E/F)_{n}&\to& 0\cr
&&\downarrow\cong &&\downarrow&&\downarrow\cong&&\cr
0&\to& (\ _{n}E(F'))^{{\rm Gal}(F'/F)}&\to& {\rm Sel}_{n}(E/F')^{{\rm Gal}(F'/F)}&\to&
{\textstyle\coprod\!\!\!\coprod}(E/F')_{n}^{{\rm Gal}(F'/F)}&\to& 0\cr
}$$ 
as required. ${\sqcap \!\!\!\!\sqcup}$

\vskip0.2in
\noindent {\it 2.2.6. Remark.}  This section is a generalised form of [1,
Chapter 7, \S7.9].

\vfil\eject
\vskip0.4in
\noindent{\bf 2.3. The Cassels pairing} 

\vskip0.2in\noindent(2.3.1)  Let $E/F$ be an elliptic curve over the global field $F$
of characteristic $p>0$ as in \S1.5. 
The Tate-Shafarevich group $\coprod\!\!\!\coprod (E/F)$ of $E/F$ is equipped with the  anti-symmetric Cassels pairing 
$$
<,>_{\rm Cassels}\ :\ \coprod\!\!\!\!\coprod(E/F)
\times
\coprod\!\!\!\!\coprod(E/F)\to {\mathds Q}/{\mathds Z}$$ 
which is 
non-degenerate if $\coprod\!\!\!\coprod$$(E/F)$ is finite.  The Tate-Shafarevich group $\coprod\!\!\!\coprod$$(E/F)$ is a torsion group of cofinite type (see [10]).

In this section, the Cassels pairing is defined for $E/F$ for the non-$p$ part 
$\coprod\!\!\!\coprod (E/F)_{   ({\rm non}-p)  }$ of 
 $\coprod\!\!\!\coprod (E/F)$, that is to say the subgroup of the Tate-Shafarevich 
group
of order coprime to the characteristic $p$.

\vskip0.2in\noindent(2.3.2) For any place $v\in \Sigma_F$ of $F$, we have the commutative diagram with exact rows obtained from restriction from
$F$ to $ F_v$  
for any integer $m$ where $\partial_m$ is the connecting homomorphism induced
from  the morphism of multiplication by 
$m$ on $E$ 
$$\matrix{
0&\to& E(F)_m&\to &E(F)&\to&E(F)&{\buildrel \partial_m\over \to} & H^1(F,E_m)&\to & H^1(F,E)_m\cr
&&\downarrow &&\downarrow 
&&\downarrow 
&&\downarrow 
&&\downarrow 
\cr
0&\to& E(F_v)_m&\to &E(F_v)&\to&E(F_v)&{\buildrel \partial_m\over \to} & H^1(F_v,E_m)&\to & H^1(F_v,E)_m\cr}$$

\vskip0.2in\noindent (2.3.3)
Let $a,b\in$ $\coprod\!\!\!\coprod (E/F)_{   ({\rm non}-p)  }$.
Let $m\geq 1$ be the order of $a$ and $n\geq 1$ be the order of $b$
where    $m,n$ are coprime to the characteristic 
$p$ of $F$.  Then we have 
\vskip0.2in
$$a\in \coprod\!\!\!\!\coprod(E/F)_m{\rm \ \ and \ \  } b\in \coprod\!\!\!\!\coprod(E/F)_{n}.$$ 

\vskip0.2in
\noindent Select elements 
$$a^{(1)}\in  H^1(F, E_m)  {\rm \ \ and \ \ } b^{(1)}\in H^1(F, E_{n})$$
mapping to $a$ and $b$ respectively in the commutative diagram of (2.3.2).

For any element $h\in H^1(F,E)$ denote by $h_v$ the restriction of $h$ to 
$H^1(F_v,E)$ for any place $v$ of $F$ and similarly for cochains.

\vskip0.2in\noindent(2.3.4) Suppose first that $a$ is divisible by $n$ in $H^1(F,E)$, say $a=na_1$ where $a_1\in H^1(F,E)_{mn}$. We may select local points $y_v\in {}_nE(F_v)$ such that 
$$\partial_n(y_v)=b^{(1)}_v,{\rm \ \ for \ all  \ places\ }v\in \Sigma_F,$$
as $b^{(1)}_v$ maps to zero in $H^1(F_v,E)$. Let $a_{1,v}$ denote the  localization in $H^1(F_v,E)_{mn}$ of $a_1$. Note that since $a\in $ $\coprod\!\!\!\coprod$$(E/F)_m$ we have $a_{1,v}\in H^1(F_v,E)_{n}$ for all $v$.  
For any $v\in \Sigma_F$, denote by
$$[,]_v:  H^1(F_v,E)_n \times {}_n E(F_v)\to {{\mathds Z}\over n{\mathds Z}}.$$
the local pairing as in Theorem 2.1.6.
Then the Cassels pairing is given in terms of the local  pairing 
 by the formula
$$<a,b>_{\rm Cassels}=\sum_{v\in \Sigma_F}[a_{1,v},y_v]_v  \leqno{(2.3.5)}$$
where the sum runs over all places of $F$.

\vskip0.2in\noindent(2.3.6)  We have 
$$a_{1,v}=i_*(c_{1,v})$$
for some $c_{1,v}\in H^1(F_v,E_{n})$ for all $v$
where $i$  is the inclusion of group schemes $E_n\hookrightarrow E$.
We then have  that the Cassels pairing is also given in terms of the cup product pairing $<,>_v$ of Theorem 2.1.4 by 
$$<a,b>_{\rm Cassels}=\sum_v
<c_{1,v},b^{(1)}_v>_v  \leqno{(2.3.7)}$$
as we have for all places $v$
$$[a_{1,v},y_v]_v  =<c_{1,v},\partial_n(y_v)>_v .$$

\vskip0.2in\noindent(2.3.8) In the previous two paragraphs, the global element $a_1\in H^1(F,E)_{mn}$ such that $na_1=a$ may not exist. Nevertheless, in a suitable sense it always exists 
locally and   in general the  Cassels pairing is defined as follows.

Select elements $a^{(1)}$ and $b^{(1)}$ of $H^1(F, E_m)$ and $H^1(F, E_{n})$ mapping to $a$ and $b$ respectively.
For each valuation $v$ of $F$, let $a^{(1)}_v$ be the localisation in $H^1(F_v,E_m)$ 
of $a^{(1)}\in H^1(F,E_m)$. For each valuation $v$ of $F$, $a$ maps to zero in 
$H^1(F_v,E)$ and hence $a^{(1)}_v$ lies in the image of $E(F_v)$ under  
$\partial_m$.
Then we can lift, by the diagram  where id is the identity map,
$$\matrix{ E(F_v)&{\buildrel \partial_m\over \to} & H^1(F_v,E_m) \cr
{\rm id}\uparrow&&\phantom{n}\uparrow n\cr
E(F_v)&{\buildrel \partial_{mn}\over \to} & H^1(F_v,E_{mn})\cr}$$
 $a^{(1)}_v\in H^1(F_v,E_m)$ to an element $a^{(1)}_{v,1}\in H^1(F_v,E_{mn})$ so that $na^{(1)}_{v,1}=a^{(1)}_v$ and $a^{(1)}_{v,1}$ is in the image of $E(F_v)$ under
$\partial_{mn}$.

 Let $\beta$ be a cocycle in ${\rm Cocy}^1(F,E_m)$ representing $a^{(1)}\in H^1(F,E_m)$ and lift it to a cochain $\beta_1\in {\rm Coch}^1(F,E_{mn})$. Select a cocycle $\beta_{v,1}\in {\rm Cocy}^1(F_v,E_{mn})$ representing the element $a^{(1)}_{v,1}\in H^1(F_v, E_{mn})$ and a cocycle $\beta'\in {\rm Cocy}^1(F,E_{n})$ representing $b^{(1)}\in H^1(F,E_n)$. 

The coboundary $d\beta_1$ of $\beta_1$ takes values in $E_{n}$ as $\beta_1$ is a cochain lifting the cocycle $\beta$ with values in $E_m$. The cup product $d\beta_1\cup \beta'$ represents an element of $H^3({\rm Gal}(F^{\rm sep}/F),
{\mathds G}_m)$ where ${\mathds G}_m$ is the multiplicative 
group scheme over $F$ and where this cup product is induced by the 
Weil pairing $E_n\times E_n\to {\mathds G}_m$ (see (2.1.3)).  But this last cohomology group 
$H^3({\rm Gal}(F^{\rm sep}/F),
{\mathds G}_m)$ is zero  (see [9, Chap. 1, 4.18 or 4.21]). Hence we have 
$$d\beta_1\cup \beta'=d\epsilon$$
 for some 2-cochain $\epsilon\in {\rm Coch}^2(F,{\mathds G}_{\bf m})$. 

The cochain
 localized at $v$
$$((\beta_{v,1}-\beta_{1,v})\cup \beta_v')+\epsilon_v$$
is a 2-cocycle in ${\rm Cocy}^2(F_v, {\mathds G}_{ m})$. 
Denote by 
$${\rm inv}_v: {\rm Br}( F_v)=H^2(F_v, {\mathds G}_m) \to {\mathds Q}/{\mathds Z}$$
the canonical isomorphism (the invariant map) given by local class field theory. Define
$$<a,b>_{\rm Cassels}=\sum_{v\in \Sigma_F}
{\rm inv}_v(((\beta_{v,1}-\beta_{1,v})\cup \beta_v')+\epsilon_v)\in {\mathds Q}/{\mathds Z}.$$
This can be checked to be independent of the selections made and defines the pairing on $\coprod\!\!\!\coprod (E/F)_{   ({\rm non}-p)  }$.

\vskip0.2in

\noindent {\it 2.3.9. Remarks.} (1) For this paper, only that part of the construction of the Cassels pairing in 
paragraphs (2.3.2)-(2.3.6) is required. This is  because, under the hypotheses of  the main theorems of this paper
stated in \S4.1, the subgroup
$\coprod\!\!\!\coprod (E/K)$ of $H^1(K,E)$ has the following divisibility property. 

Under the hypotheses and notations of  the main theorems 
of this paper
stated in \S4.1, for all prime numbers $l\in \cal P$ where $l$ is coprime
to Pic$(A)$ and for  any integer $n\geq 0$ there is 
a subgroup $H$ of 
$ H^1(K,E)$ such that $l^nH=\coprod\!\!\!\coprod (E/K)_{l^\infty}$.

 This 
divisibility holds because
 for any sufficiently large positive integer $a$  the group 
$\coprod\!\!\!\coprod (E/K)_{l^\infty}$  is generated by the cohomology
 classes $\delta_{M_0}(c)$ 
where $c$ ranges of the elements of 
$\Lambda(a)$ (Theorem 5.6.2) and that
$l^n\delta_{M_0+n}(c)=\delta_{M_0}(c)$  if $c\in \Lambda(M_0+n)$ 
(lemma 4.2.1(iii)).

It would be interesting to have examples of elliptic curves $E/K$ and prime numbers $l$ for which  $\coprod\!\!\!\coprod (E/K)_{l^\infty}$ is non-zero but does not have this  divisibility 
property as a  subgroup of $H^1(K,E)$.

\noindent (2) The Cassels pairing for abelian varieties over global fields can be 
defined in several ways. For a more geometric construction of the pairing than that 
given above, see [9,  Chap. I, Remark 6.11, p. 98]. For the special case of Jacobians
of curves, see [9, Chap. I, Remark 6.12, p. 100]. For the construction of the 
pairing  including the $p$-torsion part of the 
Tate-Shafarevich group, where $p$ is the characteristic of the base field, see
[9, Chap. II, Theorem 5.6, pp. 247-248].

\vskip0.2in

\vfil\eject

\noindent {\bf Chapter 3. The cohomology classes $\gamma_n(c),\delta_n(c)$}
\vskip0.4in
 \noindent \noindent {\bf 3.1.
 The set $\cal P$ of prime numbers}

\vskip0.2in We define a set $\cal P$ of prime numbers by arithmetic conditions. For prime numbers $l$ of $\cal P$ we shall consider in Chapters 4 and 5
the structure of  the 
$l$-primary component  of the Tate-Shafarevich group 
${\coprod\!\!\!\coprod}(E/F)$ of  elliptic curves $E/F$.
\vskip0.2in
\noindent (3.1.1) For a place $v\in \Sigma_F$ of the global field $F$, let 

\vskip0.1in

$F_v$ denote the completion of $F$ at $v$;

$F_v^{\rm nr}$ denote the maximal
unramified extension of the local field $F_v$ 

\quad \ \ (that is, $F_v^{\rm nr}$ is  the
field of fractions of the strict henselisation of the 

\quad \ \ valuation ring of $F_v$);

 $O_v$ denote the discrete valuation ring of the local field $F_v$;
 
 $\infty\in \Sigma_F$ be a place of $F$;
 
 $K$ be a separable imaginary quadratic extension field of $F$ with respect to the 

\qquad place $\infty$ as in \S1.2; 
 
 $E/F$ be an elliptic curve (as in \S1.5);

 $\cal
E$ denote the N\'eron model  over $O_v$ of the elliptic curve $E\times_FF_v/F_v$;

 ${\cal E}_0$ denote the closed fibre of ${\cal E}/O_v$;

$\pi_0({\cal E}_0)$ be the group of connected components of ${\cal E}_0$ 

\quad \ \ as a
Gal$(F_v^{\rm nr}/F_v)$-module.

\vskip0.2in

\noindent Define similarly $K_w$, $K_w^{\rm nr}$ for a place $w\in \Sigma_K$ of 
the imaginary quadratic extension field $K$ of $F$.

\vskip0.2in\noindent {\bf 3.1.2. Theorem.} ([9, Chap. I, Prop. 3.8, p.57]). {\sl  Write 
$G={\rm Gal}(F_v^{\rm nr}/F_v)$. 
There is an isomorphism
$$H^1(G, E(F_v^{\rm nr}))\cong H^1(
G, \pi_0({\cal E}_0)).$$
In particular, $H^1(G, E(F_v^{\rm nr}))$ is a finite
group for all $v$ and if the elliptic curve $E$ has good reduction  at $v$  then  $H^1(G,
E(F_v^{\rm nr}))=0$.}\qquad ${\sqcap \!\!\!\!\sqcup}$

\vskip0.2in
\noindent {\bf 3.1.3. Definition.} 
Let $\cal P$  be the set of all prime numbers 
such that for all $l\in \cal P$ we have \vskip0.2in \noindent (a) $p$, $2$, and
the prime factors of $\vert B^*\vert/\vert A^*\vert$  are not in $\cal P$;

\noindent (b) $H^i(K(E_{l^n})/K,E_{l^n})=0$ for all integers $n\geq
1$ and for all  $i\geq 0$

\qquad (see Proposition 1.12.1);

\eject
\noindent (c) the natural injection Gal$(F(E_{l^\infty})/F)\to
{\widehat \Gamma}_l$ is an isomorphism 

\qquad (see \S1.11 and Igusa's Theorem 1.11.7);

\noindent (d) $H^1(K^{\rm nr}_z/K_z,E)_{l^\infty}=0$ for all places
$z$ of $K$ (see Theorem 3.1.2 above); 

\noindent (e) $K$ and $F(E_{l^\infty})$ are linearly disjoint over $F$
(see Proposition 1.12.3 or [1, Chapter 7, 7.3.10]);

\noindent (f) $\cal P$ excludes the prime numbers of the finite set 
$\cal E$ of Proposition 1.10.1, that is 

\qquad to say, we have $E(K[A])_{l^m}=0$ 
for all $l\in {\cal P}$, for all $m\geq 1$, where $K[A]$ 

\qquad is defined in \S1.10;

\noindent (g) $\pmatrix{1&0\cr0&-1\cr}\in {\rm Gal}(F(E_{l^\infty})/F)$ (see
Proposition 1.12.3, which is a 

\qquad consequence of Igusa's Theorem 1.11.6);

\vskip0.2in
\noindent {\it 3.1.4. Remarks.} (1)  The first 6 conditions (a),(b),(c),(d),(e),(f) of
this definition hold for all except finitely many prime numbers $l$. Only the
last
condition (g) fails to hold in general for all but finitely many prime numbers.
The set $\cal P$ is infinite and  has positive Dirichlet density.

The set $\cal P$ can therefore be obtained 
from the set $S$ of prime numbers provided by  Proposition 1.12.3 by deleting a
finite number of elements. More precisely,  (see [1, Chapter 7, Remarks 7.3.14(1)] and also [1, Chapter 7, Remarks 7.7.6(1)]) the set $\cal
P$ consists of all but finitely many prime numbers $l$ of the form 
$2^sn+1$ where $s\geq 1$ and $n$ is odd such that $q=\vert k\vert$ is a $2^s$th
power non-residue modulo $l$.

\vskip0.1in

\noindent (2) The set $\cal P$ of prime numbers of Definition 3.1.3 above
coincides with the set of prime numbers written ${\cal P}\setminus {\cal F}$ of 
[1, Lemma 7.14.11]. The only  difference between the Definition 3.1.3 above of the 
set of prime numbers $\cal P$ and the similar definition [1, Definition 7.10.3] is the extra hypothesis  (f) 
of 3.1.3 above which excludes from $\cal P$ the finitely many prime numbers of the exceptional set $\cal E$ of Proposition 1.10.1.

\vskip0.2in\noindent{\bf 3.1.5. Lemma.} {\sl For any integers $0\leq m\leq n  $ 
and for all prime numbers $l\in \cal P$ the inclusion of group schemes $E_{l^m}\to E_{l^n}$ induces an injection of cohomology groups  
$$H^1(K, E_{l^m}){\longrightarrow} H^1(K, E_{l^{n}}).$$
 }
\vskip0.2in
\noindent {\it Proof of Lemma 3.1.5.} This follows from
the long exact sequence induced by the isogeny on the finite group scheme $E_n$ of 
multiplication by $l^{m}$
$$0\to E_{l^m}(K)\to E_{l^{n}}(K)\to E_{l^{n-m}}(K)\to$$
$$
H^1(K, E_{l^m})\to H^1(K, E_{l^{n}})\to  H^1(K, E_{l^{n-m}})\to\ldots$$
together with the non-existence of $K$-rational $l$-torsion on $E$ (by
Proposition 1.10.1 and the definition of $\cal P$, Definition 3.1.3(f)).
${\sqcap \!\!\!\!\sqcup}$ 

\vskip0.2in
\noindent (3.1.6)
We write, where the limits are set-theoretic unions by the previous Lemma 3.1.5,
$$H^1(K, E_{l^\infty})=\lim_{\longrightarrow\atop n}H^1(K, E_{l^n})$$
and
$$
{\rm Sel}_{l^\infty}(E/K)=\lim_{\longrightarrow\atop n} {\rm Sel}_{l^n}(E/K).$$
\vskip0.2in
\noindent (3.1.7)  Let ${\mathds Q}(l)$, for any prime number $l$,  be the additive group of rational numbers with denominators a power of $l$; the quotient group ${\mathds Q}(l)/{\mathds Z}$  is a divisible   abelian group where every element is annihilated by a power of $l$.

We have the exact sequence of abelian groups for all 
prime numbers $l\in \cal P$ where $E(K)_{{\rm tors}}$ denotes the torsion subgroup of $E(K)$ 
$$0\longrightarrow {E(K) \over E(K)_{\rm tors}}\otimes_{\mathds Z} {{\mathds Q}(l)
\over {\mathds Z}}
\longrightarrow {\rm Sel}_{l^\infty}(E/K)\longrightarrow 
{\textstyle \coprod\!\!\!\coprod}(E/K)_{l^\infty} \longrightarrow 0\leqno{(3.1.8)}$$
 This exact sequence (3.1.8) is obtained from Lemma 3.1.5, the exact sequence of 
paragraph (2.2.3), and because $E(K)$ has no $l$-torsion (Definition 3.1.3(f)).

The exact sequence (3.1.8) splits and gives  the isomorphism 
$$ {\rm Sel}_{l^\infty}(E/K)\cong {E(K) \over E(K)_{\rm tors}}\otimes_{\mathds Z} {{\mathds Q}(l)
\over {\mathds Z}}\ \ \oplus \ \ 
{\textstyle \coprod\!\!\!\coprod}(E/K)_{l^\infty}.\leqno{(3.1.9)}$$
This holds as the abelian group  ${E(K) \over E(K)_{\rm tors}}\otimes_{\mathds Z} {{\mathds Q}(l)
\over {\mathds Z}}$, where $E(K)$ is a finitely generated group, is injective in the category of abelian groups.

It follows from the isomorphism (3.1.9) that Sel$_{l^m}(E/K)$ is precisely the subgroup of $ {\rm Sel}_{l^\infty}(E/K)$ annihilated by $l^m$ that is to say we have for all $m\geq 0$ and all $l\in {\cal P}$
$$ {\rm Sel}_{l^m}(E/K)\cong ( {\rm Sel}_{l^\infty}(E/K))_{l^m}.\leqno{(3.1.10)}$$
\vfil\eject

\noindent {\bf 3.2. Frobenius elements and the set  
${\Lambda}(n)$ of divisors}

\vskip0.2in\noindent (3.2.1) Let 

\vskip0.2in $F$ be the function field of the curve $C/k$ as in \S1.2;

$\infty\in \Sigma_F$ be a closed point of $C/k$;

$K$ be a separable imaginary quadratic extension field of $F$ with respect to $\infty$

\qquad as in \S1.2;

$\tau\in {\rm Gal}(K/F)$ be the non-trivial element of the Galois group of $K/F$:

$E/F$ be an elliptic curve with conductor $I$;

$\cal P$ be  the set of prime numbers associated to $E,F,K$ as in \S3.1;

$l\in \cal P$ be a prime number in $\cal P$.

\vskip0.2in
\noindent As in paragraph (1.12.2), for every  prime number $l\not=p$, an
isomorphism is fixed between ${\rm
Gal}(F(E_{l^\infty})/F)$ and a subgroup of ${\rm GL}_2({\mathds Z}_l)$ by fixing a basis of the corresponding Tate module.

\vskip0.2in\noindent (3.2.2) 
For each  integer $n\geq 1$ 
and the chosen $l\in {\cal P}$, let $\tau_\infty\in {\rm Gal}(K(E_{l^n})/F)$ be the
unique element of this Galois group satisfying the two conditions:
\vskip0.1in
\qquad  {\rm (a)} ${\displaystyle{\ \tau_\infty\vert_{F(E_{l^n})}=\pmatrix{1&0\cr0&-1\cr} }}$ that is to say 
the restriction of $\tau_\infty$ to the field 

\qquad \qquad extension $F(E_{l^n})/F$
is ${\displaystyle{\pmatrix{1&0\cr0&-1\cr} }};$
\vskip0.1in
\qquad {\rm (b)} ${\displaystyle{\tau_\infty\vert_{K}=\tau }}$  is the  non-trivial  element of 
Gal$(K/F).$
\vskip0.1in
 
\noindent The elements $\tau$ and $\tau_\infty$ have exact  order 2.

\vskip0.2in\noindent (3.2.3) For any ${\mathds Z}[{\rm
Gal}(K/F)]$-module $M$ on which multiplication by $2$ is an 
isomorphism, we have a decomposition  of $M$ as a sum of eigenspaces under the action of $\tau$, the non-trivial element of Gal$(K/F)$,
$$M\cong M^+\oplus M^-$$
where $M^+$ is the submodule of $M$ on which $\tau $
acts as $1$ and  where $M^-$ is similarly the submodule of $M$ on which $\tau $
acts as $-1$.

\vskip0.2in\noindent {\bf 3.2.4. Definition.} (i) For a prime divisor $z\in \Sigma_F$, unramified in the field extension $K(E_{l^n})/F$, let 
$${\rm Frob}(z)$$ denote the conjugacy
class of Gal$(K(E_{l^n})/F)$ containing the  Frobenius substitutions of the prime
divisors above $z$.

 \noindent (ii) For $l\in \cal P$, let $\Lambda^1(n)$ be the set of
prime divisors $z$ in $\Sigma_F$, of support coprime to $\infty$ and Supp$(I)$ and the
discriminant of $K/F$,  which
satisfy $${\rm Frob}(z)=[\tau_\infty]$$where $[\tau_\infty]$
denotes the conjugacy class in Gal$(K(E_{l^n})/F)$ of $\tau_\infty$. 

\noindent (iii) For $r\geq 0$, let $\Lambda^r(n)$ be the set of
effective divisors $z_1+\ldots+z_r$ on the affine curve Spec $A$ which have $r$ prime components $z_i$ all of which  have multiplicity 1 and
belong to $\Lambda^1(n)$.

We conventionally put for all $n\geq 1$
$$\Lambda^0(n)=\{0\}$$
which is the set consisting of the zero divisor on $C/k$.

\noindent (iv) Put  
$$\Lambda^r=\Lambda^r(1)$$
$$\Lambda(n)=\bigcup_{r\geq 0} \Lambda^r( n)$$
$$\Lambda=\bigcup_{n\geq 1} \Lambda( n).$$
We have
$$\Lambda^r(n)=\left\{\matrix{ z_1+\ldots +z_r\in {\rm Div }_+(A)\ &\vert& \ z_1,\ldots,z_r {\rm \ are \ distinct\ prime \ divisors\ such \  }\cr
&\vert&
 {\rm that \ for \ all \ }i, \ z_i{\rm \ is \ prime \ to \ } \infty, \ {\rm Supp}( I)  
\cr
&\vert& {\rm and \ the \ discriminant\ of \ }K/F{\rm \ and \ }\cr
&\vert&{\rm Frob}
(z_i)=[\tau_\infty]{\rm \  on \ }
{ K}(E_{l^n})\cr}\right\} .$$
The set $\Lambda^r$ has a decreasing filtration
$$\Lambda^r=\Lambda^r(1)\supseteq \Lambda^r(2)\supseteq\Lambda^r(3)\supseteq\ldots.$$
\vskip0.2in 
\noindent {\it 3.2.5. Remarks.} (1)  The prime divisors in $\Lambda(n)$ are infinite in number, by the Chebotarev density theorem, remain
prime in the field extension $K/F$, and their liftings to $K$ split completely
in $K(E_{l^n})/K$. Furthermore, $E$ has good reduction at all prime divisors of
$\Lambda(n)$.

Note that the prime number $l\in \cal P$ is considered to be fixed and the sets
$\Lambda^r(n)$ depend on $l$.

\vskip0.1in\noindent (2) The set $\Lambda(n)$ is defined for any prime number $l$ in $\cal P$ and  
 contains only  effective divisors on Spec $A$ consisting of 
sums of distinct prime divisors whose corresponding Frobenius conjugacy classes in 
${\rm Gal}(K(E_{l^n})/F)$ are all the same and equal to $[\tau_\infty]$. This unique Frobenius conjugacy class is of a special
kind, in particular its elements have order $2$.

\vskip0.2in\noindent \noindent (3.2.6) For any prime number $\lambda$ distinct from the
characteristic of $F$,  let $$\rho:{\rm Gal}(F^{\rm sep}/F)\to {\rm Aut}_{{\mathds
Q}_\lambda}(T_\lambda(E)\otimes_{{\mathds Z}_\lambda}{\mathds Q}_\lambda)$$
denote the galois representation on the $\lambda$-adic Tate module
$T_\lambda(E)$  of the  elliptic curve $E$; that is to say
$$T_\lambda(E)\otimes_{{\mathds
Z}_\lambda}{\mathds Q}_\lambda=H^1_{\rm \acute et}(E\otimes_FF^{\rm sep},{\mathds 
Q}_\lambda)^{*}$$
where $*$ denotes the dual ${\mathds Q}_\lambda$-vector space.

If $z$ is a  prime of $F$ let $I_z$ be an inertia subgroup of 
${\rm Gal}(F^{\rm
sep}/F)$
at $z$; let $$a_z={\rm Tr}(\rho({\rm Frob}(z))\
\vert \ (T_\lambda(E)\otimes_{{\mathds Z}_\lambda}{\mathds Q}_\lambda)^{I_z}).$$
That is to say, $a_z$ is the trace of the Frobenius at $z$ on the part of the
Tate module invariant under $I_z$. Then we have $a_z\in {\mathds Z}$ (see [1,  Chapter 5, Examples
5.3.18(1)]).

 \vskip0.2in
\noindent {\bf  3.2.7. Lemma.} {\sl Suppose that $z\in \Sigma_F$ 
is  a
prime divisor of $F$ and  $l\in {\cal P}$; write $\kappa(z)$ for the
residue field at $z$. . 
}

\noindent (i)  {\sl If $z\in \Lambda^1(n)$  then we have  $a_z\equiv \vert \kappa(z)\vert +1\equiv 0$ {\rm (mod
$l^n$)}.}

\noindent (ii)  {\sl If $z\in \Lambda^1(n)$  and  ${\cal E}_{0,z}$ denotes the 
closed fibre
over $z$ of the N\'eron model of $E/F$ and $y$ is the prime 
divisor of $K$ lying over $z$ then we have  group isomorphisms, for
$\delta=+1$ or $-1$, }
 $${\cal
E}_{0,z}(\kappa(y))_{l^n}^\delta\cong {\cal
E}_{0,z}(\kappa(z))_{l^n}\cong {\mathds Z}/l^n{\mathds Z}.$$

\vskip0.2in
\noindent {\it Proof.} (i)  As $a_z$ is the trace of a Frobenius
above $z$ on the Tate module of $E/F$, we
have by the Grothendieck-Lefschetz trace formula, where $\kappa(z)$ is the residue field of $F$ at $z$
and  ${\cal E}_{0,z}$ is the 
closed fibre
over $z$ of the N\'eron model of $E/F$,
$$a_z=\vert\kappa(z)\vert+1-\vert{\cal E}_{0,z}(\kappa(z))\vert.$$

On the one hand, let $K'$ be the field $K(E_{l^n})$. The characteristic
polynomial of the Frobenius Frob$_{K'/F}(z)$ above $z$ acting on the $l$-adic Tate module of $E/F$  is
equal to 
$$X^2-a_zX+\vert\kappa(z)\vert $$

 On the other hand  as $z\in \Lambda^1(n)$, the prime $z$  is
a place of good reduction of $E/F$ (see remark 3.2.5(1)). Therefore the characteristic
polynomial of the Frobenius Frob$_{K'/F}(z)=[\tau_\infty]$ above $z$ acting on $E_{l^n}$, the 
group scheme of $l^n$-torsion points of $E$, is
equal to by (3.2.2)(a)
$$X^2-1
\ {\rm (mod\ \ } l^n  {\rm )}.$$

Comparing these two quadratic polynomials modulo $l^n$
proves the congruences in part (i) of the lemma. 

\noindent (ii) Let $y $ be the unique place of $K$ lying over the place $z$
of $F$ where $z\in \Lambda^1(n)$. Then $\kappa(y)$ is a quadratic
extension of $\kappa(z)$. Furthermore, as  Frob$_{K'/F}(z)=[\tau_\infty]$
where $\tau_\infty$ has order 2 (see (3.2.2), remark 3.2.5(2)), the prime $y$
splits completely in the extension $K'/K$. The map of reduction
modulo a prime of $K'$ over $z$ 
$$E(K')_{l^n}\to {\cal E}_{0,z}(\kappa(y))_{l^n}$$
is an isomorphism. Hence we have $${\cal E}_{0,z}(\kappa(y))_{l^n}\cong
\bigg({{\mathds Z}\over l^n{\mathds Z}}\bigg)^2.$$ 
As $l$ is an odd prime number and the roots $\pm1$
of the characteristic polynomial $X^2-1$ of $\tau_\infty$ on $E_{l^n}$ are rational
over the prime field ${\mathds Z}/l{\mathds Z}$,
the
action of $\tau_\infty $ on 
${\cal E}_{0,z}(\kappa(y))_{l^n}$ decomposes into a sum over the eigenspaces
of $\tau_\infty$.  Hence we have  for $\delta=\pm1$
 $${\cal E}_{0,z}(\kappa(y))_{l^n}^{\delta}\cong
{\mathds Z}/l^n{\mathds Z}.$$ 
Furthermore,  we have $${\cal
E}_{0,z}(\kappa(y))^{+}={\cal E}_{0,z}(\kappa(z)).$$ The
result follows from this. ${\sqcap \!\!\!\!\sqcup}$

\vskip0.2in\noindent {\it 3.2.8.  Remarks.} (1)  Let $z\in \Sigma_F$ be a prime divisor which is coprime to $\infty$, Supp$(I)$ and the discriminant of $K/F$ and which is inert in the field extension $K/F$. Then it can be shown, in a similar way to the proof of  Lemma 3.2.7, that $z\in \Lambda^1(n)$ if and only if  
$$a_z\equiv \vert\kappa(z)\vert+1\equiv 0 \ ({\rm mod}\ l^n)$$
and the Frobenius Frob$(z)$ does  {\it not} act as a homothety on 
$T_l(E)\otimes_{{\mathds Z}_l} {\mathds Z}/l^n{\mathds Z}$
where $a_z$ is the trace of the Frobenius Frob$(z)$ on the Tate module as in \S3.2.6.

\vskip0.1in
\noindent (2) Fix a prime $l\in \cal P$. For  two integers $a,b\in 
\mathds Z$ denote by $v_l(a,b)$ the valuation at $l$ of the greatest common 
divisor of $a$ and $b$. For a prime divisor $z$ of Spec $A$
define $\alpha(z)$ by the equation
$$\alpha(z)=v_l(\vert \kappa(z)\vert +1, a_z)$$
where as in \S3.2.6
$$a_z=\vert \kappa(z)\vert +1-\vert {\cal E}_0(\kappa(z))\vert$$
and where ${\cal E}_0$ is the closed fibre of the N\'eron model of $E$ at $z$.

The prime divisors $z$ of Spec $A$ in the set $\Lambda^1$ 
have the property that they are coprime to $I$, remain prime in $K/F$, and satisfy
$$\alpha(z)=v_l(\vert \kappa(z)\vert +1, a_z)\geq 1.$$

\vskip0.1in
\noindent (3) If $r>0$ and $c\in \Lambda^r$ we write
$$\alpha(c)=\min_{z\in {\rm Supp}( c)} \alpha(z), \ \ \alpha(0)=+\infty.$$
We have the following equation
$$\Lambda^r(n)=\{c\in \Lambda^r\ \vert \ \alpha(c)\geq n\}$$
where $n$ is any integer $\geq 1$; $\alpha(c)$ is the greatest integer $n$ such that
$c\in \Lambda^r(n)$.

The set $\Lambda^r$ is then equipped with a decreasing filtration
$$\Lambda^r=\Lambda^r(1)\supseteq \Lambda^r(2)\supseteq \Lambda^r(3)\supseteq \ldots$$

\vskip0.1in
\noindent (4) In the notation of the monograph [1, Chapter 7, \S7.11], we have $\Lambda(n)={\cal D}_{l^n}$

\vfil
\vskip0.4in\eject 
\noindent {\bf 3.3. A refined Hasse principle for  finite group schemes}

\vskip0.2in The finite group schemes in question are $E_{l^n}/F$ and the Hasse principle 
concerns the localization at places of $F$ of finite groups of their  principal homogenous spaces. 
\vskip0.2in\noindent(3.3.1) Suppose  that

\vskip0.2in

 $E$ is an elliptic curve defined over $F$;

 $K$ is  an imaginary quadratic field over $F$ with respect to $\infty$:

 $l\in {\cal P}$ is a, necessarily odd,  prime number in the set $\cal P$ (see \S3.1);

$n \geq 1$ is an integer;

 $L_n$ is the field $K(E_{l^n})$ which is galois over $F$.

 \vskip0.2in
\noindent {\bf 3.3.2 Lemma.} {\sl  The restriction
map from $K$ to $L_n$ induces an isomorphism $$H^1(K,E_{l^n}(L_n))\ 
{\buildrel \cong \over \longrightarrow }\ {\rm Hom}_{{\rm Gal}(L_n/K)}
({\rm Gal}(L_n^{\rm ab}/L_n), E_{l^n}(L_n))$$
 where $L_n^{\rm ab}$ is the maximal separable abelian extension of $L_n$.}
\vskip0.2in
\noindent {\it Proof.} The Hochschild-Serre spectral sequence, where
as in (3.3.1)  $L_n=K(E_{l^n})$,
$$H^p(L_n/K,H^q(L_n,E_{l^n}(L_n)))\Rightarrow H^{p+q}(K,E_{l^n}(L_n))$$
gives rise to a short exact sequence of low degree terms 
$$0\to H^1(L_n/K,E_{l^n}(L_n))\to H^1(K,E_{l^n}(L_n))\to
H^1(L_n,E_{l^n}(L_n))^{{\rm Gal}(L_n/K)}\to$$
$$\to H^2(L_n/K,E_{l^n}(L_n))\to H^2(K,E_{l^n}(L_n)).$$
The
two cohomology groups $H^1(L_n/K,E_{l^n}(L_n))$ and
$H^2(L_n/K,E_{l^n}(L_n))$  are zero by the definition of $\cal P$ (Definition 3.1.3(b)); hence this exact
sequence shows that restriction map of the lemma is an isomorphism.
$$H^1(K,E_{l^n}(L_n))\to H^1(L_n,
E_{l^n}(L_n))^{{\rm
Gal}(L_n/K)}$$
is an isomorphism of {\rm Gal}$(L_n/F)$-modules. The isomorphism of the lemma
follows immediately. ${\sqcap \!\!\!\!\sqcup}$

\vskip0.2in\noindent {\bf 3.3.3. Proposition.}
{\sl  If $S$ is a finite subgroup of $H^1(K,E_{l^n})$ then there is a finite abelian extension $L_{S,n}$ of $L_n$ and an isomorphism of Gal$(L_n/{K})$-modules
$${\rm Gal}(L_{S,n}/L_n)
\cong {\rm Hom}(S,E_{l^n}(L_n)) \leqno{(3.3.4)}$$
$$\!\!\!\!\!\!\!\sigma\mapsto \phi_\sigma$$
and an isomorphism of abelian groups (if $S$ is a ${\mathds Z}[{\rm Gal}(K/F)]$-module then an isomorphism of ${\mathds Z}[{\rm Gal}(K/F)]$-modules)
$$S\cong {\rm Hom}_{{\rm Gal}(L_n/K)}({\rm Gal}(L_{S,n}/L_n), E_{l^n}(L_n)).$$  }

\vskip0.2in For the proof of Proposition 3.3.3, see [1, Chapter 7, Corollary 7.18.10]. 
Note that the set of prime numbers ${\cal P}\setminus F$ of [1, Chapter 7, (7.18.1) 
and corollary 7.18.10] coincides with the set $\cal P$ of prime numbers defined in 
\S3.1
above. ${\sqcap \!\!\!\!\sqcup}$

\vskip0.2in\noindent (3.3.5) For a prime divisor $v$ of $K$,
 denote by  $${\rm res}_{v}: H^1(K,E_{l^n})\longrightarrow  H^1(K_{v},E_{l^n})$$  the restriction homomorphism  at $v$ where 
$K_v$ is the completion of $K$ at $v$.

 Let $S$ be a finite subgroup of $H^1(K, E_{l^n})$. For $s\in S$ and $v$ a place of $K$,  we have in the notation of Proposition 3.3.3 applied to the subgroup $S$
$${\rm res}_v(s)=0\Leftrightarrow 
\phi_\sigma(s)=0{\rm \ for \ all
\ } \sigma\in D_{v'}$$
where $v'$ is a prime divisor of $L_n$ above $v$ and 
$D_{v'}$ is the decomposition group of a prime divisor of $L_{S,n}$ above $v'$ and where $\phi_\sigma\in {\rm Hom}(S,E_{l^n}(L_n))$ as in (3.3.4).

\vskip0.2in\noindent(3.3.6) Let $\tau_\infty\in  $ Gal$(L_n/{F})$ be the 
element of order $2$ of paragraph (3.2.2). The element  $\tau_\infty$ acts as $-1$ on the $l^n$th   roots of unity and preserves the Weil pairing on $E_{l^n}$ (referring to the exact sequence (1.11.3), det$(\tau_\infty)=-1$ acts on the 
$l^n$th roots of unity). 

We obtain (see paragraph (3.2.3)) as $l$ is odd the decompositions into eigenspaces under the action of 
$\tau_\infty$
$$H^1(K, E_{l^n})\cong H^1(K, E_{l^n})^+\oplus H^1(K, E_{l^n})^-$$
$$E_{l^n}(L_n) \cong E_{l^n}(L_n)^+\oplus  E_{l^n}(L_n)^-.$$
We obtain that the $+1$ eigenspace
$${\rm Hom}(H^1(K, E_{l^n}), E_{l^n}(L_n))^{+}$$
is isomorphic to the profinite group of ${\mathds Z}[\tau_\infty]$-homomorphims from 
$H^1(K, E_{l^n})$ to $E_{l^n}(L_n)$ and hence we obtain an isomorphism between
this $\tau_\infty$-invariant subgroup
and the profinite Pontrjagin dual of the discrete torsion abelian group $H^1(K,E_{l^n})$
namely we have an isomorphism, where a basis of $T_l(E)$ is fixed
as in (1.12.2),
$${\rm Hom}(H^1(K,E_{l^n}), E_{l^n}(L_n))^{+}\ \cong\  H^1(K,E_{l^n})^{*{}}$$
where the Pontrjagin dual is given by $$H^1(K,E_{l^n})^{*{}}={\rm Hom}(H^1(K,E_{l^n}), {\mathds Q}/{\mathds Z})
.$$
For a finite subgroup $S$ of $H^1(K, E_{l^n})$ there is similarly an isomorphism
$${\widehat S}\cong  {\rm Hom}(S, E_{l^n}(L_n))^+.$$
\vskip0.2in
\noindent {\bf 3.3.7. Proposition.} {\sl  Let $S$ be 
a finite subgroup of $H^1(K,E_{l^n})$ and let $\chi\in {\widehat S}$. For
any integer $t\geq n$, there 
is a set of positive Dirichlet density of prime divisors  $z\in \Lambda^1(t)$ of $F$ 
unramified in 
$L_t=K(E_{l^t})$ such that
$$\chi=\phi_{{\rm Frob}(z^\times)}$$ 
 for some prime divisor $z^\times$ of 
$L_t$ lying above $z$, where $\phi_{{\rm Frob}(z^\times)}
\in {\rm Hom}(S,E_{l^t}(L_t))^+$ is as in (3.3.4).
}

\vskip0.2in
\noindent {\it Proof.} By Lemma 3.1.5, there is an injection 
of cohomology groups for all $t\geq n$
$$H^1(K, E_{l^n})\to H^1(K,E_{l^t})$$
obtained from the inclusion of finite group schemes $E_{l^n}\subseteq E_{l^t}$. 
The finite subgroup $S$ is then a subgroup of $H^1(K,E_{l^t})$ for all 
$t\geq n$.  
Applying Proposition 3.3.3 to
 $S$ as a subgroup of $H^1(K,E_{l^t})$ we obtain the abelian extension $L_{S,t}/L_t$ where we write $L_t=K(E_{l^t})$.

By 
Proposition 3.3.3 there  is an element $\sigma\in $ Gal$(L_{S,t}/L_t)$ such that 
$$\chi=\phi_\sigma.$$
We have 
$${\rm Gal}(L_{S,t}/L_t)^+\cong {\widehat S}$$
from the isomorphism of (3.3.4) and the isomorphism of paragraph (3.3.6)
$${\widehat S}\cong  {\rm Hom}(S, E_{l^t}(L_t))^+.$$
In particular we have 
   $\phi_\sigma^{\tau_\infty}=\phi_\sigma$. As the order of Gal$(L_{S,t}/L_t)$ is odd, we 
then have $\sigma=\rho^{\tau_\infty}.\rho$ for some $\rho\in $
 Gal$(L_{S,t}/L_t)$. By the Chebotarev density theorem there 
is a set of positive Dirichlet density  of prime divisors $z\in \Sigma_F$ of $F$ 
such that Frob$(z)$ in Gal$(L_{S,t}/{F})$ contains $\tau_\infty\rho$ and where
$z$ unramified in  $L_{S,t}/{F}$. Note the 
the finitely many  prime divisors ramified in the field extension $L_{S,t}/{F}$ 
depend only on $S$ and not on $t$.

Since the 
restriction of $\tau_\infty\rho$ to $L_t$ is $\tau_\infty$, 
we have $z\in \Lambda^1(t)$. We then have that $z$ has residue class 
extension degree  $2$ in $L_t/{F}$ for any prime divisor above $z$ in $L_t$. Hence for 
any $z^\times$, a prime divisor of $L_t$ lying above $z$, 
we have
$${\rm Frob}(z^\times)=(\tau_\infty\rho)^2=\rho^{\tau_\infty}.\rho=\sigma.\ \ \ {\sqcap \!\!\!\!\sqcup}$$
\vskip0.2in

\noindent {\bf 3.3.8. Proposition.} {\sl Let $S$ be a finite subgroup of  $H^1(K,E_{l^n})$ and let $\chi\in {\widehat S}$. Then for any integer $t\geq 1$ there 
is a set of positive Dirichlet density  of prime divisors $z\in \Lambda^1(t)$  such that
  for the prime divisor $y$ of $K$ lying over $z$,  there is  a commutative diagram of group homomorphisms
$$\matrix{ S&{\buildrel {\rm res}_{y}\over \longrightarrow}& 
{\rm res}_{y}(S)\cr
 &\chi \searrow &\cong \downarrow \psi\cr
&&\chi(S)\cr}$$
where $\psi$ is an isomorphism of finite cyclic groups.}

\vskip0.2in
\noindent {\it Proof.}  As $\Lambda^1(m)$, $m\geq 1$, is a decreasing filtration of 
$\Lambda^1(1)$, we may assume that $t\geq n$. 
By Lemma 3.1.5, there is an injection 
of cohomology groups for all $t\geq n$
$$H^1(K, E_{l^n})\to H^1(K,E_{l^t})$$
obtained from the inclusion of finite group schemes $E_{l^n}\subseteq E_{l^t}$. 
The finite subgroup $S$ is then a subgroup of $H^1(K,E_{l^t})$ for all 
$t\geq n$.  

Select the prime divisor $z\in \Sigma_F$ as in Proposition 3.3.7 applied to $S,\chi$ and where the divisors $z$ are different from the finitely many prime divisors of $F$ where the cohomology classes of the finite group $S$ ramify, that is to say the 
finitely many 
prime divisors where the field extension $L_{S,t}/{F}$  of (3.3.4) is ramified. For such a $z$, 
let $y$ be the prime divisor of $K$ above $z$.

 The decomposition group of $y$ in Gal$(L_{S,t}/K)$ is generated by the Frobenius element Frob$({y})$ as the field  extension $L_{S,t}/K$ is unramified at $y$. By paragraph (3.3.5) we then have 
$${\rm ker}(S\ {\buildrel {\rm res}_{y} \over \longrightarrow}\
 H^1(K_{y},E_{l^t}))={\rm ker}(\chi). $$
That is to say, the kernel in $S$ of the restriction homomorphism
at $y$ of 
the classes of $S$ is equal to ker$(\chi)$.
Hence the two finite cyclic groups $ {\rm res}_{y}(S)$, $ \chi(S)$ 
are isomorphic and we obtain the isomorphism $\psi$ of the  commutative diagram
of the corollary. \ \ ${\sqcap \!\!\!\!\sqcup}$

\vskip0.2in
\noindent {\bf 3.3.9. Corollary.} {\sl We have  $$\matrix{\Lambda^r(n)
\phantom{=\{0\}}&{\sl  is 
\ infinite \ for \ all \ integers   \ \ }r, n\geq 1;\cr
 \Lambda^0(n)=\{0\} &{\sl \!\!\!\!\!\!\!\!\!\!\!\!\!\!\!\!\!\!\!\!\!\!\!\!\!\!\!\!\!\!\!\!\!\!\!\!\!\!\!\!\!\!\!\!\!\!\!\!\!\!\!\!\!\!\!\!\!\!\!\! for \ all \ }n\geq 1.\cr}$$}  
\vskip0.2in
\noindent {\it Proof.}  We have from Proposition 3.3.8 that $\Lambda^1(n)$ has infinitely many elements for all integers $n\geq 1$.  As $\Lambda^r(n) $ consists of all sums of $r$ distinct prime divisors of $\Lambda^1(n)$ (see Definition 3.2.4) and as $\Lambda^0(n)=\{0\}$,  the corollary then follows.  ${\sqcap \!\!\!\!\sqcup}$

\vskip0.2in
\noindent {\it  3.3.10. Remark.}  Proposition 3.3.8 implies that, under the 
hypotheses of the proposition, the 
kernel of the natural homomorphism
$$\psi:H^1(K,E_{l^n})\to \prod_{z\in \Lambda^1(n)} H^1(K_{y(z)}, E_{l^n})$$
is zero, where $y(z)$ is the place of $K$ over $z\in \Lambda^1(n)$. To show this, it is sufficient to apply the proposition to a faithful character $\chi$ of the finite cyclic group $S$
generated by any element of $H^1(K,E_{l^n})$. 

The {\it Hasse principle} for the 
group scheme $E_{l^n}$ and the set of prime divisors $\Lambda^1(n)$ 
is precisely that the kernel of the homomorphism $\psi$ is zero. It says that a prinicpal homogeneous space of $E_{l^n}$ which is locally trivial at all places of
$F$ in $\Lambda^1(n)$ must be globally trival. 

In this way, the proposition is a
refined form of the Hasse principle
 for the finite group scheme $E_{l^n}$.
 It would be interesting to know to what other group schemes the localization
properties of finite subgroups of cohomology groups 
 given in Proposition 3.3.8 also hold.
 
 [For more details on the Hasse principle, 
see [9, pp.142-150].]
\vskip0.4in

\count2=1
\noindent {\bf 3.4. Drinfeld-Heegner points and the cohomology classes $\gamma_n(c), \delta_n(c)$}

\vskip0.2in
\noindent (3.4.\the\count2)  \advance\count2 by 1 Let 
\vskip0.2in
$E/F$ be an elliptic curve equipped with an origin, that
is to say $E/F$ is a 

\qquad 1-dimensional abelian variety; 

 $I$ be the ideal of $A$
which is the conductor of $E/F$ without the component 

\qquad at $\infty$;

 $\epsilon=\pm1$ be the sign in the functional equation of the $L$-function 
\label{Lfunction2}
of the  elliptic 

\qquad curve $E/F$.
 \vskip0.2in

\noindent Assume that (see also [1,  Chapter 4, \S\S4.3, 4.7, or Chapter 7, \S7.6]):

\vskip0.2in

($\alpha $)  \qquad  $E/F$    has  split Tate  multiplicative reduction at 
$\infty$   (see  \S1.7

\qquad \qquad or [1, Chapter 4, \S4.7]);

  ($\beta $) \qquad   $K$
  is an (separable) imaginary quadratic field  extension of $F$, 

\qquad \qquad   with
respect to $\infty$, 
    such  that  
 all primes dividing  the  conductor $I$ 

\qquad \qquad  split completely in  $K$.

\vskip0.2in

 \noindent The hypothesis ($\alpha$) implies that $E/F$ is covered by the 
 Drinfeld modular curve $X_0^{\rm Drin}(I)$ (see (3.4.4) below). The hypothesis
 ($\beta$) ensures the existence of Drinfeld-Heegner points on $X_0^{\rm Drin}(I)$. The
 two hypotheses together ensure that there are Drinfeld-Heegner points on the elliptic curve  $E/F$.  Note that from (3.4.7) to the end of this section \S3.4, it
 is assumed that $K\not=F\otimes_{{\mathds F}_q}{\mathds F}_{q^2}$ where
 ${\mathds F}_q$ is the exact field of constants of $F$.

In this section we detail this construction of Drinfeld-Heegner points as well as
 defining the cohomology classes  $\gamma_n(c), \delta_n(c)$.

\vskip0.2in
\noindent (3.4.\the\count2) \advance\count2 by 1  Let $\lambda$ be any prime number of $\mathds Z$ distinct from the characteristic
of $F$. 
 Let $\rho$ be the
2-dimensional $\lambda$-adic representation of Gal$(F^{\rm sep}/F)$ corresponding to $E$ where $F^{\rm sep
}$ is the separable closure of $F$; that is to say $\rho $ is the continuous homomorphism
 $$\rho:{\rm Gal}(F^{\rm sep}/F)\to {\rm End}_{{\mathds Q}_\lambda}(H^1_{\rm \acute
et}(E\otimes_FF^{\rm sep},{\mathds Q}_\lambda)).$$ For each place $v$ of
$F$ put$$a_v={\rm Tr}(\rho({\rm Frob}(v))\vert \ H^1_{\rm \acute
et}(E\otimes_FF^{\rm sep},{\mathds Q}_\lambda)^{T_v})$$ 
where  $T_v$ is the inertia
subgroup of Gal$(F^{\rm sep}/F)$ over $v$. 
The representation $\rho$ satisfies
 $$a_v\in {\mathds Z} {\rm \ for \ all \ } v\in \Sigma_F$$
(see [1, Chapter 5, example 5.3.18]).

\vskip0.2in\noindent(3.4.\the\count2) \advance\count2 by 1 Let 

\vskip0.2in
  $B$ be the integral closure of $A$ in $K$; 
  
  $O_c$ be the order of $K$, with respect to $A$, and of conductor $c$
for any divisor 

\qquad $c\in$ Div$_+(A)$,  (see [1, Chapter 2, \S2.2]);

 $\tau$ be the non-trivial
element of Gal$(K/F)$;

$IB=I_1I_2$ be a factorisation of ideals of $B$ where $I_1$, $I_2$ are ideals of $B$ such 

\qquad that $I_2^\tau=I_1$ where  $I$ is the   ideal of $A$ which is the conductor of $E/F$ 

\qquad without the 
 place at $\infty$, as 
 in (3.4.1); such a factorisation exists

\qquad  because of the hypothesis that 
the prime ideal components of $I$ split 

\qquad completely in $K/F$ (by (3.4.1)($\beta$)).

\vskip0.2in

As in \S1.4, $K[c]$ denotes the ring class field of $K$ of conductor $c\in $ Div$_+(A)$. In particular $K[0]$ is the Hilbert class field of $K$, that is to say $K[0]$ is the maximal unramified abelian extension of $K$ which is 
split completely at $\infty$.

\vskip0.2in

\noindent (3.4.\the\count2) \advance\count2 by 1
Let ${\cal H}(\rho)$  be the Heegner module of $\rho$ and $K/F$  with
exceptional set\label{hmod7} of primes those dividing $I$ and the place $\infty$ with
coefficients in $\mathds Z$ (see [1, Chapter 5, \S5.3]), where $I$ is the conductor of $E$, without the component at $\infty$.

As in \S1.7 (see also  [1, Chapter 4, \S4.7] and [1, Appendix B]),  there is a finite surjective morphism of curves \label{modsch10} 
over $F$ 
under the hypothesis ($\alpha$)
$$\pi:X_0^{\rm Drin}(I)\to E.$$We may translate $\pi$ in the group scheme $E$ so that 
$\pi^{-1}(0)$ consists of at least one cusp of $X_0^{\rm Drin}(I)$ (as in [1, 
Chapter 4,
(4.8.1)]); this rigidifies the map $\pi$. The cusps of the modular curve $X_0^{\rm 
Drin}(I)$ generate a torsion subgroup of the 
jacobian of this curve [1, Chapter 2, Theorem 2.4.9].  

Let $a$ be a divisor class in the Picard group Pic$(O_c)$ of the order $O_c$. Assume that $c$ and $I$ are coprime. Then there is a Drinfeld-Heegner point
$$(a,I_1,c)\in X_0^{\rm Drin}(I)(K[c])$$
which is a non-cuspidal point of $X_0^{\rm Drin}(I)$ and is rational over the
ring class field $K[c]$.

This point $(a,I_1,c)$ is constructed as follows. Fix an embedding 
$K\to {\widehat{\overline{F}}}_\infty$ where ${\widehat{\overline{F}}}_\infty$
is the completion of the algebraic closure of $F_\infty$ which is the completion of 
$F$ at $\infty$. Let $L$ be a projective $O_c$-module of rank 1 in the class $a$ and 
contained as a lattice in ${\widehat{\overline{F}}}_\infty$. Then 
$I_1(O_c)=I_1\cap O_c$ is an invertible ideal of $O_c$ and 
$L'=I_1(O_c)^{-1}L$ is a projective $O_c$ module of rank 1 contained as a 
lattice in ${\widehat{\overline{F}}}_\infty$. Let $D$ and $D'$ be the rank 2
Drinfeld modules for $A$ over the field ${\widehat{\overline{F}}}_\infty$ corresponding
respectively to the lattices $L$ and $L'$. Then $D$ and $D'$ have general
characteristic and complex multiplication by $O_c$. The inclusion of 
$O_c$-modules $L\subset L'$ corresponds to an $I$-cyclic isogeny
$f:D\to D'$, as its kernel is isomorphic as an $A$-module to 
$O_c/I_1(O_c)\cong A/I$. The pair $(D, {\rm ker}(f))$ defines the point
$(a,I_1,c)$ on  $X_0^{\rm Drin}(I)( {\widehat{\overline{F}}}_\infty)$. That this 
point $(a,I_1,c)$  is defined over $K[c]$ results from the main theorem of complex multiplication  (see [1, Chapter 4, \S4.3] for more details).

  \vskip0.2in
\noindent (3.4.\the\count2) \advance\count2 by 1 The image
$$\pi(a,I_1,c)\in E(K[c])$$
is a Drinfeld-Heegner point of $E$ rational over the ring class field $K[c]$ and is written in the notation of [1, Chap. 4, \S4.8]
$$(a,I_1,c,\pi)\in E(K[c]).$$

By  [1, Chap. 5, example 5.3.18] there is a homomorphism of  discrete Gal$(K^{\rm
sep}/K)$-modules
\begin{eqnarray*}
{\cal H}(\pi):\quad {\cal H}(\rho)^{(0)}&\to E(F^{\rm sep})\\
<a,c>&\mapsto (a,I_1,c,\pi)\\
\end{eqnarray*}
for all $c\in {\rm Div}_+(A)$ coprime to $I$,  $a\in {\rm Pic}(O_c)$. The image of this homomorphism ${\cal H}(\pi)$ consists of
the $\mathds Z$-linear combinations of Drinfeld-Heegner points of $E$ rational over
all the ring class fields $K[c]$ for all $c$.

Let $<0,0>$ be the element of the Heegner module ${\cal
H}(\rho)^{(0)}$ given by the principal class of Pic$(B)$, where $B$ is the integral
closure of $A$ in $K$. Let $$(0,I_1,0,\pi)={\cal H}(\pi)(<0,0>)$$ be the
corresponding Drinfeld-Heegner point of $E(K[0])$ (see [1, Chapter 4, (4.8.2)]). Let  $$P_0={\rm
Tr}_{K[0]/K}(0,I_1,0,\pi)\in E(K).\leqno{(3.4.\the\count2)}$$\advance\count2 by 1 
That is to say, $P_0$ is the trace from $K[0]$ to $K$ of the 
point $(0,I_1,0,\pi)$; the point $P_0$ belongs to $E(K)$ and the point
$(0,I_1,0,\pi)$ belongs to $E(K[0])$.

\vskip0.2in\noindent(3.4.\the\count2) \advance\count2 by 1 
We now impose for the rest of this section 
the hypothesis that $K\not=F\otimes_{{\mathds F}_q}{\mathds F}_{q^2}$ where
 ${\mathds F}_q$ is the exact field of constants of $F$, that is to say 
 $K$ is not obtained from $F$ by ground field extension.

Let $c\in \Lambda(1)$ that is to say $c$ is a sum of distinct prime
divisors of $\Lambda^1(1)$ with multiplicity $1$. Let
$$y_c=\pi(0,I_1,c)\in E(K[c]).$$
so that $y_c=(0,I_1,c,\pi)$.
The field $K[0]$ is the Hilbert class field of $K$. Let
$$G_c={\rm Gal}(K[c]/K[0]).$$
As $ K\not=F\otimes_{{\mathds F}_q}{\mathds F}_{q^2}$ by hypothesis, we
have that $B^*=A^*$, where $A^*,B^*$ are  the unit groups of $A,B$, and hence (by [1, 
Chap. 2, (2.3.8), p.19])
 there is a group isomorphism
$$G_c\cong \prod_z G_{z}$$
where the product runs over the prime divisors $z$ in the support of $c$ and 
$$G_{z}={\rm Gal}(K[c]/K[c-z])$$
is a cyclic group of order $\vert\kappa(z)\vert+1$ 
as $z$ is inert in $K/F$ (by [1, Chapter 2, (2.3.12), p.20]). Fix a generator $\sigma_z$ of the cyclic group $G_z$ for all prime divisors $z\in \Lambda^1(1)$.

\vskip0.2in\noindent(3.4.\the\count2) \advance\count2 by 1 Write
for any prime divisor $z\in \Lambda^1(1)$ of $F$, where $D_z\in {\mathds Z}[G_z]$,
$$D_z=-\sum_{i=1}^{\vert\kappa(z)\vert} i.\sigma_z^i.$$
Here $D_z$ is the Kolyvagin element of the map
$h_z: G_z\to {\mathds Z}$ where $\sigma_z^{-i}\mapsto -i$ for all 
$0\leq i\leq \vert \kappa(z)\vert$.
Note the minus sign here which agrees with the Kolyvagin elements of 
[1, Chap. 5, pp. 175-178]. 

\vskip0.2in

\noindent (3.4.\the\count2) \advance\count2 by 1 For any divisor $c\in \Lambda(1)$, put
$$D_c=\prod_{z\in \ {\rm Supp}( c)} D_z$$
where $D_c\in {\mathds Z}[G_c]$. Then $D_c$ is the Kolyvagin element of the map
$h:  G_c\to {\mathds Z}$ given by  $h=\prod_{z\in {\rm Supp}(z)} h_z$.
Let
$${\cal G}_c={\rm Gal}(K[c]/K)$$ where there is an exact sequence of finite 
abelian groups
$$0\longrightarrow G_c\longrightarrow {\cal G}_c\longrightarrow {\rm Gal}(K[0]/K)\longrightarrow 0.$$
 Let $\cal S$ be a set of coset representatives for $G_c$ in ${\cal G}_c$. Define the point $P_c\in E(K[c])$ by, where $y_c=\pi(0,I_1,c)$ is the element of
$E(K[c])$ in (3.4.7), 
$$P_c=\sum_{s\in {\cal S}} sD_cy_c.$$

\vskip0.2in

\noindent (3.4.\the\count2) \advance\count2 by 1 Suppose now that $c\in \Lambda(n)$. We  write $P_c$ (mod $l^n$) for the image of 
$P_c$ in the quotient group ${}_{l^n}E(K[c])$.
Then we have  that $P_c$ (mod $l^n$) 
belongs to
$$P_c \ ({\rm mod\ } l^n)\  \in \Bigg({}_{l^n}E(K[c])\Bigg)^{{\cal G}_c}.$$
This inclusion follows immediately from the formula in ${\mathds Z}[G_c]$
$$(\sigma_z-1)D_z=-\vert G(c/c-z)\vert +\sum_{g\in G(c/c-z)} g$$
for all $z\in $ Supp$(c)$ and that $\vert G(c/c-z)\vert=\vert\kappa(z)\vert+1$
is divisible by $l^n$ for all $z\in $ Supp$(c)$  (for a detailed proof of this inclusion 
$P_c \ ({\rm mod\ } l^n)\  \in ({}_{l^n}E(K[c]))^{{\cal G}_c}$ see [1, Chapter 7, Lemma 7.14.9
and Lemma 7.14.11]).
Furthermore, we have  where $y_0$ is defined in (3.4.7) and where this notation $P_0$ agrees with that of (3.4.6)
$$P_0={\rm Tr}_{K[0]/K}(y_0)\in E(K).$$

[The point $P_c$ (mod $l^n$) in $E(K[c])/ l^nE(K[c])$ coincides with the point denoted $P_c$ in the book  [1, Chapter 7, 
Notation 7.14.10(iii) and Lemma 7.14.11].]

\vskip0.2in\noindent(3.4.\the\count2) \advance\count2 by 1 Let $S$ be the ring ${\mathds Z}/m{\mathds Z}$ where $m$ is any non-zero integer.
The morphism of multiplication by $m$ on the elliptic curve $E$ then provides    
for any divisor $c$ of Div$_+(A)$ prime to $I$
the following commutative diagram with exact rows and an exact right-hand column (see [1, Chap. 7, \S 7.14.5]):\vskip0.2in
$$\matrix{&&&&&0&&\cr
&&&&&\downarrow&&\cr
&&&&&\hidewidth H^1(K[c]/K,E(K[c]))_{m}\hidewidth&&\cr
&&&&{_j}&{\rm inf}\downarrow&&\cr
0\to&{_{m}}\!E(K)&\to&H^1(K,E_{m})&\to&H^1(K,E)_{m}&\to&0\cr
&\downarrow&&{\rm res}\downarrow{\rm quasi-isom.}&&{\rm res}\downarrow&&\cr
0\to&(_{m}\!E(K[c]))^{{\cal
G}_c}&{\buildrel\partial\over\to}&H^1(K[c],E_{m})^{{\cal
G}_c}&\to&H^1(K[c],E)_{m}^{{\cal G}_c}&&\cr
&f\uparrow&&&&&&\cr
&({\cal H}_{c,S}^{(0)})^{{\cal
G}_c}&&&&&&\cr}\leqno{(3.4.\the\count2) }$$ \advance\count2 by 1
Here $({\cal H}_{c,S}^{(0)})$ is the $c$-component Heegner module with coefficients
in the ring $S$ (see [1, Chap. 7, (7.11.3), and Chap. 5, \S5.3]).

Let $\cal E$ be the finite exceptional set of prime numbers 
of Propositions 1.10.1 and 1.10.2. 
The middle restriction homomorphism here in (3.4.12) is a quasi-isomorphism of quasi-groups in $[{\bf
N}^{(p)}]_{\mathds Z}$ where the finite exceptional set of prime numbers is $\cal E$ and is independent of $c$ (Propositions 1.10.1 and 1.10.2, or alternatively [1, Chap. 7, Proposition 7.14.2]). In particular this middle restriction homomorphism is an isomorphism for all integers $m$ which are powers of prime numbers of $\cal P$ 
(Definition 3.1.3) as $\cal P$ excludes the finitely many prime numbers of 
$\cal E$.

The map $f$ is obtained by 
sending a 
generator $<a,c>$ of ${\cal H}^{(0)}_{c,S}$ to its image  $(a,I_1,c,\pi)\in E(K[c])$
as in (3.4.4) and 
(3.4.5).

         This diagram (3.4.12) then provides the { Heegner homomorphism}, for all
integers $m\in \mathds N$ prime
to $\cal E$,$$({\cal H}_{c,S}^{(0)})^{{\cal
G}_c}\to H^1(K,E)$$whose image belongs to $ H^1(K[c]/K,E(K[c]))_{m}$. 

\vskip0.2in\noindent(3.4.\the\count2) \advance\count2 by 1 Take $m=l^n$ where $l\in \cal P$ so that the middle quasi-isomorphism in diagram (3.4.12) is an isomorphism. 
From the diagram (3.4.12),  define $\gamma_n(c)$ and $\delta_n(c)$ by the formulae
$$P_c \ ({\rm mod}\ l^n)\ \in \ \Bigg({_{l^n}E(K[c]))}\Bigg)^{{\cal G}_c};
$$ 
$$\gamma_n(c){\rm \ is \ the \ image\ of \ } P_c\ {\rm ( mod\ } l^n{\rm )\ 
in \ } H^1(K, E_{l^n});
$$
$$\delta_n(c){\rm \ is \ the \ image\ of \ } \gamma_n(c)\ \
{\rm in \ } H^1(K, E)_{l^n}.
$$
[This is the same as the construction  in [1, Chap. 7, \S7.14.10], where we write here  $\gamma_n(c), \delta_n(c)$ in place of  $\gamma(c), \delta(c)$ to 
 indicate their dependence on the integer $n$.]

\vskip0.2in\noindent {\bf 3.4.\the\count2 \advance\count2 by 1. Proposition.} (i) {\sl  The order of $\gamma_{n}(c)$ is equal to the order of 
$P_c$ {\rm  (mod $l^n$)} in ${}_{l^n}E(K[c])$.   }

\noindent (ii) {\sl The exponent $t$ of the order $l^t$ of $\delta_{n}(c)$
is the least integer $t$ such that 
$$l^t(P_c\ ({\rm mod} \ l^n))\in { l^nE(K[c])+E(K)\over l^nE([c])}.$$}

\vskip0.2in
\noindent{\it Proof.} Here $P_c$ (mod $l^n$) denotes the image of 
$P_c\in E(K[c])$ in ${}_{l^n} E(K[c])$. These results on the orders of $\gamma_n(c)$ and $\delta_n(c)$ follow immediately from their definition and the  commutative diagram 
(3.4.12). ${\sqcap \!\!\!\!\sqcup}$
\vskip0.2in

\vfil\eject
\noindent {\bf Chapter 4. Structure of the Tate-Shafarevich group and the Selmer group}
\vskip0.4in
\noindent{\bf 4.1. Statement of the main theorems} 
\vskip0.2in
This section contains no proofs. The  main Theorems 4.1.9, 4.1.13, 4.1.14 and 4.1.15
of this paper are finally proved in \S\S5.3-5.5. 
\vskip0.2in\noindent(4.1.1) Throughout this chapter 4,  we  assume that 

\vskip0.2in

$E/F$ is an elliptic curve where $F$ is a global field of characteristic $p>0$;

$I$ is the ideal of $A$
which is the conductor of $E/F$ without the component 

\qquad at $\infty$;

$K$ is an imaginary quadratic extension field of $F$ with respect to $\infty$;

 $l\in \cal P$ is a, necessarily odd, prime number in the set of prime numbers $\cal P$ 

\qquad (see \S3.1);

$\epsilon=\pm1$ is the sign of the functional equation of the $L$-function $L(E/{F},s)$  

\qquad of $E/F$;

$\tau$ is the  element of order $2$ of the galois group Gal$(K/F)$;

$P_c\in E(K[c])$ are the points defined in (3.4.9) over the ring class fields
$K[c]$

\qquad for  all divisors 
$c\in \Lambda(1)$, and where $P_0$ belongs to the group $E(K)$ 

\qquad of $K$-rational points;

 $\alpha(c)$, for a divisor $c\in \Lambda(1)$,  is the largest integer $n$ such that
$c\in \Lambda(n)$ if $c\not=0$ 

\qquad and such that $\alpha(0)=+\infty$ (see remarks 3.2.8).

\vskip0.2in

\noindent Assume that  $E,K,F$ satisfy the following hypotheses (as in [1, Chapter 7, \S7.6.1]):

\vskip0.2in
(a) \qquad $\infty $  is a place  of  $F$    with
 residue field equal to  $k$;

(b)  \qquad  $E/F$    has  split Tate  multiplicative reduction at 
$\infty$   (see  \S1.7

\qquad \qquad or [1, Chapter 4, \S4.7]);

  (c) \qquad   $K$
  is an (separable) imaginary quadratic field  extension of $F$, 

\qquad \qquad   with
respect to $\infty$, 
    such  that  
 all primes dividing  the  conductor $I$ 

\qquad \qquad  split completely in  $K$
 and  $ K\not=F\otimes_{{\mathds F}_q}{\mathds F}_{q^2}$.

\vskip0.2in

 These hypotheses (a),(b),(c) are
assumed for the rest of  this paper.

\vskip0.2in
\noindent(4.1.2) If $M$ is a ${\rm Gal}(K/F)$-module denote by $M^+$  the submodule of $M$ on which $\tau $ acts by $+1$ and $M^-$ is similarly the submodule on which $\tau $ acts by $-1$.

\vskip0.2in\noindent(4.1.3) Let $G$ be a finite abelian $l$-group. 
The {\it invariants} of $G$ are the integers $${r_1}\geq {r_2}\geq {r_3}\geq \ldots$$such that $G$ decomposes into elementary components
$$G\cong {{\mathds Z}\over l^{r_1}{\mathds Z}}\oplus 
{{\mathds Z}\over l^{r_2}{\mathds Z}}\oplus  {{\mathds Z}\over l^{r_3}{\mathds Z}}\oplus \ldots$$
The integers  $r_1,r_2,r_3\ldots$  are uniquely determined by $G$; the integers 
$l^{r_1}, l^{r_2},l^{r_3}, \ldots$ are  also sometimes called the invariants of $G$.

 \vskip0.2in\noindent(4.1.4)  As the elliptic curve $E$ is defined over $F$, the $l$-power torsion subgroup $\coprod\!\!\!\coprod$$(E/K)_{l^\infty}$ of the Tate-Shafarevich group $\coprod\!\!\!\coprod$$(E/K)$ of the elliptic curve $E\times_FK$ over $K$ decomposes  into eigenspaces
$$\coprod\!\!\!\!\coprod(E/K)_{l^\infty} \cong \coprod\!\!\!\!\coprod(E/K)^+_{l^\infty}  \oplus \coprod\!\!\!\!\coprod(E/K)^-_{l^\infty} $$
under the action of the element  $\tau\in $ Gal$(K/F)$.

 The Cassels pairing on $\coprod\!\!\!\coprod$$(E/K)$ is anti-symmetric and respects this decomposition into eigenspaces; furthermore, the Cassels pairing is non-degenerate if this Tate-Shafarevich group $\coprod\!\!\!\coprod$$(E/K)$ is finite; therefore if $\coprod\!\!\!\coprod$$(E/K)$ is finite, the invariants  of the finite abelian $l$-group $\coprod\!\!\!\coprod$$(E/K)_{l^\infty} $ have even  multiplicity.

\vskip0.2in\noindent (4.1.5) 
Under the hypothesis that the Tate-Shafarevich group $\coprod\!\!\!\coprod(E/K)_{l^\infty}$ be finite, let
$$N_1\geq N_3\geq N_5\ldots$$
and $$N_2\geq N_4\geq N_6\ldots$$
be integers such that 
$${N_1},{N_1},{N_3},{N_3}, N_5,N_5\ldots$$ are the invariants of the finite abelian $l$-group $\coprod\!\!\!\coprod(E/K)_{l^\infty}^{\epsilon} $, which is the 
$\epsilon$-eigenspace, and 
$${N_2},{N_2},{N_4},{N_4}, N_6,N_6\ldots$$ are the invariants of the finite abelian $l$-group $\coprod\!\!\!\coprod(E/K)_{l^\infty}^{-\epsilon} $, which is the 
$-\epsilon$-eigenspace.

That is to say, putting $\nu(r)=(-1)^r\epsilon$ then we  have 
isomorphisms for  $r=1$ or $r=2$  
$${\coprod\!\!\!\!\coprod}(E/K)^{-\nu(r)}_{l^\infty}\cong \bigoplus_{s\geq 0}({{\mathds Z}/l^{N_{2s+r}}{\mathds Z}})^2.$$
\vskip0.2in
\noindent {\bf 4.1.6. Definition.} Let $m\geq 1$ be an integer and $c\in \Lambda(m)$. Then the point $P_c$, defined in (3.4.9),  belongs to $E(K[c])$  and the point $P_0$ belongs to $E(K)$.

\noindent  (1) Write $$l^m\vert P_c{\rm \ \ if \ \ } P_c\in l^mE(K[c])$$ and $$l^m\vert\vert P_c{\rm \ \ if \ \
 }P_c\in l^mE(K[c]) \setminus  l^{m+1}E(K[c]).$$

\noindent (2)  For any integer $r\geq 0$,  let $E(r)$ be the  abelian group
$$E(r)=\bigoplus_{c\in \Lambda^r(1)}   {}_{A(c)}E(K[c])$$
where the sum runs over all divisors $c$ of $\Lambda^r(1)$ and where 
$A(c)=l^{\alpha(c)}$. If \break $r>0$ the group $E(r)$ is then a direct sum of finite abelian $l$-groups of
the form $E(K[c])/l^{\alpha(c)}E(K[c])$. If $r=0$ then $E(0)$ is
defined conventionally to be $E(K[0])$.

A point $P_c\in E(K[c])$, where $c\in \Lambda^r(1)$, induces an element in ${}_{A(c)}E(K[c])$ and hence an element of $E(r)$ where all its components in $E(r)$ are zero except possibly for that in ${}_{A(c)}E(K[c])$.

\noindent (3) Let $P(r)$ be the subgroup of $E(r)$ generated by the images in $E(r)$ of the 
points $P_c$ for all $c\in \Lambda^r(1)$. 

Define $M_r$ for any $r\geq 0$ to be the 
largest integer $n\geq 0$ such that $$P(r)\subseteq l^nE(r).$$
If $P(r)=0$ then there is no such largest integer $M_r$ and we then put
$$M_r=+\infty.$$

\vskip0.2in

\noindent {\bf 4.1.7. Conjecture.} For some $r\geq 0$ we have $M_r<+\infty$. 

\vskip0.2in

[See [7, Conjecture A], and also [1, Conjecture 7.14.19] for an explanation of this conjecture.]
\vskip0.2in
\noindent {\it 4.1.8. Remarks.} (1)  
The order of $\gamma_n(c)$ is $l^{n-m}$ for all 
$n\geq m$ if and only if $l^m\vert\vert P_c$ by Proposition 3.4.14.

If $M_r$ is finite then  there is a divisor $c\in \Lambda^r(1)$ such that $\gamma_{{{M_r}+1}}(c)\not=0$ and 
$\gamma_{{{M_r}}}(c)=0$; furthermore, $\gamma_{M_r}(c)=0$ 
for any  $c\in \Lambda^r(1)$ from Proposition 3.4.14(i).

\noindent(2) We have $M_0<+\infty$ if and only if $P_0$ has infinite order 
in $E(K)$.

[This equivalence requires that $E(K)$ has no $l$-torsion; but this requirement holds by the restriction on the prime number $l$  given by  Definition 3.1.3(b).
See Lemma 5.1.2 below for details.] 

\noindent (3) We have that  $\alpha(c)$ is finite if and only if  $c\not=0$ (see Remarks 
3.2.8(2) and (3)).
 Write if $c\not=0$
$${\rm ord}_l(P_c)=\cases{{\displaystyle \max\{ N: l^N\vert P_c\} }&
provided that this maximum is $< \alpha(c)$  \cr
&\cr
{\displaystyle +\infty}& if $\max\{ N: l^N\vert P_c\}$ is $\geq \alpha(c)$.\cr
}$$
Write $${\rm ord}_l(P_0)= \max\{ N: l^N\vert P_0\}$$
which is either an integer or $+\infty$.
If $c\not=0$, we have ord$_l(P_c)<+\infty$ if and only if $l^{{\rm ord}_l(P_c)}\vert\vert P_c$.
and 
ord$_l(P_c)<\alpha (c)$.

\noindent  (4) By definition $M_r$, for all integers $r\geq 0$, 
is given by 
$$M_r=\min\{ {\rm ord}_l(P_c)\vert\ \ 
c\in \Lambda^r(1)\}.$$

Note that  $M_r<+\infty$ if and only if there is  $c\in \Lambda^r(1)$ 
such that  
$${\rm ord}_l(P_c)< \alpha(c)$$with an evident interpretation in the case where $r=0$,
 $c=0$ and
$\alpha(0)=+\infty$.
Furthermore, $M_r<+\infty$ implies that for this selection of 
$c\in \Lambda^r(1)$ we would have $l^{M_r}\vert\vert P_c$ and $c\in \Lambda^r(\alpha(c))$
and $M_r<\alpha(c)$.

Alternatively,  if $r>0$ then we have $M_r<+\infty$ if and only if for some $c\in \Lambda^r(1)$ we have $$c\in \Lambda^r(m+1)\setminus \Lambda^r(m+2)$$ and $${\rm ord}_l(P_c)\leq m.$$

\vskip0.2in
\noindent{\bf 4.1.9. Theorem.} {\sl Suppose that $P_0$ has infinite order in $E(K)$, the group of $K$-rational points of $E$. Let $l$ be 
a prime number in $\cal P$ coprime to the order of {\rm Pic}$(A)$ the Picard group
of $A$. Then the Tate-Shafarevich group ${\coprod\!\!\!\coprod}(E/K)$ is finite and the invariants ${N_i}$, with
multiplicity $2$, of the subgroup ${\coprod\!\!\!\coprod}(E/K)_{l^\infty}$ are given by
$$N_i=M_{i-1}-M_i, {\sl \ \ for \ \ } i\geq 1.$$ That is to say, we have 
the isomorphisms of eigenspaces
$${\coprod\!\!\!\!\coprod}(E/K)^\epsilon_{l^\infty}\cong \prod_{i{\rm \  even}\atop i\geq 0}
({\mathds Z}/l^{M_i-M_{i+1}}{\mathds Z})^2$$
$${\coprod\!\!\!\!\coprod}(E/K)^{-\epsilon}_{l^\infty}\cong \prod_{i{\rm \  odd}\atop
i\geq 0}
({\mathds Z}/l^{M_i-M_{i+1}}{\mathds Z})^2.$$}

\vskip0.2in
\noindent {\it 4.1.10. Remarks.} (1)  The image of the Drinfeld-Heegner point $P_0$ in $E(K)/E(K)_{\rm tors}$ is an eigenvector for $\tau$ and its eigenvalue is equal to $-\epsilon$ (by [1, Chapter 4, Theorem 4.8.6, p. 98], see also [1, Chapter 7, Lemma 7.14.11, p. 388]) where $E(K)_{\rm tors}$ is the torsion subgroup of 
$E(K)$. 

\noindent (2) If the Drinfeld-Heegner point
$P_0$ has infinite order in $E(K)$  then  by [1, Chapter 7, Theorem 7.7.5 and
Remarks 7.7.6(3)] the Tate-Shafarevich group $
\coprod\!\!\!\coprod$$(E/K)$ is finite and therefore the alternating Cassels pairing on 
$\coprod\!\!\!\coprod$$(E/K)$  is non-degenerate and by \S\S4.1.4,4.1.5 the invariants of the $\tau$-eigenspaces of  $\coprod\!\!\!\coprod$$(E/K)_{l^\infty}$ have even multiplicity.

\vskip0.2in
The next  corollary is an immediate consequence of Theorem 4.1.9.\vskip0.2in
\noindent{\bf 4.1.11. Corollary.} {\sl  Under the hypotheses of 
Theorem 4.1.9,   the integers $M_i$ satisfy
$$M_i-M_{i+1}\geq M_{i+2}-M_{i+3}\geq 0, \ {\sl for \ all\ } i\geq 0,$$
and if $j$ is  such that
$$M_{j}=M_{j+1}=M_{j+2}$$ then the sequence $M_j, M_{j+1}, M_{j+2},M_{j+3}\ldots $ is constant.}
${\sqcap \!\!\!\!\sqcup}$

\vskip0.2in
\noindent {\bf 4.1.12. Definition.} Let
$$G=G_1\times\ldots\times G_r$$
be a direct product of finite cyclic groups $G_i$. The 
characters $\chi_1,\ldots,\chi_r$  of $G$ form a {\it triangular basis for the dual  ${\widehat G}$}, relative to the product  $G=\prod_iG_i$, if
they generate ${\widehat G}$ and 
$$\chi_i(G_j)=0, {\rm \ \ for \ all \ }j>i.$$

\vskip0.2in
\noindent{\bf 4.1.13. Theorem.} {\sl   Assume that 
$P_0$ has infinite order in $E(K)$ and $l\in \cal P$ is
a prime number  coprime to the order of the Picard group {\rm Pic}$(A)$. Suppose that
the direct product $$D=\prod_{i\geq 1} D_i$$
is a maximal isotropic subgroup of $\coprod\!\!\!\coprod$$(E/K)_{l^\infty}$ for the Cassels pairing   and $D_i$ is a finite cyclic group of order $l^{N_i}$ for 
all $i\geq 1$ and we have direct products
$$D^{\epsilon }=\prod_{i{\rm \ odd}} D_i$$
$$D^{-\epsilon} =\prod_{i{\rm \ even}} D_i.$$
Then there are effective divisors $c_1\leq c_2\leq \ldots$ on $F$ such that $c_i\in \Lambda^i(M_{i-1})$  for all $i\geq 1$ and the characters
$$d\mapsto <d,\delta_{M_{i-1}}(c_i)>_{\rm Cassels}, \ {\sl \ \ for \ all \ } i$$
form a triangular basis of characters of $D$ relative to the product 
$D=\prod_iD_i$.}
\vskip0.2in
\noindent{\bf 4.1.14. Theorem.} {\sl  Assume that 
$P_0$ has infinite order in $E(K)$ and $l\in {\cal P}$ is  coprime to the order of {\rm Pic}$(A)$.
Let $$M_\infty=\min_{i\in {\bf N}} M_i. $$ 
 Then the group ${\mathds Z}P_0$ has finite index in $E(K)$ and the highest power of $l$ dividing
the index $[E(K):{\bf Z}P_0]$ equals $l^{M_0}$ that is to say
$$\vert (E(K)/{\bf Z}P_0)_{l^\infty}\vert =l^{M_0}.\ \ \ 
$$Furthermore,  we have 
$$\vert\  { {\coprod\!\!\!\!\coprod}(E/K)_{l^\infty}}\ \vert  =l^{2(M_0-M_\infty)}.$$
}

\vskip0.2in\noindent {\bf 4.1.15. Theorem.}
{\sl   Assume that 
$P_0$ has infinite order in $E(K)$ and $l\in {\cal P}$ is  coprime to the order of {\rm Pic}$(A)$. Then for all integers $m\geq 0$, the natural surjection from the Selmer
group to the Tate-Shafarevich group
$$\pi_m: {\rm Sel}_{l^m}(E/K)\longrightarrow 
{\coprod\! \!\!\!\coprod}(E/K) 
_{l^m}\leqno{(4.1.16)}$$
splits and we have isomorphisms of eigenspaces
$${\rm Sel}_{l^m}(E/K)^{\pm} \cong
\big({}_{l^m}E(K) \big)^{\pm}\oplus
 {\coprod\! \!\!\!\coprod}(E/K) 
_{l^m}^\pm {\rm \ \ 
for \ all \ } m\geq 0.\leqno{(4.1.17)}$$
 Furthermore, we have isomorphisms for all $m\geq 0$
$$\matrix{\big(  {}_{l^m}{E(K)}\big)^{-\epsilon}&\cong&
{\mathds Z}/l^m{\mathds Z} \cr
\big(  {}_{l^m}{E(K)}\big)^{\epsilon}&\cong&
 0.\cr}
\leqno{(4.1.18)}$$
}

\vskip0.2in

The next corollary follows immediately from Theorems 4.1.9 and 4.1.15.

\vskip0.2in\noindent {\bf 4.1.19. Corollary.}
{\sl   Under the hypotheses of Theorem 4.1.15,  put
 $${\widehat
N_i}=\min(m,M_{i-1}-M_i)$$ for all $i$ and for all integers $m$. Then for all integers $m\geq 0$, the 
invariants of the   Selmer
$-\epsilon$-eigenspace 
${\rm Sel}_{l^m}(E/K)^{-\epsilon} $ are
$m, {\widehat N_2}, {\widehat N_2},
{\widehat N_4},{\widehat N_4},\ldots$
and the 
invariants of the   Selmer
$\epsilon$-eigenspace 
${\rm Sel}_{l^m}(E/K)^{\epsilon} $ are   
 ${\widehat N_1},{\widehat N_1},{\widehat N_3},{\widehat N_3},\ldots$}
 ${\sqcap \!\!\!\!\sqcup}$

\vskip0.2in

\noindent {\it 4.1.20. Remark. } The main result on the finiteness of Tate-Shafarevich groups proved in the monograph [1, Chapter 7, Theorem 7.6.5 and Theorem 7.7.5] is required for the proofs of the main results of this paper stated in this section; in particular, this  paper does not provide a different proof of  finiteness
independent of the book [1]. This is also the prinicpal reason why the hypotheses (a), (b), (c) of \S4.1.1 are required.

\vskip0.4in 

\eject

\noindent {\bf 4.2. Cochains for the cohomology classes
$\gamma_n(c), \delta_n(c)$ }

\vskip0.2in
The notation and hypotheses of (4.1.1) and (4.1.2) hold in this section.
\vskip0.2in

\noindent {\bf 4.2.1. Lemma.}  \advance\count2 by 1 
{\sl Let $c\in \Lambda(n)$.}

\noindent (i) {\sl The cohomology class $\gamma_n(c)$ is represented by the cocycle
$$\sigma\mapsto -{(\sigma-1)P_c\over l^n}+\sigma{P_c\over l^n}-{P_c\over l^n}, \ \ {\rm Gal}(K^{\rm sep}/K)\to E(K^{\rm sep})_{l^n}$$ where 
${(\sigma-1)P_c\over l^n}$ is the unique $l^n$-division point of $(\sigma-1)P_c$ in $E(K[c])$ and $P_c/l^n$ is a fixed $l^n$-division point of $P_c$.}

\noindent (ii) {\sl The cohomology class $\delta_n(c)$ is represented by the cocycle,
 where 
${(\sigma-1)P_c\over l^n}$ is the unique $l^n$-division point of $(\sigma-1)P_c$ in $E(K[c])$, 
$$\sigma\mapsto -{(\sigma-1)P_c\over l^n}, 
\ \ {\rm Gal}(K^{\rm sep}/K)\to E(K^{\rm sep}) .$$ }

\noindent (iii) {\sl  Let $n\geq m$ be positive integers and let $c\in \Lambda(n)$. Then we have 
$$l^{m} \delta_n(c)=\delta_{n-m}(c)$$
$$l^{m}\gamma_n(c)=\gamma_{n-m}(c).$$}

\vskip0.2in
\noindent{\it Proof of Lemma 4.2.1.} 
(i) and (ii) These two formulae for cocycles representing the cohomology classes $\gamma_n(c)$ and $\delta_n(c)$ can be extracted from  Step 2 of the proof of [1, Chapter 7, Lemma 7.14.14].  We reprove these formulae here.

From the diagram (3.4.12) the restriction homomorphism,
where ${\cal G}_c=$ Gal$(K[c]/K)$,
$$H^1(K,E_{l^n}(K^{\rm sep}))\to H^1(K[c], E_{l^n}(K^{\rm
sep}))^{{\cal G}_c}$$
is an isomorphism as the prime number $l$ belongs to 
$\cal P$. The point ${ P}_c$ belongs to $ E(K[c])$.

Let ${{ P}_c\over l^n}\in E(K^{\rm sep})$ be a fixed  $l^n$th division point of ${ P}_c$,
that is  to say ${{ P}_c\over l^n}$
is any point which satisfies $l^n({{ P}_c\over l^n})={ P}_c$. Then the cocycle
$$\phi:g\mapsto g\big({{ P}_c\over l^n}\big)-{{ P}_c\over l^n}, \ {\rm Gal}(K^{\rm sep}/K[c])\to
E_{l^n}(K^{\rm sep}),$$ represents a cohomology class in $H^1(K[c], E_{l^n}(K^{\rm sep}))^{{\cal
G}_c}$ which is the image of \break $P_c\ ({\rm mod}
\ l^n)\in ({}_{l^n}E(K[c]))^{{\cal G}_c}$ under 
the coboundary map $$\partial_{l^n} : 
({}_{l^n}E(K[c]))^{{\cal G}_c}\to H^1(K[c], E_{l^n}(K^{\rm sep}))^{{\cal G}_c}$$ (see the diagram (3.4.12)). 
The inflation of $\phi$
to ${\rm Gal}(K^{\rm sep}/K)$ is given by the cocycle $$\phi^\sharp :{\rm Gal}(K^{\rm sep}/K)\to
E(K^{\rm sep}), \ \ g\mapsto g\big({{ P}_c\over l^n}\big)-{{ P}_c\over l^n} $$ which need
not necessarily be annihilated by $l^n$.

For any element $g\in {\rm Gal}(K^{\rm sep}/K)$, denote by $${(g-1){
P}_c\over l^n}$$ the unique $l^n$th root of $(g-1){ P}_c$ in $E(K[c])$. This
root exists because $P_c$ belongs to $({}_{l^n}E(K[c])^{{\cal G}_c}$; furthermore,
it is unique  because  $E(K[c])_{l^\infty}=0$ (by Definition 3.1.3(f) and 
Proposition 1.10.1,  [1, Lemma   7.14.11(i)]. 

The cochain 
$$\psi:{\rm Gal}(K^{\rm sep}/K)\to E(K[c]), \ \ g\mapsto -{(g-1){ P}_c\over
l^n}, $$
is a cocycle whose restriction to the subgroup ${\rm Gal}(K^{\rm sep}/K[c])$ is the zero cochain. But
$\psi$ need not be annihilated by 
$l^n$. The cochain $$\phi^\sharp+\psi:{\rm Gal}(K^{\rm sep}/K)\to E(K^{\rm sep}), 
 g\mapsto g\big({{ P}_c\over
l^n}\big)-{{ P}_c\over l^n} -{(g-1){ P}_c\over
l^n}, $$
is a cocycle which is annihilated by $l^n$ and whose restriction to
${\rm Gal}(K^{\rm sep}/K[c])$ is the cocycle 
$$\phi:{\rm Gal}(K^{\rm sep}/K[c])\to E_{l^n}(K^{\rm
sep}).$$
Hence the cochain $\phi^\sharp+\psi$ is a cocycle 
$$\phi^\sharp+\psi:{\rm Gal}(K^{\rm sep}/K)\to E_{l^n}(K^{\rm sep})$$ and this cochain represents
the cohomology class $\gamma_n(c)$ in $H^1(K,E_{l^n})$. Therefore the cohomology
class $\delta_n(c)$ of  $H^1(K[c]/K,E)_{l^n}$ is represented by the cocycle
$$\psi:{\rm Gal}(K[c]/K)\to E(K^{\rm sep}), \ \ g\mapsto -{(g-1){ P}_c\over
l^n}$$
 as required.

\noindent (iii) This follows immediately from the explicit cocycle formulae of parts (i) \break and (ii) of this lemma.
${\sqcap \!\!\!\!\sqcup}$

\vskip0.2in
\noindent{\bf 4.2.2. Lemma.}
 \advance\count2 by 1 {\sl Denote by a superscript $\pm1$ the eigenspaces under the action of   the non-trivial element of 
${\rm Gal}(K/F)$. If $c\in \Lambda^r(n)$ then  we have 
$$P_c \ ({\rm mod\ } l^n)\  \in \Bigg(({}_{l^n}E(K[c]))^{{\cal G}_c}\Bigg)^{-\nu(c)}$$
$$\gamma_n(c)\in H^1(K, E_{l^n})^{-\nu(c)}$$
$$\delta_n(c)\in H^1(K[c]/K, E)_{l^n}^{-\nu(c)}$$
 where ${\cal G}_c=$ Gal$(K[c]/K)$,
$$\nu(c)=(-1)^{r}\epsilon,$$
$r$ is the number of distinct prime divisors in the support of $c$, and $\epsilon$ is the sign in the functional equation of the $L$-function of $E/F$.}

[For the proof, see [1, Lemma 7.14.11]. Note that the set of prime numbers $\cal P$ 
is contained in  the set of prime numbers denoted ${\cal P}\setminus {\cal F}$ of [1, Lemma 7.14.11].] ${\sqcap \!\!\!\!\sqcup}$

\vskip0.4in

\noindent {\bf 4.3. Points $P_c$  defined over local fields}

\vskip0.2in
The notation and hypotheses of (4.1.1)  hold in this section. The 
Drinfeld-Heegner points
$(a,I_1,c)$, $(a,I_1, c,\pi)$ are defined in (3.4.4) and (3.4.5).

\vskip0.2in\noindent{\bf 4.3.1. Proposition.} 
{\sl Suppose that $c\in \Lambda(1)$ and $z\in \Lambda^1(1)$
 is a prime divisor  disjoint from the support of the divisor $c$. Let $y$ be the unique prime of 
$K$ lying over $z$ and
$K_{y}$ be the completion of $K$ at $y$. Then the point $(a,I_1,c)
\in X_0^{\rm Drin}( K[c])$   is definable over $K_{y}$
that is to say
$$(a,I_1,c)\in X_0^{\rm Drin}( K_{y}).$$
Furthermore, we  have that 
$$(a,I_1,c,\pi)\in E(K_{y})$$
is a point of the elliptic curve $E$ definable over $K_{y}$. }

\vskip0.2in
\noindent{\it Proof of Proposition 4.3.1.} The prime $z$ is inert and unramified in $K/F$ by Definition 3.2.4 and remark 3.2.5(i). Furthermore the elliptic curve $E/F$ has good reduction at $z$ by remark 3.2.5(i). As $I$ is the conductor of 
$E/F$ without the component at $\infty$, by (4.1.1), we have that $z$ is disjoint from the support of $I$ and hence
the curve  
$X_0^{\rm Drin}(I)/F
$ also has good reduction  $z$ where there may be several disoint components in 
the closed fibre over $z$.

   By [1, Theorem 4.6.19(ii)] as $z$ is inert and
unramified in $K/F$ the reduction $(a,I_1,c)$ mod $z$ is defined over 
the quadratic extension field $\kappa(y)$ of $\kappa(z)$. That is to say  $(a,I_1,c)$ mod $z$ is a point of the reduction at $z$ of 
$X_0^{\rm Drin}(I)$ which is defined over $\kappa(y)$. 

Let $F_z$ be the completion of $F$ at $z$. As $X_0^{\rm Drin}(I)$
has good reduction at $z$, it follows that the point $(a,I_1,c)$
is defined over the field $K_y$ as this is the unique quadratic extension field
of the local field $F_z$ which is unramified over $z$. It then immediately follows that 
$(a,I_1,c,\pi)$, which is the image of  $(a,I_1,c)$
under the morphism $\pi : X_0^{\rm Drin}\to E$ of $F$-schemes (see (3.4.4) and (3.4.5)),
is a point of the elliptic curve $E$ defined over $K_{y}$.
${\sqcap \!\!\!\!\sqcup}$

\vskip0.2in
\noindent {\bf 4.3.2. Proposition.}  {\sl If $z\in \Lambda^1(1)$ is  a prime divisor disjoint from the support of the divisor $c\in \Lambda(1)$
and $y$ is the prime of $K$ lying over $z$, then the image of $P_c$ in ${}_{l^n}E(K_y)$  via Proposition 4.3.1 is uniquely determined by $P_c$.}

\vskip0.2in
\noindent {\it Proof of Proposition 4.3.2.}  We have the isomorphism
$$K[c]\otimes K_y\cong \prod_i K[c]_{x_i}$$
where $x_i$ are the places of $K[c]$ over the place $y$ of $K$ and $K[c]_{x_i}$
is the completion of $K[c]$ at $x_i$. The galois 
group ${\cal G}_c=$ Gal$(K[c]/K)$ permutes transitively the places $x_i$ and  the completions $K[c]_{x_i}$. We then have
$${}_{l^n}E(K[c]\otimes_KK_y)\cong \prod_i {}_{l^n}E(K[c]_{x_i}).$$

The point $P_c$ belongs to $E(K[c])$ by construction. Hence the point $P_c$ (mod $l^n$)
of ${}_{l^n}E(K[c])$ induces an element $(q_1,q_2,\ldots)$ of 
$
{}_{l^n}E(K[c]\otimes_KK_y)$
where
$$(q_1,q_2,\ldots)\in\prod_i  {}_{l^n}E(K[c]_{x_i})
$$
and $$q_i\in {}_{l^n}E(K[c]_{x_i})\ \ {\rm \ for \ all \ }i.$$

By Proposition 4.3.1, $P_c$ is definable over $K_y$ that is to say 
$P_c\in E(K_y)$ and so we have $q_i\in {}_{l^n}E(K_y)$ for all $i$. But
$P_c$ (mod  $l^n$)  
$ \in \Bigg(  {}_{l^n}E(K[c])\Bigg)^{{\cal G}_c}$
by (3.4.13) (see also Lemma 4.2.2). It follows that $(q_1,q_2,\ldots)$ is invariant
under $
{\cal G}_c$. But the elements $q_i\in {}_{l^n}E(K_y)$ are permuted transitively
by ${\cal G}_c$. Hence the elements $q_i$ are all equal and 
 $P_c$ (mod $l^n$) is the point $(q_1,q_1,\ldots)\in \prod_i{}_{l^n}E([c]_{x_i})$ which is in the image of the diagonal map ${}_{l^n}E(K_y)\to \prod_i {}_{l^n}E(K[c]_{x_i})$. Hence
 the components of the  point $P_c$ (mod $l^n$) in ${}_{l^n} E(K[c]\otimes_KK_y)$ are independent
 of the place $x_i$ and depend only on $P_c$ as required. ${\sqcap \!\!\!\!\sqcup}$

\vskip0.4in

\noindent {\bf 4.4. The map $\chi_z$} \count2=1

\vskip0.2in
The notation and hypotheses of (4.1.1)  hold in this section.

\vskip0.2in\noindent {\bf 4.4.\the\count2. Proposition.}
 \advance\count2 by 1
{\sl  Let $z\in \Lambda^1(n)$ be a prime divisor and let 
${\cal E}_0/\kappa(z)$ be the closed fibre of the N\'eron model of $E/F$ at the place $z$. 
Let $a_z\in {\mathds Z}$ be the trace of the Frobenius at $z$ on the 
inertia invariant part of the Tate module 
of $E$ as in (3.2.6).
Let $y$ be the unique place of $K$ over $z$. 
Let 
$c\in \Lambda(n)$ be a divisor whose support contains   $z$ and put $c'=c-z\in \Lambda(n)$ so that $c'$ has support coprime to $z$.}

\noindent (i) {\sl The endomorphism $\vert G(c/c')\vert {\rm Frob}(z) -a_z
$ of the elliptic curve ${\cal E}_0/\kappa(z)$ annihilates the abelian group 
${\cal E}_0(\kappa(y))$.}

\noindent (ii) {\sl The group homomorphism 
$$h: {}_{l^n}{\cal E}_0(\kappa(y))\ \to 
\ {\cal E}_0(\kappa(y))_{l^n}, \ \ 
x\mapsto \bigg( {  
\vert G(c/c')\vert {\rm Frob}(z) -a_z
\over l^n}\bigg)x$$
is an isomorphism which commutes with $\tau$.}

\vskip0.2in
\noindent{\it Proof of Proposition 4.4.1.}  (i) As in \S1.4, $G(c/c')$ denotes the galois group Gal$(K[c]/K[c'])$. The Frobenius ${\rm Frob}(z)$ acts on the Tate module
of the elliptic curve
${\cal E}_0/\kappa(z)$ with characteristic polynomial
$$X^2-Xa_z+\vert \kappa(z)\vert .$$
Writing $F$ for ${\rm Frob}({z})$ it follows that 
$F^2-Fa_z+\vert\kappa(z)\vert$ annihilates the abelian group
 ${\cal E}_0(\kappa(y))$. But $F^2$ is the identity automorphism on 
${\cal E}_0(\kappa(y))$ as $\kappa(y)/\kappa(z)$ is a quadratic
extension  of finite fields.  Hence
$F^2-Fa_z+F^2\vert\kappa(z)\vert$ annihilates
 ${\cal E}_0(\kappa(y))$. That is to say
$F(F(1+\vert\kappa(z)\vert)-a_z)$
annihilates
 ${\cal E}_0(\kappa(y))$. As $F$ is an automorphism of 
${\cal E}_0({\overline {\kappa(y)}  })$ we obtain that 
$F(1+\vert\kappa(z)\vert)-a_z$
annihilates
 ${\cal E}_0(\kappa(y))$. 
 
As $K\not= F\otimes_{{\mathds F}_q}{\mathds F}_{q^2}$ (by the hypothesis (c) of paragraph (4.1.1))
we have that the unit group
$B^*/A^*$ is the trivial group. The Galois group $$G(c/c')={\rm Gal}(K[c]/K[c'])$$
therefore  has order (see [1, Chap. 2, (2.3.8) and  (2.3.12)])
\begin{eqnarray*}\vert G(c/c')\vert =\vert \kappa(z)\vert+1.\\
\end{eqnarray*}
Hence the endomorphism $\vert G(c/c')\vert {\rm Frob}({z})-a_z$ of 
${\cal
E}_0$  annihilates 
the abelian group ${\cal
E}_0(\kappa(y))$ as required.

\noindent (ii) Let $\alpha,\beta\in \mathds C$ be the complex roots of the characteristic polynomial of
Frobenius ${\rm Frob}(z)$ acting on the Tate module of ${\cal E}_0/\kappa(z)$
$$X^2-Xa_z+\vert \kappa(z)\vert .$$
Then we have by the trace formula for the Frobenius automorphism, where $G(c/c')$ is cyclic of order $\vert \kappa(z)\vert+1$ as in the proof of part (i),
$$\vert {\cal E}_0(\kappa(z))\vert= \vert G(c/c')\vert-\alpha-\beta=
\vert G(c/c')\vert-a_z$$
and
$$\vert {\cal E}_0(\kappa(y))\vert= \vert\kappa(z)\vert^2+1-\alpha^2-\beta^2=(\vert G(c/c')\vert-a_z)(\vert G(c/c')\vert+a_z).$$
We then obtain the decomposition into eigenspaces under the action of 
the involution $\tau$, the non-trivial element
of Gal$(K/F)$, $${\cal E}_0(\kappa(y))\cong {\cal E}_0(\kappa(z))\oplus 
{\cal E}_0(\kappa(y))^{-}$$
where
$$\vert {\cal E}_0(\kappa(y))^{\delta}\vert= \vert
G(c/c')\vert-\delta a_z{\rm \ \ for \ }\delta=\pm1.$$
Hence the order of the $l^\infty$-torsion is given by
$$\vert {\cal E}_0(\kappa(y))^{\delta}_{l^\infty}\vert=
l^{s(\delta)}$$
where $l^{s(\delta)}$ is the highest power of $l$ dividing 
$\delta \vert G(c/c')\vert-a_z$.

By Lemma 3.2.7(ii), we have group isomorphisms for the $l^n$-torsion
$${\cal E}_0(\kappa(y))^{\delta}_{l^n}\cong{{\mathds Z}\over l^n{\mathds Z}}
{\rm \ \ for \ } \delta=\pm1$$
that  is to say these groups are cyclic of order $l^n$.
It follows from this  that there are isomorphisms for the $l^\infty$-torsion
$${\cal E}_0(\kappa(y))^{\delta}_{l^\infty}\cong{{\mathds Z}\over l^{s(\delta)}{\mathds Z}}
{\rm \ \ for \ } \delta=\pm1.$$
We obtain that the $l^\infty$-torsion subgroups of 
${\cal E}_0(\kappa(y))^{-}_{l^\infty}$ and ${\cal
E}_0(\kappa(y))^{+}_{l^\infty}$ are both cyclic and denoting by
$\vert.\vert_l$ the normalised $l$-adic absolute value on ${\mathds Q}$ we have
 $$\vert {\cal
E}_0(\kappa(y))^{\delta}_{l^\infty}\vert = \vert
\ \  \delta \vert G(c/c')\vert-a_z\vert_l^{-1}{\rm \ \ for \ } \delta=\pm1.$$

Let $g$ be the homomorphism
$$g= {\vert G(c/c')\vert {\rm Frob}({z})-a_z\over
l^n}: \ {\cal E}_0\to {\cal E}_0$$
where the integers $\vert G(c/c')\vert$, $a_z$ are both divisible by $l^n$
by Lemma 3.2.7(i).
The homomorphism induced by $g$ on the $\kappa(y)$-rational points
of ${\cal E}_0$
 $$g:{\cal E}_0(\kappa(y))\to {\cal E}_0(\kappa(y))$$
is annihilated by $l^n$ by part (i); furthermore, the subgroup $l^n{\cal E}_0(\kappa(y))$ of ${\cal E}_0(\kappa(y))$ belongs to the
kernel of $g$ again by part (i). Hence $g$ induces a homomorphism
$$h: {}_{l^n}{\cal E}_0(\kappa(y))\to {\cal
E}_0(\kappa(y))_{l^n}.$$

On each eigencomponent $({}_{l^n}{\cal E}_0(\kappa(y)))^\delta$ under the action of
$\tau$, the non-trivial element of Gal$(K/F)$ and 
where $\delta=\pm1$, the map $g$ is multiplication by  the integer
$$ N(\delta)={\delta \vert G(c/c')\vert- a_z\over l^n}
={\delta \vert {\cal E}_0(\kappa(y))^\delta\vert\over l^n}.
$$
It follows from this formula for $N(\delta)$ that the restriction of $h$ to each eigencomponent 
 $({}_{l^n}{\cal E}_0(\kappa(y)))^\delta$ is an injection. 
As $\tau
$ commutes with $h$, the homomorphism $h$ preserves the $\tau$-eigencomponents
and therefore $h$ is an injection.  As
${}_{l^n}{\cal E}_0(\kappa(y))$ and $ {\cal
E}_0(\kappa(y))_{l^n}$ have the same number $l^{2n}$ of elements  it follows
that $h$ is an isomorphism.   ${\sqcap \!\!\!\!\sqcup}$

\vskip0.2in\noindent

\vskip0.2in
\noindent {\bf {\bf 4.4.\the\count2. Proposition}.}  \advance\count2 by 1  {\sl  Let $z$ be a prime divisor of $F$ which belongs to $\Lambda^1(n)$. Let $y\in
\Sigma_K$ be the unique place of $K$ lying over $z$. Assume that the prime number 
$l\in {\cal P}$ is coprime to the order of ${\rm Pic}(A)$, the Picard group of the ring $A$. Then  there is a homomorphism
$$\chi_z:{}_{l^n}E(K_{y})\to H^1(K_{y}, E_{l^n})$$ 
with the following properties.}

\noindent (1) {\sl Let $x$ be a  place of the ring class field $K[z]$ above the place 
$z$ of $K$. Then the image of $\chi_z$ is contained in the subgroup
$H^1(K[z]_x/K_y,E_{l^n}(K[z]_x))$.}

\noindent (2) {\sl The homomorphism  $\chi_z$ is injective.
}

\noindent (3) {\sl The composition of $\chi_z$ with the homomorphism, obtained from
the inclusion of group schemes $E_{l^n}\subset E$,
$$H^1(K_{y}, E_{l^n})\to H^1(K_{y},E)_{l^n}$$is an isomorphism
$$ {}_{l^n}E(K_{y})\cong H^1(K_{y}, E)_{l^n}.$$}

\noindent (4) {\sl For all divisors $c\in \Lambda(n)$ such that $z$ belongs to  ${\rm Supp}(c)$, the cohomology class $\gamma_n(c) $ satisfies}
$$\gamma_n(c)_{y}=\chi_z(P_{c-z}\ ({\rm  mod\ } l^n))$$
{\sl where $P_{c-z}$ ({\rm mod $l^n$}) is an element of 
${}_{l^n}E(K_y)$ (Proposition 4.3.2) and $c-z$ has support coprime to $z$.}

\vskip0.2in
\noindent {\it Proof of Proposition 4.4.2.} For all divisors $d\in \Lambda(n)$ with support
coprime to $z$,  by Propositions 4.3.1 and  4.3.2 the point 
$P_{d}\in E(K[c])$ induces an 
 element of ${}_{l^n}E(K_y)$ which is uniquely 
determined by 
$P_{d}$ where  $y$ is the place of $K$ over $z$ (remark 3.2.5(1)).

Select  a divisor $c\in \Lambda(n)$ with support
containing  $z$ and put
$$c'=c-z\in \Lambda(n).$$
where Supp$(c')$ is coprime to $z$.

Let $z'$ be a prime divisor of $K[c']$ over the place $y$, which is the place of 
$K$ over $z$. The prime divisor
$z^\prime$ is totally ramified in the field extension $K[c]/K[c']$; this follows
via class field theory from the definition of the ring class field $K[c]$
(see \S1.4; more details are given in [1, Chapter 2, (2.3.13)]). 
Let $z^\times$ be the unique prime divisor of $K[c]$ lying over $z^\prime$. 

The 
cohomology class 
$\delta_n(c)$ belongs to $H^1(K[c]/K, E(K[c]))_{l^n}$ which is contained via the inflation map in
$H^1(K, E(K^{\rm sep}))_{l^n}$ (diagram (3.4.12)).

Let ${\cal E}_0/\kappa(z)$ be the closed fibre above $z$ of the N\'eron model of 
$E/F$. We then define a composite isomorphism $\Phi$ as follows where the maps $i,h,j$ are explained below
$$\Phi:{}_{l^n}E(K_y)\ {\buildrel  i\over \longrightarrow}\ {}_{l^n}{\cal E}_0(\kappa(y))
\ {\buildrel
h\over \longrightarrow} \ {\cal E}_0(\kappa(y))_{l^n}\ {\buildrel
j\over \longrightarrow} \ { E}(K_y)_{l^n}.
 \leqno{(4.4.\the\count2)}$$  \advance\count2 by 1
Here $i: {}_{l^n}E(K_y) \longrightarrow{}_{l^n}{\cal E}_0(\kappa(y))$ is the isomorphism obtained from the 
surjective homomorphism $E(K_y) \longrightarrow{\cal E}_0(\kappa(y))$ of reduction at $z$ 
 whose kernel is a pro-$p$-group.  The map
 $h:{}_{l^n}{\cal E}_0(\kappa(y))
 {\longrightarrow}  {\cal E}_0(\kappa(y))_{l^n}$ is the isomorphism of Proposition 4.4.1(ii). The map $j: {\cal E}_0(\kappa(y))_{l^n} { \longrightarrow}  { E}(K_y)_{l^n}$
 is the isomorphism obtained from reduction modulo $y$; the map $j$ is an isomorphism because $E$ has good reduction at $y$ and the prime number $l$ is distinct
 from the characteristic of $F$. The map $\Phi$ is an isomorphism as $i,h ,j$ are
 isomorphisms.

  Let $x'$ be a prime of $K[0]$ lying over $y$ and let $x$ be the unique prime of $K[z]$ over $x'$, where $y$ is 
the prime of $  K$ over $z$. The prime $x$ over $x'$ is uniquely determined by $x'$ because $K[z]/K[0]$
  is totally ramified at $x'$ (see \S1.4).
 The restriction of elements of the Galois group 
$G(c/c')$ to the fields $K[z]$ and $K[z]_x$ induces isomorphisms  $${\rm Gal}(K[z]_x/K[0]_{x'})\cong G(z/0)\cong G(c/c-z)\leqno{(4.4.\the\count2)}$$  \advance\count2 by 1 
where these groups are  cyclic of order $\vert \kappa(z)\vert+1$ (by [1, Chapter 2, (2.3.12), p.20], and as 
$B^*/A^*$ is the trivial group by hypothesis (c) of (4.1.1)). 
Put
$$G={\rm Gal}(K[z]_x/K[0]_{x'}).$$
We have the exact sequence of galois groups
$$0\to G\to {\rm Gal}(K[z]_x/K_{y})
\to {\rm Gal}(K[0]_{x'}/K_{y})\to 0.$$
As the field extension $K[z]/K[0]$ is totally ramified at the place $x'$ of $K[0]$
over $z$, the restriction homomorphism  of $G(c/0)$ to $K[0]_{x^\prime}$ gives the isomorphism
$$G\cong G(z/0)$$
of (4.4.4).

As the field extension $K[z]_x/K[0]_{x'}$ is totally ramified, we have that
the group $G$ acts trivially on $E_{l^n}(K[z]_x)$ and that
$$E_{l^n}(K[z]_x)=E_{l^n}(K[0]_{x'}).$$
Hence we have an isomorphism
$$H^1(G, E_{l^n}(E(K[z]_x))\cong {\rm Hom}(G, E_{l^n}(K[0]_{x'})).\leqno{(4.4.\the\count2)}$$  \advance\count2 by 1

We have already fixed in (3.4.7)  a generator $\sigma_z$ of $G(c/c')$. Let 
 $\sigma\in G$ be the generator induced by $\sigma_z$ under the isomorphism \break $G(c/c')\cong G$ of (4.4.4).

For any $P\in {}_{l^n} E(K_y)$, define a homomorphism
$$f_P: G\to E_{l^n}(K_y)\leqno{(4.4.\the\count2)}$$  \advance\count2 by 1 
as follows. Put, where $\sigma$ is the chosen generator of 
the cyclic group $G$,
$$f_P(\sigma)= \Phi(P).$$
As $\Phi(P)$ is a point of $E_{l^n}(K_y)$  and as the order of the cyclic group 
$G$ is equal to $\vert\kappa(z)\vert+1$ which is divisible by $l^n$, the homomorphism 
$f_P$ is well defined.
Hence $f_P$ defines a cohomology class
in $H^1(G,E_{l^n}(K_y))$ where $G$ acts trivially on
$E_{l^n}(K_y)$. 

As $\Phi: {}_{l^n}E(K_y)\to E(K_y)_{l^n}$
is an isomorphism (see (4.4.3)) and $G$ is cyclic of order divisible by $l^n$,
 this map $P\mapsto f_P$ defines a group isomorphism
$$f: {}_{l^n}E(K_y)\ {\buildrel \cong \over \longrightarrow }\  H^1(G,E_{l^n}(K_y)), \ \ P\mapsto f_P.
\leqno{(4.4.\the\count2)}$$  \advance\count2 by 1

We have the Hochschild-Serre spectral sequence
$$ E^{i,j}_2
\Rightarrow H^{i+j}(K[c]_{z^\times}/K_y, E_{l^n}(K[c]_{z^\times}))$$
where we write
$$G_{z^\times/z'}={\rm Gal}(K[c]_{z^\times}/K[c']_{z'}).$$
and 
$$E^{i,j}_2=H^{i}(K[c']_{z'}/K_y, H^j(G_{z^\times/z'}, E_{l^n}(K[c]_{z^\times}))).$$
The short exact sequence of low degree terms attached to this spectral
sequence in part takes the form
$$0\to E^{1,0}_2\to H^{1}(K[c]_{z^\times}/K_y, E_{l^n}(K[c]_{z^\times}))\to E^{0,1}_2\to
E^{2,0}_2.$$

As the order of Pic$(A)$ is coprime to $l$ we have that the order of the place 
$y$ above $z$ in Pic$(B)$ is coprime to $l$, as $y$ is the unique place of 
$K$ over $z$. It follows that the degree of the 
field extension $K[c']_{z'}/K_y$, which is equal to the degree of the residue
field extensions, is coprime to $l$. Hence we have 
$$H^{i}(K[c']_{z'}/K_y, H^j(G_{z^\times/z'}, E_{l^n}(K[c]_{z^\times})))\cong 0
{\rm \ \ for \ all \ } i\geq 1$$
as the group $H^j(G_{z^\times/z'}, E_{l^n}(K[c]_{z^\times}))$  is  $l$-power torsion
for all $i$. The above short exact sequence of low degree terms then becomes
an isomorphism
$$H^{1}(K[c]_{z^\times}/K_y, E_{l^n}(K[c]_{z^\times}))\cong 
H^{0}(K[c']_{z'}/K_y, H^1(G_{z^\times/z'}, E_{l^n}(K[c]_{z^\times}))).
\leqno{(4.4.\the\count2)}$$  \advance\count2 by 1 
As the field extension $K[c]_{z^\times}/K[c']_{z'}$ is totally ramified, we have 
$ E_{l^n}(K[c]_{z^\times})= E_{l^n}(K[c']_{z'})$ and the group  $G_{z^\times/z'}$ 
acts trivially on $ E_{l^n}(K[c]_{z^\times})$. This last isomorphism of (4.4.8) then  becomes the isomorphism
$$H^{1}(K[c]_{z^\times}/K_y, E_{l^n}(K[c]_{z^\times}))\cong 
H^{0}(K[c']_{z'}/K_y,{\rm  Hom}(G_{z^\times/z'},   E_{l^n}(K[c']_{z'})    ))$$
where $G_{z^\times/z'}$ acts trivially on $E_{l^n}(K[c']_{z'})    $.
This then provides the isomorphism
$$H^{0}(K[c']_{z'}/K_y, H^1(G_{z^\times/z'}, E_{l^n}(K[c]_{z^\times})))
\cong  {\rm Hom}(G_{z^\times/z'}, E_{l^n}(K_y)).
$$ This isomorphism combined with the isomorphism of (4.4.8) gives the isomorphism, where
$G_{z^\times/z'}$ acts trivially on $ E_{l^n}(K_y)$,
$$H^{1}(K[c]_{z^\times}/K_y, E_{l^n}(K[c]_{z^\times}))\cong H^1(G_{z^\times/z'}, E_{l^n}(K_y)).
 \leqno{(4.4.\the\count2)}$$  \advance\count2 by 1

 Take now $c=z$ in the isomorphism of (4.4.9) where $x'$ is a place of $K[0]$ over $y$ and $x$ is the unique place of 
$K[z]$ over $x'$;  we have that 
the composition of this isomorphism (4.4.9) with the isomorphism $f$ of (4.4.7)
gives the isomorphism
$$\phi:{}_{l^n}E(K_y)\ {\buildrel \cong \over \longrightarrow} \ H^1(K[z]_x/K_y, E_{l^n}(K[z]_x)). \leqno{(4.4.\the\count2)}$$  \advance\count2 by 1

The inflation map from Gal$(K[z]_x/K_y)$ to Gal$(K_y^{\rm sep}/K_y)$ is an 
injective homomorphism
$$H^1(K[z]_x/K_y, E_{l^n}(K[z]_x))\
 {\buildrel {\rm inf}\over   \longrightarrow}\ H^1(K_y,E_{l^n}(K^{\rm sep}_y)).
  \leqno{(4.4.\the\count2)}$$  \advance\count2 by 1
The composite of  this 
injective inflation map  of (4.4.11) with the isomorphism $\phi$ of (4.4.10) is then 
 defined to be  the injective homomorphism 
 $\chi_z$  
 $$\matrix{\chi_z:&{}_{l^n}E(K_{y})&\to &H^1(K_{y}, E_{l^n})\cr
 &P&\mapsto& \{\sigma\mapsto \Phi(P)\}.\cr}
 \leqno{(4.4.\the\count2)}$$  \advance\count2 by 1 
 Here $\{\sigma\mapsto \Phi(P)\}$ denotes the cohomology class $f_P$ of 
 (4.4.6) and (4.4.7)
 defined  by 
 $\sigma\mapsto \Phi(P)$ where  $\sigma\in G$ is the chosen generator of $G$; this 
$f_P$ then defines a cohomology class in $H^{1}(K[z]_x/K_y, E_{l^n}(K[z]_x))$ by the isomorphism of (4.4.9) for $c=z$. 
By inflation, this $f_P$ gives a cohomology class 
 of Gal$(K_y^{\rm sep}/K_y)$  and hence an element of the cohomology group
 $H^1(K_{y}, E_{l^n})$ and this defines the homomorphism 
 $\chi_z$. This map $\chi_z$ of (4.4.12) is injective by construction.

  To prove property (1),  by construction the homomorphism
 $$\chi_z:{}_{l^n}E(K_{y})\longrightarrow H^1(K_{y}, E_{l^n})$$
 takes a point $P\in {}_{l^n}E(K_{y})$ to an element $f_P$ of 
$H^1(G, E_{l^n}(K[z]_x))$ which is inflated to an element of 
$H^1(K_{y}, E_{l^n}(K^{\rm sep}))$
and this proves the property  (1) of the map $\chi_z$.

  To prove property (2), 
 the map  $$\chi_z:   {}_{l^n}   E(K_{y})\to
H^1(K_{y}, E_{l^n})^{}, P\mapsto \{\sigma\mapsto \Phi(P)\}$$
is injective by construction, as already noted. 



 
 For property  (3), we have by definition that the map $\chi_z$ of (4.4.12) is the composition of the 
 isomorphism $\phi$ of (4.4.10) with the inflation map of (4.4.11). That is to say,  the map 
 $\chi_z$ factors as 
 $${}_{l^n}E(K_y)\ {\buildrel \phi \over  \longrightarrow } \ H^1(K[z]_x/K_y, E_{l^n}(K[z]_x)) \
 {\buildrel {\rm inf}\over   \longrightarrow}\ H^1(K_y,E_{l^n}(K^{\rm sep}_y))
  \leqno{(4.4.\the\count2)}$$  \advance\count2 by 1 
where $\phi$ is an isomorphism and inf is an injection. 

We have an exact sequence obtained from 
the inclusion of group schemes $E_{l^n}\subset E$,
$$0\ \longrightarrow\  {}_{l^n}E(K_y)\  {\buildrel \partial_{l^n}\over 
\longrightarrow}\ H^1(K_{y}, E_{l^n})\ {\buildrel \psi\over \longrightarrow} \ H^1(K_{y},E)_{l^n}.$$
The morphism of multiplication by $l^n$ on  the N\'eron model of $E$ over the ring of valuation integers of $F$ 
at the place
$z$ is \'etale; therefore the image
of $\partial_{l^n}$ consists of unramified cohomology classes and more precisely the 
image of $\partial_{l^n}$ belongs to the subgroup $H^1(K_y^{\rm nr}/K_y, E_{l^n})$
of $H^1(K_y, E_{l^n})$ where $K_y^{nr}$ is the maximal separable unramified extension
of $K_y$. The intersection of $H^1(K_y^{\rm nr}/K_y, E_{l^n})$ with 
$
H^1(K[z]_x/K_y, E_{l^n}(K[z]_x))$, which are both subgroups of $H^1(K_y, E_{l^n})$, 
is therefore contained in $H^1(K[0]_{x'}/K_y,E_{l^n}(K[0]_{x'}))$; this follows by considering representative 
cocycles and because  
the field extension $K[z]_x/K[0]_{x'}$ is totally ramified and 
$K[0]_{x'}/K_y$ is unramified. But  the order of the Picard group Pic$(A)$ is coprime to $l$ by hypothesis; therefore the order of the place  $y$ over $z$ in Pic$(B)$
is coprime to $l$ as $y$ is the unique place of $K$ over $z$. It follows that the degree of the field extension
$K[0]_{x'}/K_y$ is coprime to $l$ and therefore we have
$$H^1(K[0]_{x'}/K_y,E_{l^n}(K[0]_{x'}))=0.$$
It follows that the composition of the injective inflation map 
$$ H^1(K[z]_x/K_y, E_{l^n}(K[z]_x)) \
 {\buildrel {\rm inf}\over   \longrightarrow}\ H^1(K_y,E_{l^n}(K^{\rm sep}_y))
  $$ 
  with $\psi: H^1(K_{y}, E_{l^n})\to  H^1(K_{y},E)_{l^n}$ is an injection
  $$ H^1(K[z]_x/K_y, E_{l^n}(K[z]_x)) \
 {\buildrel \over   \longrightarrow}\ H^1(K_{y},E)_{l^n}.
 \leqno{(4.4.\the\count2)}$$  \advance\count2 by 1

By arithmetic flat local duality, there is an isomorphism of discrete groups
$$H^1(K_y,E)\cong {\rm Hom}( E(K_y),{\mathds Q}/{\mathds Z})
 \leqno{(4.4.\the\count2)}$$  \advance\count2 by 1 
for the proof of which, using flat cohomology and local class field theory, see [9, Chapter 3, Theorem 7.8]; we only require the 
prime-to-$p$ version of this duality (4.4.15) and whose simpler proof is explained in [9, Chapter 1, Remark 3.6]. From (4.4.15) we obtain the 
isomorphism of discrete groups
$$H^1(K_y,E)_{l^n}\cong {\rm Hom}( {}_{l^n}E(K_y),{\mathds Z}/l^n{\mathds Z}).
 \leqno{(4.4.\the\count2)}$$  \advance\count2 by 1 
 It follows from (4.4.16) that $H^1(K_y,E)_{l^n}$ and ${}_{l^n}E(K_y)$ have the same number of elements. As $E_{l^n}(K_y)$ is isomorphic to ${}_{l^n}E(K_y)$
 by the isomorphism $\Phi$ of (4.4.3), we have that 
 $H^1(K_y,E)_{l^n}$ and $E_{l^n}(K_y)$ have the same number of elements.
 
 From (4.4.9) for $c=z$,  we have the isomorphism, where $G$ is cyclic of order 
 divisible by $l^n$ and which acts trivially on $E_{l^n}(K_y)$,
 $$H^{1}(K[z]_x/K_y, E_{l^n}(K[z]_x))\cong H^1(G, E_{l^n}(K_y))
 \cong E_{l^n}(K_y).$$ 
 
 It now follows that the finite groups 
 $H^{1}(K[z]_x/K_y, E_{l^n}(K[z]_x))$, $E_{l^n}(K_y)$ and $H^1(K_y,E)_{l^n}$ have the same number of elements and hence the injective homomorphism  of (4.4.14)
 $$H^1(K[z]_x/K_y, E_{l^n}(K[z]_x))\longrightarrow H^1(K_{y}, E)_{l^n}$$
is an isomorphism. It follows that the natural 
homomorphism, obtained from
the inclusion of group schemes $E_{l^n}\subset E$,
$$H^1(K_{y}, E_{l^n})\ {\buildrel \psi\over \longrightarrow}\  H^1(K_{y},E)_{l^n}$$
composed with   $\chi_z$, where $\chi_z$ factors through the group $H^{1}(K[z]_x/K_y, E_{l^n}(K[z]_x))$ as in (4.4.13),
$${}_{l^n}E(K_y)\ {\buildrel \chi_z\over \longrightarrow} \ 
H^1(K_y,E_{l^n})\ {\buildrel \psi\over \longrightarrow }\
H^1(K_y,E)_{l^n}$$
is an isomorphism
$$ {}_{l^n}E(K_{y})\cong H^1(K_{y}, E)_{l^n}.$$
This proves the property (3).

 It only remains to prove property (4).
For all divisors $d\in \Lambda(n)$,  as in (3.4.10) we shall write $P_d$ (mod $l^n$) for the image of $P_d\in E(K[d])$ in
 ${}_{l^n}E(K[d])$.

 From Lemma 4.2.1, we have  that the cohomology class $\gamma_n(c)$ in $H^1(K, E_{l^n})$  is represented by the cocycle, where 
${(\sigma-1)P_c\over l^n}$ is the unique $l^n$-division point of $(\sigma-1)P_c$ in $E(K[c])$ and $P_c/l^n$ is a fixed $l^n$-division point of $P_c$,
$$\Gamma_n(c):\sigma\mapsto -{(\sigma-1)P_c\over l^n}+\sigma{P_c\over l^n}-{P_c\over l^n}, \ \ {\rm Gal}(K^{\rm sep}/K)\to E(K^{\rm sep})_{l^n}
 \leqno{(4.4.\the\count2)}$$  \advance\count2 by 1 
and the cohomology class $\delta_n(c)$ in $H^1(K, E)_{l^n}$  is represented by the cocycle
$$\Delta_n(c): \sigma\mapsto -{(\sigma-1)P_c\over l^n}, \ \ 
{\rm Gal}(K^{\rm sep}/K)\to E(K^{\rm sep}). \leqno{(4.4.\the\count2)}$$  \advance\count2 by 1

 The cohomology class $\delta_n(c)_{y}\in H^1(K_{y},
 E(K_{y}^{\rm sep}))$ is the restriction at $y$ of the class
$\delta_n(c)\in H^1(K,E(K^{\rm sep}))$. Let $Q\in E(K[c])$ be the 
element  given by  
$$Q=-{(\sigma_z-1){ P}_c\over l^n}  \leqno{(4.4.\the\count2)}$$  \advance\count2 by 1 
and where  $\sigma_z$ is the fixed generator of the cyclic group $G(c/c')$.

The
definition of  $P_c$ (as in [1, Chap. 7, \S7.14.10 pp.386-387] and (3.4.9) above) is the following. Let 
$c=\sum_{i=1}^r z_i$ be the decomposition  of $c$ as
a sum of distinct  prime divisors, where we write $z=z_1$. 
Then we have
 $$P_c=\sum_{s\in {\cal S}} sD_cy_c$$
 where
 $$y_c=(0,I_1,c,\pi)\in E(K[c])\leqno{(4.4.\the\count2)}$$  \advance\count2 by 1 
 as in (3.4.7) and
 where $\cal S$ is a set of coset representatives for $G(c/0)$ in 
 Gal$(K[c]/K)$ and $D_c$ is the Kolyvagin element of (3.4.9).

We have an exact sequence of abelian groups
$$0\to G(c/c')\to G(c/0)\to G(c'/0)\to 0.$$
We obtain a corresponding decomposition in the group algebra
${\mathds Z}[G(c/0)]$ of the Kolyvagin element $D_c$
$$D_c=D_{c-z}D_z$$

We have already selected in (3.4.7) a generator  $\sigma_{z}$  of the cyclic group
$G(c/c-z)$. 
We may then define a map of sets $$h_{z}:G(c/c-z)\to \mathds Z$$ by
$$\sigma_{z}^{-s}\mapsto -s,{\rm \ \ for \ } s=0,1,\ldots,\vert 
G(c/c-z)\vert-1.$$
The Kolyvagin element of $h_{z}$ is then as in (3.4.8)
$$D_{z}=-\sum_{r=1}^{\vert G(c/c-z)\vert-1} r\sigma_{z}^r.$$
We have
$$(\sigma_{z}-1)D_{z}=-\vert G(c/c-z)\vert+ e_{G(c/c-z)}\leqno{(4.4.
\the\count2)}$$ \advance\count2 by 1
where $e_{G(c/c-z)}$ is the element of the group 
algebra ${\mathds Z}[G(c/c-z)]$ given by
$$ e_{G(c/c-z)}=\sum_{g\in  G(c/c-z)}g.$$

The element ${ P}_c\in E(K[c])$ may then be written as,  where we write $y_{c}=(0,I_1,c,\pi)$ as in (3.4.7) and (4.4.20),  
$${ P}_c=\sum_{s\in {\cal S}} s D_{c-z}D_z  y_{c}.\leqno{(4.4.
\the\count2)}$$ \advance\count2 by 1

We have (see {[1, Chap. 4, (4.8.3) and table
4.8.5}])  $${\vert O_{c'}^*\vert\over \vert
A^*\vert}{\rm Tr}_{K[c]/K[c']}y_{c}=a_zy_{c'}$$
where $a_z\in {\mathds Z}$ is as in Proposition 4.4.1 and (3.2.6) and where $y_{c'}=(0,I_1,c',\pi)$.
By definition $Q=-((\sigma_z-1){ P}_c)/l^n$ (see (4.4.19)); hence  we have from (4.4.21)
$$Q=-\sum_{s\in{\cal S}} s D_{c-z}  \big( 
{(\sigma_z-1)D_z\over l^n}\big)y_{c} $$
$$= \sum_{s\in{\cal S}} s D_{c-z}  \big( 
{\vert G(c/c')\vert \over l^n}y_{c}-{ \vert
A^*\vert\over \vert O_{c'}^*\vert}{a_z\over
l^n}y_{c'}\big).$$
As $K$ is not obtained from $F$ by ground field extension (hypothesis (c) of paragraph (4.1.1)), we
have $
A^*= O_{c'}^*$; hence we obtain
 $$Q= \sum_{s\in{\cal S}} s D_{c-z}  \big( 
{\vert G(c/c')\vert \over l^n}y_{c}-{a_z\over
l^n}y_{c'}\big).\leqno{(4.4.\the\count2)}$$  \advance\count2 by 1

As at the beginning of this proof, let $y$ be the unique place of $K$ over $z$, $z'$ be a prime of $K[c']$ over the place $y$ of $K$, and 
let $z^\times$ be the place of $K[c]$ over the prime $z'$ of $K[c']$ where the field extension 
$K[c]/K[c']$ is totally ramified at $z'$. Also ${\cal E}_0/\kappa(z)$ denotes the closed fibre above 
$z$ of the N\'eron model of $E/F$.

We write ${ Q}_0$ for the image of $Q$ modulo $z^\times$ in ${\cal
E}_0(\kappa(z^\times))$ by passage to the residue field
$\kappa(z^\times)$. From the isomorphism (4.4.9), where $\delta_n(c)_y$ belongs to 
$H^{1}(K[c]_{z^\times}/K_y, E_{l^n}(K[c]_{z^\times}))$, we have that the reduction 
 of $\delta_n(c)_y$ at $y$ belongs to   ${\rm Hom}(G(c/c'), {\cal
E}_0(\kappa(y)))_{l^n}$ and  is given by the
  cocycle (see (4.4.18))
 $$g\mapsto -{(g-1){ P}_c\over l^n}\ \ ({\rm mod} \ z^\times),
 \ \ {\rm Gal}(K[c]_{z^\times}/K[c']_{z'})\cong G(c/c')\to {\cal
E}_0(\kappa(y))$$ 
and we have that $-{(g-1){ P}_c\over l^n}$ modulo  $z^\times$
belongs to the $l^n$-torsion subgroup  ${\cal
E}_0(\kappa(y))_{l^n}$ rational over $\kappa(y)$. Hence the point $ Q_0$, the reduction of 
$Q$ modulo $z^\times$, belongs to 
${\cal
E}_0(\kappa(y))_{l^n}$. Note that 
$-{(g-1){ P}_c}$ modulo  $z^\times$
reduces to zero, for all $g\in G(c/c')$,  as $K[c]/K[c']$ is totally
ramified at $z'$.

Denote by  Frob$(z)$  the Frobenius
automorphism $x\mapsto x^{\vert\kappa(z)\vert}$ of the closed fibre  ${\cal E}_0/\kappa(z)$
 over $z$ of the
N\'eron model of $E/F$. The theorem {[1, Chap. 4, Theorem 4.8.9] } gives that for the prime
$z^{\prime}$ of $K[c']$ above $z$ we have, where $y_{c}=(0,I_1,c,\pi)$
 as in (4.4.20), $${\rm
Frob}(z)y_{c}\equiv y_{c'}\ \ ({\rm mod} \ \
z^{\prime}).\leqno{(4.4.\the\count2)}$$  \advance\count2 by 1

     We obtain from Proposition 4.3.1 that $y_{c}$ mod $z^\times$ is defined
over the subfield $\kappa(y)$  of $\kappa(z^\times)$ where $\kappa(y)$ is the quadratic extension field 
of the finite field $\kappa(z)$.
 Hence we have from (4.4.24)
 $${\vert
G(c/c')\vert \over l^n}y_{c}-{a_z\over l^n}y_{c'}\equiv{\vert G(c/c')\vert {\rm
Frob}({z})-a_z\over l^n}y_{c'}\ ({\rm mod}\ z^{\prime}).
\leqno{(4.4.\the\count2)}$$  \advance\count2 by 1

The point ${ P}_{c'}\in E(K[c'])$ is given by (see (3.4.9))
$${ P}_{c'}=\sum_{s\in {\cal S}}s
D_{c-z} y_{c'}.$$ 
We then have from (4.4.23) and (4.4.25)
$$Q=\sum_{s\in{\cal S}} s D_{c-z}  \big( 
{\vert G(c/c')\vert \over l^n}y_{c}-{a_z\over
l^n}y_{c'}\big) $$
$$\equiv \ {Q_0}\ \equiv{\vert G(c/c')\vert {\rm Frob}({z})-a_z\over l^n}{
P}_{c'}\ \ ({\rm mod}\ z^{\prime}). \leqno{(4.4.\the\count2)}$$\advance\count2 by 1

 From  Lemma 4.2.2 or [1, Lemma 7.14.11(ii)]), we have that  ${ P}_{c'}$ (mod  $l^n$) belongs to the
$-\nu({c'})$-eigenspace for $\tau$ on ${}_{l^n}{
E}(K[c'])$, where $\nu(c')=(-1)^r\epsilon$ and $r$ is the number of prime divisors in the support of $c'$ and $\epsilon$ is the sign in the functional equation of the 
$L$-function of $E/F$ (as in Lemma 4.2.2). As $z$ is inert in $K/F$, it follows that the
image of the  reduction
$P^\flat_{c'} $  of ${ P}_{c'}$ modulo  $z^{\prime}$ belongs to the
$-\nu({c'})$-eigenspace for $\tau$ on ${}_{l^n}{\cal
E}_0(\kappa(z^{\prime}))$. Let 
$$h: {}_{l^n}{\cal E}_0(\kappa(y))\ \to 
\ {\cal E}_0(\kappa(y))_{l^n}, \ \ 
x\mapsto \bigg( {  
\vert G(c/c')\vert {\rm Frob}(z) -a_z
\over l^n}\bigg)x$$
be the isomorphism  which commutes with $\tau$ of Proposition 4.4.1. As   ${ Q}_0$ belongs to the subgroup ${\cal
E}_0(\kappa(y))$ as already noted and also as $P^\flat_{c'}$ (mod $l^n$) 
is an element of ${}_{l^n} {\cal E}_0(\kappa(y))$ by Proposition 4.3.2, we have 
by (4.4.26) $$  {
Q}_0=h(P^\flat_{c'} \ ({\rm mod}\ l^n)). $$
With the notation of (4.4.3) we then have
$$  {
Q}_0=h\circ i(P_{c'} \ ({\rm mod}\ l^n))$$
where $i$ is the reduction isomorphism ${}_{l^n}E(K_y)\to {}_{l^n}{\cal E}_0(\kappa(y))$
and where $P_{c'}$ (mod $l^n$) belongs to ${}_{l^n} E(K_y)$ by Proposition 4.3.2.
Furthermore, by (4.4.3) and as $\Phi=j\circ h\circ i$ we have
$$\Phi(P_{c'} \ ({\rm mod}\ l^n)) =j(Q_0) \leqno{(4.4.\the\count2)}$$\advance\count2 by 1
is the unique $l^n$-torsion point of $E(K_y)$   whose reduction at $y$ is $Q_0$.

Let ${\overline {z} }$ be any prime of ${K^{\rm sep}}$ above ${z^\times} $. Restrict the cocycle 
$$\Gamma_n(c): 
{\rm Gal}(K^{\rm sep}/K)\to E(K^{\rm sep})_{l^n}$$
of (4.4.17), which represents $\gamma_n(c)$,
to the decomposition group of ${\overline {z} }$. As $K[c]/K[c']$ is totally ramified at $z^\prime$ we have that
 the Kolyvagin element $D_z$ restricted to the 
the residue field of $z^{\times}$ satisfies
$$D_z=-\vert\kappa(z)\vert(\vert\kappa(z)\vert+1)/2.$$ 
Furthermore,
$l^n$ divides $\vert\kappa(z)\vert+1$ and $l$ is different from $2$; hence we have that (from (4.4.22))
the reduction $P_c^\flat$ of $P_c$ at $z^{\times}$ satisfies
$P_c^\flat \in l^n{\cal E}_0(\kappa(z^{\times}))$ and hence as $E$ has good reduction at $z^{\times}$ we have that
$$P_c\in l^nE(K[c]_{z^{\times}}).$$

 We then  have that  the 
cocycle $\Gamma_n(c) $ of (4.4.17) satisfies $$\Gamma_n(c)(\rho)=0 {\rm \ \ for \ all \ } \rho\in 
{\rm Gal} (K[c]_{z^{\times}}^{\rm sep}/K[c]_{z^{\times}})$$  and this cocycle, restricted to the decomposition group
of ${\overline z}$, factors through the subgroup
Gal$(K[c]_{z^{\times}}/K_{y} )$. Furthermore, since $\sigma_z$, which is the selected generator of $G(c/c')$, is in the 
inertia group of ${y} $ we have that
$$\sigma_z{P_c\over l^n}-{P_c\over l^n}$$ reduces to zero modulo ${z^{\times}}$. Hence
$$\Gamma_n(c)(\sigma_z)= -{(\sigma_z-1)P_c\over l^n}+\sigma_z{P_c\over l^n}-{P_c\over l^n}$$
is the unique $l^n$-torsion point congruent to $Q$ (mod ${z^{\times}}$) where $Q$ is given by 
(4.4.19)
that is to say 
$$\Gamma_n(c)(\sigma_z)\equiv -{(\sigma_z-1)P_c\over l^n}\ \ 
({\rm mod}\ \  {z^{\times}}). \leqno{(4.4.\the\count2)}$$\advance\count2 by 1
 But from  (4.4.27),  we then obtain 
$$ \Gamma_n(c)(\sigma_z)= \Phi(P_{c'} \ ({\rm mod}\ l^n)).
 \leqno{(4.4.\the\count2)}$$\advance\count2 by 1
Then  (4.4.29) shows, as $\Gamma_n(c)(\sigma_z)$ is the unique $l^n$-torsion point 
with reduction at $y$ coinciding with the reduction of   $Q=-{(\sigma_z-1)P_c\over l^n}$ at $y$, that we have an equality of cohomology classes in $H^1(K_y,E_{l^n})$
$$\gamma_n(c)_y=\chi_z(P_{c-z}\ ({\rm mod} \ l^n))$$
where $P_{c-z}$ (mod $l^n$) is an element of ${}_{l^n}E(K_y)$.
This proves 
property (4) and completes the proof of 
Proposition 4.4.2.  ${\sqcap \!\!\!\!\sqcup}$

\vskip0.2in

\vskip0.2in
\noindent {\it 4.4.\the\count2. Remarks.} (1) \advance\count2 by 1
This Proposition 4.4.2 and its consequence stated  in Proposition 4.5.1(iv) below are 
an extension of [1, Chapter 7, Lemma 7.14.14(ii)]. It would be 
of
interest to eliminate the hypothesis
that the prime number $l$ be coprime to the order of the Picard group Pic$(A)$
from these Propositions 4.4.2 and 4.5.1(iv).

\noindent (2)  The homomorphism  $\chi_z$  interchanges the 
$\tau$-eigenspaces in that 
for $\delta=\pm1$ we have
$$\chi_z({}_{l^n}   E(K_{y}))^\delta)\subseteq
H^1(K_{y}, E_{l^n})^{-\delta}.$$
This property, which  is not required for this paper, follows from the group Gal$(K[c]/F)$ being generalised dihedral.

\vskip0.4in

\noindent {\bf 4.5. Localizations of the classes $\gamma_n(c)$ and $\delta_n(c)$}
\vskip0.2in
The notation and hypotheses of (4.1.1) hold in this section. Let 
$\coprod\!\!\!\coprod(E/K)$ be the Tate-Shafarevich group of the elliptic curve 
$E\times_FK$ over $K$.
\vskip0.2in
\noindent {\bf 4.5.1. Proposition.}   {\sl 
Let $c\in \Lambda(n)$, let $z$ be a prime divisor in 
$\Lambda^1(n)$ and let $y$ be the place of $K$  over the place $z$ of $F$.}

\noindent (i) {\sl Let $v$ be a place of $K$ coprime to $c$. Then we have $\gamma_n(c)_v\in \partial_n(E(K_v))$ that is to say we have $\delta_n(c)_v=0$.}

\noindent (ii) {\sl  If   $c$ is coprime to $z$ then we have 
$\gamma_n(c)_{y}=\partial_n((P_{c} \ ({\rm mod } \ l^n))_y).$}

\noindent (iii) {\sl If $l^n\vert P_c$ then we have $\delta_n(c)=0$. If $l^n\vert 
P_{c-w}$ for all prime divisors  $w$ in the support of $c$ then we have
 $\delta_n(c)\in \coprod\!\!\!\coprod(E/K)_{l^\infty}$.}
 
 \noindent (iv) {\sl Assume that the prime number 
$l\in {\cal P}$ is coprime to the order of ${\rm Pic}(A)$, the Picard group of the ring $A$. If $z\in {\rm Supp}(c)$  then we have
$${\rm ord}\  \delta_n(c)_{y}={\rm ord}\ \gamma_n(c)_{y}={\rm ord}\ \gamma_n(c-z)_{y}=$$
$${\rm ord}\bigg( (P_{c-z}   \ ({\rm mod} \ l^n))_y  \bigg)$$
where   $(P_{c-z}\  ({\rm mod}\  l^n))_y$
is  an element of 
${}_{l^n}E(K_y)$.}

\vskip0.2in
\noindent {\it Proof.} (i) The field $K[c]$ is a subfield of $K_\infty$ as 
$\infty$ is split completely in $K[c]/K$ (see [1, Chapter 2, (2.3.13)]). We have
that the localization $\delta_n(c)_\infty$ at $\infty$ of 
$\delta_n(c)$ satisfies
$\delta_n(c)_\infty\in H^1(K_\infty, E)_{l^n}$ and the
localization $ ( P_c\ ({\rm mod}\
l^n))_\infty$ at $\infty$ of 
 $ P_c\ ({\rm mod}\
l^n)$ satisfies  $( P_c\ ({\rm mod}\
l^n))_\infty\in
{}_{l^n}E(K_\infty)$.  It then follows from the diagram (3.4.12) that $\delta_n(c)_\infty=0$.

Suppose now that $v$ is a place of $K$ such that $v\not=\infty$ and 
$v\not\in $ Supp$(c)$. We have by construction that 
$$\delta_n(c)\in H^1(K[c]/K, E(K[c]))_{l^n}.$$
The field extension $K[c]/K$ is unramified at $v$ (see [1, Chapter 2, (2.3.13)]). Hence we have that the localization $\delta_n(c)_v$ at $v$ satisfies
$$ \delta_n(c)_v\in H^1(K_v^{nr}/K_v,E(K_v^{nr}))_{l^n}\subseteq
H^1(K_v, E(K_v^{\rm sep}))_{l^n}$$
where $K_v^{nr}$ is the maximal unramified separable extension of 
the local field $K_v$. But $H^1(K_v^{nr}/K_v,E(K_v^{nr}))_{l^n}=0$ by Definition 3.1.3(d) of the set of prime numbers $\cal P$ to which $l$ belongs. Hence we have 
$\delta_n(c)_v=0$; this last vanishing is equivalent to  $\gamma_n(c)_v\in \partial_n(E(K_v))$, as required.

\noindent (ii) By Lemma 4.2.1, the cohomology class $\gamma_n(c)$ is represented by the cocycle
$$\sigma\mapsto -{(\sigma-1)P_c\over l^n}+\sigma{P_c\over l^n}-{P_c\over l^n}, \ \ {\rm Gal}(K^{\rm sep}/K)\to E(K^{\rm sep})_{l^n}$$ where 
${(\sigma-1)P_c\over l^n}$ is the unique $l^n$-division point of $(\sigma-1)P_c$ in $E(K[c])$. Furthermore by the same Lemma 4.2.1, the cohomology class $\delta_n(c)\in H^1(K, E(K^{\rm sep}))_{l^n}$
is represented by  
the cocycle
$$\sigma\mapsto -{(\sigma-1)P_c\over l^n}, \ \ {\rm Gal}(K^{\rm sep}/K)\to E(K^{\rm sep}).$$
By part (i) as $z\not\in $ Supp$(c)$, we have $\delta_n(c)_y=0$ hence a 
cocycle representing $\delta_n(c)$ when 
localized at $v$ is cohomologous to zero. Therefore  there is an 
element 
$a\in E(K_y^{\rm sep})$ such that 
$$-{(\sigma-1)P_c\over l^n} = \sigma a-a, \ \ {\rm  for \ all \ } \sigma\in \ {\rm Gal}(K^{\rm sep}_y/K_y).$$
Hence we have
$$-(\sigma-1)(P_c+l^n a) =0, \ \ {\rm  for \ all \ } \sigma\in \ {\rm Gal}(K^{\rm sep}_y/K_y).$$
This implies that 
$$P_c+l^n a\in E(K_y).$$

By Lemma 4.3.1, the  localization at $y$ of the point $P_c$  lies in $ E(K_y)$
where $y$ is the place of $K$ over $z$. Furthermore, by Lemma 4.3.2 the point $(P_c{\rm \  (mod \ }l^n))_y $
belongs to 
${}_{l^n} E(K_y)$ and its image in this group is uniquely determined by $P_c$.

It follows from this and the above cocycle formulae for 
$\gamma_n(c)$ and $\delta_n(c)$ that the cohomology class $\gamma_n(c)_y$ is represented 
by the cocyle
$$\sigma\mapsto \sigma{P_c\over l^n}-{P_c\over l^n}, \ \ {\rm Gal}(K_y^{\rm sep}/K_y)\to E(K_y^{\rm sep})_{l^n}.$$
That is to say we have $\gamma_n(c)_{y}=\partial_n((P_{c} \ ({\rm mod } \ l^n))_y)$.

\noindent (iii) From Lemma 4.2.1(ii), the cohomology class $\delta_n(c)$ is 
represented by the cocycle
$$\sigma\mapsto -{(\sigma-1)P_c\over l^n}, 
\ \ {\rm Gal}(K^{\rm sep}/K)\to E(K^{\rm sep}) .$$ 
We then evidently have that if $l^n\vert P_c$, that is to say $P_c\in l^nE(K[c])$, then  $\delta_n(c)=0$.

  Suppose that   $l^n\vert P_{c-z}$ for all prime 
divisors $z$ in the support of $c$. Then for any prime divisor $y$ of $K$ lying over
the prime divisor $z$ dividing $c$ we have $\delta_n(c)_{y}=0$
by Proposition 4.4.2(4). It then follows from part (i) of the present Proposition 4.5.1 that 
$\delta_n(c)_{v}=0$ for all places $v$ of $K$ and hence that 
$\delta_n(c)$ belongs to the Tate-Shafarevich group $\coprod\!\!\!\coprod(E/K)_{l^\infty}$.

\noindent (iv) From
  the property (4) of the map $\chi_z$ of Proposition 4.4.2, we have
$$\gamma_n(c)_y=\chi_z((P_{c-z}\ ({\rm mod} \ l^n))_y)$$
where $(P_{c-z}\ ({\rm mod} \ l^n))_y$ is an element of ${}_{l^n}E(K_y)$ by Proposition 
4.3.2.
As the homomorphism $\chi_z$ is injective (property (2) of Proposition 4.4.2), 
  we obtain
$${\rm ord}(\gamma_n(c)_y)={\rm ord}((P_{c-z}\ ({\rm mod}\ l^n))_y)$$
where again $(P_{c-z} \ ({\rm mod}\ l^n))_y$ denotes an element of 
${}_{l^n}E(K_y)$. 

From part (ii) of the present Proposition 4.5.1, we have 
$$\gamma_n(c-z)_{y}=\partial_n((P_{c-z} \ ({\rm mod } \ l^n))_y)$$
where $z$ is coprime to the support of $c-z$.
Hence we obtain 
$${\rm ord}\big( (P_{c-z} \ ({\rm mod } \ l^n))_y\big)
={\rm ord}(\gamma_n(c-z)_{y}).$$

From property (3) of Proposition 4.4.2, we have that the composition
 $\chi_z$ with the homomorphism, obtained from
the inclusion of group schemes $E_{l^n}\subset E$,
$$\theta :H^1(K_{y}, E_{l^n})\to H^1(K_{y},E)_{l^n}$$is an isomorphism
$$ {}_{l^n}E(K_{y})\cong H^1(K_{y}, E)_{l^n}.$$
Hence the image in $ H^1(K_{y}, E)_{l^n}$ under $\theta$ of $\chi_z(P_{c-z} \ ({\rm mod}\ l^n))$
is the cohomology class $\delta_n(c)_y$ because this is the image under 
$\theta$ of 
$\gamma_n(c)_y\in H^1(K_y,E_{l^n})$ and furthermore,  because $\theta\circ\chi_z$ is an 
isomorphism,
 $\delta_n(c)_y$ has the same order as
${\rm ord}\big( (P_{c-z} \ ({\rm mod } \ l^n))_y\big)$, which completes
the  proof of the proposition.
${\sqcap \!\!\!\!\sqcup}$

\vskip1.8in
 
\noindent {\bf 4.6. The Cassels pairing with a class $\delta_n(c)$}

\vskip0.2in

The notation and hypotheses of (4.1.1) hold in this section. The torsion abelian group
${\mathds Z}/n{\mathds Z}$ for $n\not=0$ is considered to be a subgroup of ${\mathds Q}/{\mathds Z}$ via the map $1\mapsto 1/n$. Let 

\vskip0.2in
$\coprod\!\!\!\coprod(E/K)$ be the Tate-Shafarevich group of the elliptic curve 
$E\times_FK$ over $K$;

$< -, ->_{\rm Cassels}$ be the Cassels pairing on $\coprod\!\!\!\coprod(E/K)$ (see \S2.3); 

$[ , ]_{w}: H^1(K_w,E)_n\times {}_nE(K_w)\to {\mathds Z}/n{\mathds Z}$ be the local Tate pairing 
for any place 

\qquad $w\in \Sigma_K$ of $K$ and any integer $n$ coprime to the characteristic of $F$.

\vskip0.2in

\noindent{\bf 4.6.1. Proposition.} \advance\count2 by 1 {\sl  Let $m$ and $n$ be integers $\geq 1$ and let $c\in \Lambda(m+n)$. Suppose that  $\delta_m(c)$ belongs to the Tate-Shafarevich group $\coprod\!\!\!\coprod$$(E/K)_{l^\infty}$. 
 Suppose that the element $s\in\coprod\!\!\!\coprod$$(E/K)_{l^\infty}$ has order at most 
$l^{n}$. Lift the element  $s\in\coprod\!\!\!\coprod$$(E/K)_{l^\infty}$
 to an element $S\in H^1(K,E_{l^{n}})$ and select
points 
$x(w)\in E(K_{w})$  such that $S_{w}=\partial_{l^{n}}(x(w))$
for all $w\in\Sigma_K$. Then we have
$$< \delta_m(c), s>_{\rm Cassels}
= \sum_{y\in\Sigma_K {\rm \ divides\ }  {\rm Supp}( c)}
[ \delta_{m+n}(c)_{y}, x(y)]_{y}
\leqno{(4.6.2)}$$ 
where the sum runs over the places of $K$ which divide an element of 
{\rm Supp}$(c)$.}
\advance\count2 by 1

\vskip0.2in
\noindent {\it Proof.} The construction of the Cassel's pairing is given in \S2.3.
As in the statement of the proposition,
we may lift the element  $s\in\coprod\!\!\!\coprod$$(E/K)_{l^\infty}$
of order $l^{n}$ to an element $S\in H^1(K,E_{l^{n}})$. The 
points 
$x(w)\in E(K_{w})$ are then selected such that $S_{w}=\partial_{l^{n}}(x(w))$
for all $w\in\Sigma_K$.

  By Lemma 4.2.1(iii) we have $$l^{n} \delta_{m+n}(c)=\delta_m(c).$$ From the formula (2.3.5) for the Cassels pairing, we then obtain 
$$< \delta_m(c), s>_{\rm Cassels}
= \sum_{w\in \Sigma_K}
[ \delta_{m+n}(c)_{w}, x(w)]_{w}.
\leqno{(4.6.3)}$$\advance\count2 by 1
where the sum runs over all places $w$ of $K$.

By Proposition 4.5.1(i), we have 
$\delta_{m+n}(c)_w=0$ for all places $w\in \Sigma_K $ which do not divide an element 
of Supp$(c)$.  Hence there is no contribution to the Cassels pairing in the 
sum (4.6.3) when $w$ does not divide an element of Supp$(c)$. Hence
we obtain  the formula (4.6.2). ${\sqcap \!\!\!\!\sqcup}$

\vskip0.2in
\noindent{\bf 4.6.4. Proposition.}  {\sl  Let $m$ and $n$ be integers $\geq 1$ and let $c\in \Lambda(m+n)$, $d\in \Lambda(n)$. Suppose that  $\delta_m(c)$ and  $\delta_{n}(d)$ belong to the Tate-Shafarevich group $\coprod\!\!\!\coprod$$(E/K)_{l^\infty}$. 
 Assume that the prime 
number $l\in {\cal P}$ is coprime to the order of ${\rm Pic}(A)$. Then we have
$$< \delta_m(c), \delta_{n}(d)>_{\rm Cassels}
= \sum_{y\in \Sigma_K {\rm \ divides \ }{\rm Supp}( c)
\setminus {\rm Supp}(d)}
[\delta_{m+n}(c)_{y}, (P_{d} \ ({\rm mod}\ l^n))_{{y}}]_{y}\leqno{(4.6.5)}$$ 
 where the sum runs over the places of $K$ which divide an element of 
${\rm Supp}(c)\setminus {\rm Supp}(d)$.}

\vskip0.2in
\noindent {\it Proof.}  The element $\gamma_{n}(d) \in H^1(K,E_{l^{n}})$ is a 
lifting  of 
$\delta_{n}(d)\in$ $\coprod\!\!\!\coprod$$(E/K)$ which has order at most
$l^{n}$. Suppose that  $z\in$ Supp$( c)$ and $y$ is the unique  place of $K$ lying over the place $z$ of $F$.  
If $z$ also satisfies $z\in $ Supp$( d)$ then by  Proposition 4.5.1(iv),
we have
$$\gamma_{n}(d)_{y}=0$$ because $\delta_{n}(d)_{y}=0$ as we have $\delta_{n}(d)\in $ 
$\coprod\!\!\!\coprod$$(E/K)$. If on the other hand $z$ satisfies $z\not\in $ Supp$( d)$ then by Lemma 
4.5.1(ii), 
we have $$\gamma_{n}(d)_{y}=\partial_{l^{n}}((P_{d}\ ({\rm mod}\ l^n))_{y}).$$ 
where $(P_{d}\ ({\rm mod}\ l^n))_{y}$ denotes the localization at 
$y$ of $P_d$ (mod $l^n$). The formula  (4.6.5) to be proved now follows from the 
formula (4.6.2) of Proposition 4.6.1. ${\sqcap \!\!\!\!\sqcup}$

\vfil\eject
\noindent {\bf Chapter 5. Construction of cohomology classes and proofs of the
main theorems}

\vskip0.4in\noindent {\bf 5.1. $M_r$ is finite for some $r$ }
\vskip0.2in

\vskip0.2in\noindent(5.1.1) Throughout this chapter 5, the notation of (4.1.1) of the
previous chapter 4 remains valid and  the elliptic curve $E/F$ and quadratic 
field extension $K/F$ satisfy the hypotheses (a), (b), (c) of (4.1.1). 

We  further let in this chapter  

\vskip0.2in

$N^+,N^-$ be the eigenspaces under the action of the involution 
$\tau$ whenever  $N$ 

\qquad  is a 
${\mathds Z}[{\rm Gal}(K/F)]$-module on which multiplication by $2$ is an isomorphism;

$\epsilon=\pm1$ is the sign of the functional equation of the $L$-function $L(E/{F},s)$  

\qquad of $E/F$;

$\nu(r)=(-1)^r\epsilon$ for any natural number $r\in \mathds N$;

$\nu(c)=(-1)^r\epsilon$ for any divisor $c$ of $F$ with exactly $r$ distinct prime divisors in 

\qquad its support;

$P_c\in E(K[c])$ be the Drinfeld-Heegner points of paragraph (3.4.9) for

\qquad  all
$c\in \Lambda(1)$, where $P_0\in E(K)$;

$M_r$, for all $r\in{\mathds N}$, be the quantities in ${\mathds N}\cup\infty$
given in Definition 4.1.6;

$\coprod\!\!\!\coprod$$(E/K)$ be the Tate-Shafarevich group of the elliptic curve
$E\times_FK$ over

\qquad  $K$ (see \S2.2);

Sel$_n(E/K)$ be the $n$-Selmer group of $E\times_FK$ over $K$ for any integer
$n$ coprime

\qquad  to the characteristic of $F$ (see \S2.2));

$
{\displaystyle {\rm Sel}_{a^\infty}(E/K)=\lim_{\longrightarrow\atop n} {\rm Sel}_{a^n}(E/K)}$ for any number $a$
coprime to the characteristic 

\qquad of $F$ (see also (3.1.6)).
\vskip0.2in

In this section  \S5.1 some consequences are presented  of the 
hypothesis that $M_r$ be finite for some $r$. 

\vskip0.2in
\noindent{\bf 5.1.2. Lemma.} {\sl   The Drinfeld-Heegner $P_0$ has infinite order in $E(K)$ if and only if 
  $M_0$ is finite. If $P_0$ has infinite order then   we have 
$$l^{M_0}=\vert  (E(K)/{\mathds Z}P_0)_{l^\infty} \vert $$
$${\rm ord}\ \gamma_{M_0+m}(0)=  l^m, {\sl \ \ for \ all\ } m\geq 0,$$
$$\gamma_{M_0+m}(0)\in {\rm Sel}_{l^\infty}(E/K)^{-\epsilon}
 {\sl \ \ for \ all\ } m\geq 0.$$}

\vskip0.2in
\noindent {\it Proof of Lemma 5.1.2.} The definition of $M_0$ is that (Definition 4.1.6)
$$M_0={\rm ord}_l( P_0)=\max\{ m\vert\  P_0\in l^m E(K[0])\}.$$
The group  $E(K[0])$ has no $l$-torsion (as $l\in {\cal P}$, see Proposition 1.10.1 and 
Definition 3.1.3(f)) and is a finitely generated group by the Mordell-Weil 
theorem. Hence $P_0$ has infinite order in $E(K)$ if and only if 
$M_0$ is finite.

We have 
$\gamma_{M_0+m}(0)\in {\rm Sel}_{l^\infty}(E/K)^{-\epsilon}
$, for  all $ m\geq 1,$ by Lemma 4.2.2 and because 
$\gamma_{M_0+m}(0)=\partial_{l^{M_0+m}}(P_0)$. 

Assume now that $P_0$ has infinite order. 
By [1, Chap. 7, Theorem 7.6.5] the point  $P_0$ generates a subgroup of $E(K)$ of finite index (note that in the notation of this theorem [1, Chap. 7, Theorem 7.6.5], we have $P_0=x_0$).  

By definition we have
$$\vert  (E(K)/{\mathds Z}P_0)_{l^\infty} \vert ={\max\{ l^m\vert \ P_0\in l^mE(K)\}}.$$
where the group  $E(K[0])$ has no $l^\infty$-torsion (as $l\in {\cal P}$, see Proposition 1.10.1 and 
Definition 3.1.3(f)).

 We have the Hochschild-Serre spectral sequence
$$H^i({\rm Gal}(K[0]/K),H^j(K[0], E_{l^n}))\Rightarrow H^{i+j}(
K, E_{l^n}).$$
The short exact sequence of low degree terms of this spectral sequence is in part
$$0\to H^1({\rm Gal}(K[0]/K), E_{l^n}(K[0]))
\to H^1(K, E_{l^n})\to H^1(K[0], E_{l^n})^{{\rm Gal}
(K[0]/K)}$$
$$\to H^2({\rm Gal}(K[0]/K), E_{l^n}(K[0])).$$
The two extreme terms $H^1({\rm Gal}(K[0]/K), E_{l^n}(K[0]))$
and $H^2({\rm Gal}(K[0]/K), E_{l^n}(K[0]))$ are both zero as $E_{l^n}(K[0])$ is
zero as already noted. Hence this short exact sequence provides the 
isomorphism $$H^1(K, E_{l^n})\cong H^1(K[0], E_{l^n})^{{\rm Gal}
(K[0]/K)}\leqno({5.1.3)}$$
which is induced from the injection $E(K)\to E(K[0])$.
The short exact sequence of sheaves for the \'etale topology on Spec $K$
$$0\longrightarrow E_{l^n}\longrightarrow  E{\buildrel l^n\over \longrightarrow} E\longrightarrow 0$$
then provides  the commutative diagram of cohomology groups
$$\matrix{ 
0&\longrightarrow &E_{l^n}(K)&\longrightarrow & E(K) &{\buildrel l^n\over \longrightarrow} &E(K)&\longrightarrow & H^1(K, E_{l^n})&\ldots\cr
&&\downarrow&&\downarrow&&\downarrow&&\downarrow&&\cr
0&\longrightarrow& E_{l^n}(K[0])&\longrightarrow & E(K[0]) &{\buildrel l^n\over \longrightarrow}& E(K[0])& \longrightarrow & H^1(K[0], E_{l^n})&\ldots}$$
The isomorphism $H^1(K, E_{l^n})\cong H^1(K[0], E_{l^n})^{{\rm Gal}
(K[0]/K)}$ of (5.1.3)  together with this commutative diagram then shows that the homomorphism
$$E(K)/l^nE(K)\to E(K[0])/l^nE(K[0])$$
is injective.
Hence these two numbers 
$l^{M_0}$ 
and  $\vert  (E(K)/{\mathds Z}P_0)_{l^\infty} \vert$ 
are equal.

 We obtain that  the order of $\gamma_{M_0+m}(0)$ is equal to $l^m$ from  Proposition 3.4.14(i). ${\sqcap \!\!\!\!\sqcup}$

\vskip0.2in

\noindent{\bf 5.1.4. Lemma.} {\sl  Assume that the prime 
number $l\in {\cal P}$ is coprime to the order of ${\rm Pic}(A)$. Suppose that $M_r$ is finite for some integer $r\geq 0$. Then $M_s$ is finite for all $s\geq r$ and $$M_r,M_{r+1},M_{r+2}\ldots$$ is a decreasing sequence of non-negative integers.}

\vskip0.2in
\noindent {\it Proof of Lemma 5.1.4.} 
Suppose that $M_s$ is finite for some $s\geq0$. Then there is a divisor  $$c\in \Lambda^s(M_s+1)$$ which satisfies $l^{M_s}\vert\vert P_c$ and 
$M_s<\alpha(c)$; hence the cohomology class $\gamma_{M_s+1}(c)$ is non-zero (by Proposition 3.4.14(i)). From Proposition 3.3.8, where we take   $n=M_s+1$, we obtain   a prime divisor 
$$z\in \Lambda^1(M_s+1)\leqno{(5.1.5)}$$ 
coprime to $c$ such that the localization $\gamma_{M_s+1}(c)_y$ is non-zero 
where $y$ is the unique prime of $K$ lying over $z$ (see remark 3.2.5(1)). Then by Proposition 4.5.1(ii), we have $P_c\not\in l^{M_s+1}E(K_{y})$; that is to say 
$${\rm ord}((P_c\ ({\rm mod}\ l^{M_s+1}))_y)> 1.\leqno{(5.1.6)}$$
    It follows from 
Proposition 4.5.1(iv)  that ord$\ \gamma_n (c+z)_y\ > 1$, where 
$n=M_s+1$, and in particular we have
$\gamma_n (c+z)\not=0$. Therefore by
Proposition 3.4.14(i)), or the definition of $\gamma_n(c+z)$, 
  we have $P_{c+z}\not\in l^{M_s+1} E(K[c+z])$. It follows that we have $$c+z\in \Lambda^{s+1}(M_s+1)$$ and $M_{s+1}$ is finite and we have  $M_{s+1}\leq M_s$. As $M_r$ is finite by hypothesis, we then have that $M_s$ is finite and $M_s\geq M_{s+1}$ for all $s\geq r$ by induction.
${\sqcap \!\!\!\!\sqcup}$

\vskip0.2in
\noindent{\bf 5.1.7. Lemma.} {\sl  Suppose that $M_{r-1}$ is finite where $r\geq 1$. Assume that the prime 
number $l\in {\cal P}$ is coprime to the order of ${\rm Pic}(A)$. Let $c\in \Lambda^{r}(M_{r-1})$ and put $\nu(c)=(-1)^r\epsilon$. Then we have }

\noindent (i) $\delta_{M_{r-1}}(c)\in$ $\coprod\!\!\!\coprod$$(E/K)_{l^\infty}^{-\nu(c)}$ 

\noindent (ii) {\sl $\gamma_{M_{r-1}}(c)\in{\rm Sel}_{l^{\infty}}(E/K)^{-\nu(c)}$; }

\noindent (iii) {\sl the order of $\delta_{M_{r-1}}(c)$ is at most $l^{M_{r-1}-M_{r}}$.}

\vskip0.2in
\noindent {\it Proof of Lemma 5.1.7.} By the previous Lemma 5.1.4, we have
that $M_r$ is finite and $M_{r-1}\geq M_r$. The set $\Lambda^{r}(M_{r-1})$ is non-empty by corollary 3.3.9; therefore there is an element  $c\in \Lambda^{r}(M_{r-1})$.  It follows from the definition of the $M_i$ that $l^{M_{r}}\vert P_c$ and $l^{M_{r-1}}\vert P_{c-z}$ for all prime 
divisors $z$  in the support of $c$. From  Proposition 4.5.1(iii) and 
Lemma 4.2.2 we obtain that $$\delta_{M_{r-1}}(c)\in\ \coprod\!\!\!\!\coprod(E/K)^{-\nu(c)}$$ Hence $\gamma_
{M_{r-1}}(c)$ belongs to the Selmer group ${\rm Sel}_{l^{\infty}}(E/K)^{-\nu(c)}$
(see Lemma 4.2.2).
  The order of the element $\delta_{M_{r-1}}(c)$ is at most $l^{M_{r-1}-M_{r}}$ by Proposition 3.4.14(ii). ${
\sqcap \!\!\!\!\sqcup}$

\vskip0.2in
\noindent{\bf 5.1.8. Lemma.} {\sl  Let $z\in \Lambda^1(n)$  be a prime divisor. Let $y$ be the unique prime of $K$ lying above $z$. }

\noindent (i) {\sl The Tate pairing on $H^1(K_y,E_{l^n})$ of Theorem 2.1.4  induces a non-degenerate pairing, where $\delta=\pm1$,
$$({}_{l^n}E(K_{y}))^\delta \times H^1(K_{y},E)_{l^n}^{\delta }\longrightarrow   {\mathds Z}/l^n{\mathds Z}$$
of eigenspaces under the action of $\tau$ which are finite cyclic groups of order $l^n$.  }

\noindent (ii) {\sl The  image of the homomorphism 
 $\chi_z: {}_{l^n}E(K_y)\to H^1(K_y, E_{l^n})$ is an isotropic subgroup for the alternating Tate pairing on $H^1(K_{y}, E_{l^n})$ and ${\rm Im}(\chi_z)^\delta\cong 
{\mathds Z}/l^n{\mathds Z}$ for $\delta=\pm1$. }

\vskip0.2in

\noindent {\it Proof of Lemma 5.1.8.} (i) We have isomorphisms of 
Gal$(K/F)$-modules
$${}_{l^n}E(K_y)\cong  E(K_y)_{l^n},\ \ H^1(K_y, E)_{l^n}\cong 
{\rm Hom}(\mu_{l^n}(K_{l^n}), E(K_y)).$$
The first isomorphism here follows from Proposition 4.4.1. As 
$\vert\kappa(z)\vert+1\equiv  0$ (mod $l^n$)  and $\kappa(y)^*$ has 
$\vert\kappa(z)\vert^2-1$ elements we have that
$\mu_{l^n}(K_y)$ has $l^n$ elements and is contained in the $-1$ eigenspace under $\tau$ of 
$\kappa(y)^*$. Hence we have for $\delta=\pm1$
$${}_{l^n}E(K_y)^\delta\cong H^1(K_y,E)_{l^n}^{-\delta}$$
and these groups are cyclic of order $l^n$ by Lemma 3.2.7(ii).

\noindent (ii) The map $\chi_z$ is defined in Proposition 4.4.2.  We have by property   (3) of $\chi_z$ of this Proposition 4.4.2 that the composite of 
$\chi_z$ with the surjective homomorphism 
$H^1(K_y, E_{l^n})\to H^1(K_y,E)_{l^n}$ is an isomorphism. This implies that 
 that we have for 
$\delta=\pm1$, where Im$(\chi_z)$ denotes the image of $\chi_z$ in 
$H^1(K_y, E_{l^n})$, 
$$ {\rm Im}(\chi_z)^\delta \cong H^1(K_y,E)_{l^n}^\delta.$$

Hence part (i) on Tate duality implies  that 
for $\delta=\pm1$
$$ {\rm Im}(\chi_z)^\delta\cong {\mathds Z}/l^n{\mathds Z}.$$
This evidently shows that ${\rm Im}(\chi_z)^\delta$ is an isotropic subgroup of 
$H^1(K_{y}, E_{l^n})$ for the  anti-symmetric Tate pairing. Since the cup product on $H^1(K_{y}, E_{l^n})$ is Gal$(K/{ F})$-equivariant it follows that Im$(\chi_z)
\cong {\rm Im}(\chi_z)^{+1}\oplus {\rm Im}(\chi_z)^{-1}$
 is an isotropic subgroup of 
$H^1(K_{y}, E_{l^n})$. ${\sqcap \!\!\!\!\sqcup}$

\vskip0.2in
\noindent{\bf 5.1.9. Lemma.} {\sl  Let $z\in \Lambda^1(n)$  be a prime divisor
 where
$n\geq 1$. Let $y$ be the unique prime of $K$ lying above $z$.  Let $S$ be a finite set of prime 
divisors of $\Lambda^1(n)$ not containing $z$. Let $\delta=\pm1$. Then there is a non-zero element $h\in H^1(K, E_{l^n})^\delta $ in the $\delta$-eigenspace satisfying these two conditions:
}

 (a) {\sl $h_x\in \partial_{l^n}(E(K_x))$ for any place $x$ of $K$ 
not lying over a place of $S\cup \{ z\}$; }

 (b) {\sl $h_{{x}}\in {\rm Im}( \chi_w) $ for all $w\in S$ where $x$ is the unique place of $K$ over $w$ and \nopagebreak

\qquad where $\chi_w:{}_{l^n}E(K_x)\to H^1(K_x,E_{l^n})$ is the homomorphism \nopagebreak

\qquad of Proposition 4.4.2.}

\vskip0.2in
\noindent {\it Proof of Lemma 5.1.9.} The places of $S$ remain inert in the field extension $K/F$; for any
$u\in S$
 denote by $u^\sharp$ the corresponding place of 
$K$ over $u$.

For any place $v\in \Sigma_K$ of $K$ we have the exact sequence
$$0\ \longrightarrow  \ {}_{l^n}E(K_v)\ {\buildrel \partial_{l^n}\over \longrightarrow }\  H^1(K_v, E_{l^n})\ \longrightarrow  \ H^1(K_v, E)_{l^n}\ \longrightarrow\   0\leqno{(5.1.10)}$$
where the extremities of this sequence are in duality by Theorem 2.1.6.

Let  $\delta=\pm1$.  The elliptic curve $E/F$ has good reduction at all places of $S\cup\{z\}$ (remark 3.2.5). 
Let 
$U$ be the finite subset of $\Sigma_F$ of places of $F$ given by $$U=S\cup\{z\}\cup\{{\rm  bad \ reduction \ places \ of \ }E/F {\rm 
\ in \ } \Sigma_F\}.$$ Let $U_K$ be the finite set of all places of $K$ over the places of $U$.
For all $v\in U_K$ dividing an element $u$ of $U\setminus\{ z\}$
 put
$$H_{v}={\rm Im}(\chi_u)^\delta{\rm \ \ if } \ u\in S$$
$$H_v=\partial_{l^n}(E(K_v))^\delta{\rm \ \ if \ } u\in U\setminus (S\cup\{z\}).$$
Then we have
$$\vert H_{v}\vert^2=\vert H^1(K_v, E_{l^n})^\delta \vert {\rm \ \ for \ all \ } v{\rm \ dividing \ a\ place
\ } u \in U\setminus \{z\}. $$
For the case where $u\in S$ this follows from Lemma 5.1.8(ii). For the case where $u\in U\setminus (S\cup\{z\})$, that is to say a place
of bad reduction, this follows from
the exact sequence (5.1.10)
and that the extremities of this sequence are in duality.

 From Tate local and global duality, there is a self dual exact sequence (see [9, Chap. 1, Theorem
 4.10, p.70])
$$H^1(K_U/K,E_{l^n})\to \bigoplus_{v\in U_K} H^1(K_v,E_{l^n})\to H^1(K_U/K, E_{l^n})^*$$where $N^*$, for a ${\mathds Z}/l^n{\mathds Z}$-module $N$,  denotes  ${\rm Hom}(N, {\mathds Z}/l^n{\mathds Z})$, and where $K_U$ is the maximal extension of $K$ unramified outside $U_K$. Hence the image $I$ of $H^1(K_U/K, E_{l^n})$ in $\bigoplus_{v\in U_K} H^1(K_v,E_{l^n})$ is a maximal isotropic subgroup, for the Tate pairing, of the group
$$\bigoplus_{v\in U_K} H^1(K_v,E_{l^n}).$$
Since at the place $y$ lying over $z$ we have $H^1(K_y,E_{l^n})^\delta \not=0 $ by Lemma 5.1.8(i), the subgroup $I^\delta$ is of order strictly larger than that of 
$$\bigoplus_{u\in S} {H^1(K_{u^\sharp},E_{l^n})^\delta\over H_{u^\sharp}}.$$
 Hence the natural homomorphism 
 $$H^1(K_U/K,E_{l^n})^\delta \to \bigoplus_{u\in S} {H^1(K_{u^\sharp},E_{l^n})^\delta \over H_{u^\sharp}}$$
has non-zero kernel.
Therefore we may select a non-zero element  $h\in H^1(K_U/K, E_{l^n})^\delta$ in this  
kernel. Then  $h$ satisfies the condition (b), that is to say we have
$h_{u^\sharp}\in H_{u^\sharp}$ for all $u\in S$. Furthermore, $h$ satisfies the condition (a) by the selection
of $H_v$ if $v\in\Sigma_K$ is a bad reduction place of $E$ and  because we have
$H^1(K_v^{\rm un}/K_v, E_{l^n})=\partial_{l^n}(E(K_v))$ if  $v\in \Sigma_K$  is a good reduction place of $E$,  
where 
$K_v^{\rm un}$ is the maximal unramified separable extension of $K_v$, ${\sqcap \!\!\!\!\sqcup}$

\vskip0.2in

\noindent {\it 5.1.10. Remark.} The Lemma 5.1.9 is a technical result required 
for the proof of Proposition 5.2.1 in the next section.

\vskip0.4in
\noindent {\bf 5.2. A class $\gamma_n(c)$ in the Selmer group}

\vskip0.2in
The notation and hypotheses (5.1.1) of section 5.1 hold also for this section.

\vskip0.2in
\noindent{\bf 5.2.1. Proposition.} {\sl
Let $r\geq 1$ be an integer and put $\nu(r)=(-1)^r\epsilon$. Let $G$ be a subgroup of the Selmer eigenspace ${\rm Sel}_{l^\infty}(E/K)^{-\nu(r)}$ such that  $${\rm rank}( G)\leq r.$$ 
 Assume that 
 $M_{r-1}$ is finite and that $l\in {\cal P}$ is coprime to the order of {\rm Pic}$(A)$.
Then there is a cohomology class $\gamma_{M_{r-1}}(c)$ for some divisor $c\in \Lambda^r(M_{r-1})$ such that:}

\noindent (i) {\sl $\gamma_{M_{r-1}}(c)$ belongs to the Selmer eigenspace  $ {\rm Sel}_{l^\infty}(E/K)^{-\nu(r)}$; } 

\noindent (ii) {\sl $\gamma_{M_{r-1}}(c)$ has order $l^{M_{r-1}-M_r}$;}

\noindent (iii) {\sl 
${\mathds Z}\gamma_{M_{r-1}}(c)\cap G=\{0\}.$
}
\vskip0.2in

\noindent {\it Proof of Proposition 5.2.1.} The group $ {\rm Sel}_{l^\infty}(E/K)$
is an $l^\infty$-torsion group, that is to say every element is annihilated by a power of $l$. Hence the subgroup $G$ of finite rank is therefore finite.

Let exp$(G)$ be the exponent of $G$, that is to say the largest order of an element of the abelian $l$-group $G$. Select an integer $$m\geq 1$$ such that
$$l^m\geq  \max\{ {\rm exp }( G),l^{M_{r-1}}\}$$
and put $$L=K(E_{l^m}).$$ Put
$$t={\rm rank}(G)$$
where
$${\rm rank}(G)={\rm dim}_{{\mathds Z}/l{\mathds Z}}    \ G/lG$$
and where $t\leq r$ by hypothesis. Note that $G$ is a subgroup of Hom$({\rm Gal}(L^{\rm sep}/L), E_{l^m}(L))$, as  $l^m\geq $ exp$(G)$
(see Lemma 3.3.2), and the elements of $\Lambda^1(m)$ are unramified in $L/{ F}$.
As in \S1.2, for any divisor $c$ on $F$, 
${\rm Supp}(c)$ denotes the set of distinct prime divisors in the support of $c$. 

For any divisor $c$ in $\Lambda^r(1)$, 
put
$$\Xi(c)={\rm Supp}(c)\cap\Lambda^1(m)$$
that is to say $\Xi(c)$ is the set of prime divisors in the support of 
$c$ which belong to $\Lambda^1(m)$; in particular, $\Xi(c)$ depends on $m$.
 For any prime divisor $z$ of $F$ in $\Lambda^1(1)$, select a place $z^\times$ of 
$L=K(E_{l^m})$ lying over  $z$.  Let $\Gamma(c)\subseteq {\widehat {\ G\  }}$
  be the subgroup of characters of the abelian group $G$ generated by the set of 
characters
$$\big\{\phi_{{\rm Frob}(z^\times)}\ \ \vert\ \  \ z\in \Xi(c)\big\}$$
(as in Proposition 3.3.3 and (3.3.4), applied to the finite group $G$, and Proposition 3.3.7).
Put $$s(c)  ={\rm dim}_{{\mathds Z}/l{\mathds Z}}\ \  {\Gamma(c)+l
{\widehat {\ G\  }}
\over l{\widehat {\ G\  }}}.$$
That is to say, $s(c)$ is the dimension of the image of $\Gamma(c)$ in the vector space ${\widehat {\ G\  }}/l{\widehat {\ G\  }}$ of dimension $t$; the non-negative
integer $s(c)$ is then at most equal to $t$.
Put
$$n(c)=\vert \Xi(c)\vert$$
that is to say $n(c)$ is the cardinality of $\Xi(c)$ and which is a non-negative integer at most equal to $r$, the number of prime divisors in the support of $c$.

Define the {\it defect} $\Delta(c)$ of a divisor $c\in \Lambda^r(1)$ on $F$ to be 
$$\Delta(c)=\max(t-s(c), r-n(c)).$$
Then we have $$0\leq \Delta(c)\leq r$$
and we have
$$\Delta(c)=0{\rm \ \ if \ and \ only \ if \ \ } \Gamma(c)={\widehat {\ G\  }} {\rm \ and \ }
c\in \Lambda^r(m).$$
This equivalence holds because $s(c)=t$ if and only if 
$\Gamma(c)={\widehat {\ G\  }}$ by Nakayama'a Lemma.

 We have $M_{r-1}\geq M_r$ by Lemma 5.1.4 and in particular $M_r$ is finite as $M_{r-1}$ is assumed to be finite. We may then select a divisor
$$d\in \Lambda^r(M_r+1)\leqno{(5.2.2)}$$
such that 
$$l^{M_r}\vert\vert P_d.$$
Then  the cohomology class $\gamma_{M_r+1}(d)$
is defined and belongs to $H^1(K, E_{l^{M_r+1}})^{-\nu(r)}$ and  has exact order $l$ by Lemma 4.2.2 and Proposition 3.4.14(i).
 
 Suppose that $M_r=M_{r-1}$. Then $\gamma_{M_{r-1}}(d)$ is equal to zero and evidently belongs to the Selmer eigenspace ${\rm Sel}_{l^\infty}(E/K)^{-\nu(r)}$
 which proves the lemma in this trivial case where $M_r=M_{r-1}$.
 
We may suppose for the rest of the proof of this Proposition 5.2.1  that $M_{r-1}>M_r$. 

Assume that the defect of the divisor $d$ already selected in (5.2.2)
satisfies
$$\Delta(d)>0.$$
That is to say either (where $s(d)<t$) the image of $\Gamma(d)$ is a proper subspace of ${\widehat {\ G\  }}/l{\widehat {\ G\  }}$ or  (where $n(d)<r$) that $d\not\in \Lambda^r(m)$. These two possibilities $s(d)<t$ and $n(d)<r$ for the 
divisor $d$  are not
mutually exclusive. 

We 
select a character $\psi\in 
{\widehat {\ G\  }}$  and a prime divisor $z_0\in$ Supp$(d)$ by the following recipe.

\vskip0.2in\noindent {\it Selection of the character $\psi$}

If $s(d)<t$ select a character $\psi\in 
{\widehat {\ G\  }}$ 
such that 
$$ \psi(\gamma_{M_r+1}(d))\not=0, {\rm \ \ if \ } \gamma_{M_r+1}(d)\in G\leqno{(5.2.3)}$$
and
$$\psi\in {\widehat {\ G\  }}\setminus (\Gamma(d)+l{\widehat {\ G\  }})\leqno{(5.2.4)}$$
where the condition (5.2.3) is vacuous if $\gamma_{M_r+1}(d)\not\in G$.
A character $\psi\in {\widehat {\ G\  }}$ exists which satisfies the two conditions (5.2.3) and (5.2.4) because a finite group cannot be the union of  two proper subgroups.

If $s(d)=t$ and $n(d)<r$ select any character $\psi\in {\widehat {\ G\  }}$
such that 
$$ \psi(\gamma_{M_r+1}(d))\not=0, {\rm \ \ if \ } \gamma_{M_r+1}(d)\in G\leqno{(5.2.5)}$$
where the condition (5.2.5) is vacuous if $\gamma_{M_r+1}(d)\not\in G$.
\vskip0.6in\noindent {\it Selection of the divisor $z_0\in {\rm Supp}(d)$}

If $s(d)<t$ select a prime divisor $z_0\in $ Supp$(d)$ such that
$$\big\{\phi_{{\rm Frob}(z^\times)}\ \ \vert\ \  \ z\in \Xi(d-z_0)\big\}$$
is a generating set for $\Gamma(d)$ modulo $l{\widehat {\ G\  }}$.
This choice is possible because $t\leq r$ and Supp$(d)$ has $r$ distinct 
elements.

If $s(d)=t$ and $n(d)<r$ select a prime divisor $z_0\in$ Supp$(d)$ such that 
$z_0\not\in \Lambda^1(m)$.

\vskip0.2in
 We have now defined the pair $\psi,z_0$ for the divisor $d$ where the defect
 $\Delta(d)>0$.
By Lemma 5.1.9  we may  select a cohomology class $h$ in the 
$\nu(r)$-eigenspace
$$h\in H^1(K, E_l)^{\nu(r)}$$
of the group scheme $E_l$ such that
\vskip0.2in
\noindent(5.2.6) $h\not=0$;

\noindent(5.2.7) $h_v\in \partial_l(E(K_v)) $  for all places  $v$ of $K$  coprime to ${\rm Supp}(d)$;

\noindent(5.2.8) $h_{y}\in $ Im$(\chi_z)$ for all $z\in {\rm Supp}(d-z_0)$ where $y$ is the prime of $K$ lying

\qquad  over $z$ and $\chi_z$ is the homomorphism of \S4.4.
\vskip0.2in
\noindent Note that $H^1(K, E_l)^{\nu(r)}$ is a subgroup of 
$H^1(K, E_{l^m})$ by Lemma 3.1.5.

Since $h$ is in a different eigenspace from $G$ and $ \gamma_{M_r+1}(d)$, which both belong to the $-\nu(r)$-eigenspace,  we have 
$$(G+{\mathds Z}\gamma_{M_r+1}(d))\cap {\mathds Z}h=0$$
where $G+{\mathds Z}\gamma_{M_r+1}(d)$ is the subgroup of $H^1(K, E_{l^m})^{-\nu(r)}$ generated by $G$ and $\gamma_{M_r+1}(d)$ and where  ${\mathds Z}h$ is the subgroup of $H^1(K, E_{l^m})^{\nu(r)}$ generated by $h$.

Let $D$ be the subgroup of $H^1(K, E_{l^m})$ generated by 
$G$, $\gamma_{M_r+1}(d)$ and $h$
$$D=G+{\mathds Z} \gamma_{M_r+1}(d)+{\mathds Z}h.$$
Then $D$ is a finite subgroup of $H^1(K,E_{l^m})$ by Lemma 3.1.5.
As $D$ is isomorphic to the direct product of
$G+{\mathds Z}\gamma_{M_r+1}(d)$ and ${\mathds Z}h$, 
we can select a homomorphism
$$\chi: 
D\to {\mathds Q}/{\mathds Z}$$
such that (by (5.2.3), (5.2.5), and (5.2.6))
\vskip0.2in

\noindent(5.2.9) $\chi\vert_G=\psi$;

\noindent(5.2.10) $ \chi(\gamma_{M_r+1}(d))\not=0$;

\noindent(5.2.11) $\chi(h)\not=0$.

\vskip0.2in
 By Proposition 3.3.7 applied to the finite group $D$, there is a prime divisor $$z_1\in \Lambda^1(m)\leqno{(5.2.12)}$$ distinct  from the elements of Supp$(d)$ such that, 
in the notation of (3.3.4), 
$$\chi=\phi_{{\rm Frob}(z^\times_1)}\leqno{(5.2.13)}$$
where $z^\times_1$ is a prime divisor of $L=K(E_{l^m})$ above $z_1$. We then 
obtain the cohomology class $\gamma_{M_r+1}(d+z_1)$ associated to the 
divisor $d+z_1$.

For all places $v$ of $K$, denote  the alternating cup product  induced by the 
Weil pairing  of $h$ with $\gamma_{M_r+1}(d+z_1)$ localized at $v$  by (see Theorem 2.1.4)
$$<\gamma_{M_r+1}(d+z_1)_v, h_v>_v$$
which is an element of ${\mathds Q}/{\mathds Z}$.  
The sum of local pairings over all places of $K$
$$\sum_{{\rm all\ places \ } v {\rm \ of \ } K}<\gamma_{M_r+1}(d+z_1)_v, h_v>_v=0\leqno{(5.2.14)}$$
is zero by Proposition 2.1.8.

If $v$ does not divide any element of $ {\rm Supp}(d+z_1)$  then $$\gamma_{M_r+1}(d+z_1)_v\in \partial_{l^{M_r+1}}(E(K_v))$$ by Proposition 4.5.1(i) and $$h_v\in \partial_l(E(K_v))$$ by (5.2.7), that is to say both localized elements
$\gamma_{M_r+1}(d+z_1)_v$, $h_v$ are in the image of the map
$$\partial_{l^m}: E(K_v)\to H^1(K_v,E_{l^m}).$$But the 
image of this map $\partial_{l^m}$ in $H^1(K_v,E_{l^m})$  is an isotropic subgroup for the alternating pairing $<,>_v$ (see [1, p. 403, Theorem 7.15.6] or Theorem 2.1.6(i)). Therefore we have 
$$<\gamma_{M_r+1}(d+z_1)_v, h_v>_v=0{\rm \ \ for \ all \ } v{\rm \ coprime \ to \ Supp}(d+z_1).$$

If $y$ is a place of $K$ which  divides an element $z$ of $ {\rm Supp}(d-z_0)$, then we have $h_y\in $ Im$(\chi_z)$ by (5.2.8) and 
$$\gamma_{M_r+1}(d+z_1)_y\in {\rm Im}(\chi_z)$$
by property (4) of  Proposition 4.4.2, which is the main property of the
map $\chi_z$. As Im$(\chi_z)$ is isotropic for the cup product $<,>_v$ by Lemma 
5.1.8(ii), we have $$<\gamma_{M_r+1}(d+z_1)_y, h_y>_y=0 {\rm \ \ for \ all \ } y{\rm 
\ dividing\
an\ element \ of \  Supp}(d-z_0).$$

Therefore the only possible non-zero terms in the sum $\sum_v<\gamma_{M_r+1}(d+z_1)_{v}, h_{v}>_{v}$ of  (5.2.14) are for the places of $K$ lying over the places $z_0$ and $z_1$. For these places
$z_i\in \Lambda^1$ of $F$, for $i=0,1$, denote by $y_i$ the corresponding place of $K$ lying over the place $z_i$ which remains inert in $K/F$.

From (5.2.10) and (3.3.5) and that $\chi=\phi_{{\rm Frob}(z_1^\times)}$ by (5.2.13), we have that this localization $\gamma_{M_r+1}(d)_{y_1}$
at $y_1$ is non-zero; hence by  Proposition 4.5.1(iv)
the localization
$$\delta_{M_r+1}(d+z_1)_{y_1}\in H^1(K_{y_1}, E)_{l^{M_r+1}}^{\nu(r)}$$
is non-zero. Furthermore we have 
$$h_{y_1}\in\partial_l (E(K_{y_1}))^{\nu(r)}
$$ by (5.2.7) and $h_{y_1}$  is non-zero by (5.2.11) and (5.2.13). Hence we have by Lemma 5.1.8(i)
$$<\gamma_{M_r+1}(d+z_1)_{y_1}, h_{y_1}>_{y_1}\not=0.$$
Since the sum $\sum_v<\gamma_{M_r+1}(d+z_1)_{v}, h_{v}>_{v}$ of (5.2.14) is zero, this implies 
that the local term at $y_0$
$$<\gamma_{M_r+1}(d+z_1)_{y_0}, h_{y_0}>_{y_0}$$
is non-zero.
Hence we have $\gamma_{M_r+1}(d+z_1)_{y_0}\not=0$ and so by  Proposition 
4.5.1(iv)
$$P_{d+z_1-z_0}\not\in l^{M_r+1}E(K_{y_0})$$ hence we  obtain, where
$z_1\in \Lambda^1(m)$,
$$l^{M_r}\vert\vert P_{d+z_1-z_0}\leqno{(5.2.15)}$$ because  we have 
$l^{M_r}\vert P_{d+z_1-z_0}$ and $d+z_1-z_0\in \Lambda^r(M_r+1)$ as well as
$d\in \Lambda^r(M_r+1)$ and $z_1\in \Lambda^1(m)$.

Put
 $$c_1=d+z_1-z_0  \in \Lambda^r(M_r+1).$$
On the one hand, we have that if $s(d)<t$ then the group $\Gamma(c_1)$ is generated by $\Gamma(d)$ and $\psi$
 $$ \Gamma(c_1)=\Gamma(d)+{\mathds Z}\psi \leqno{(5.2.16)}$$ because $\chi=\phi_{{\rm Frob}(z^\times_1)}$
by (5.2.13)  and $\chi\vert_G=\psi$ by (5.2.9) in this case where $s(d)<t$.  The condition that $z_1\in 
\Lambda(m)$ is satisfied by the choice of $z_1$ in (5.2.12). We then have 
from (5.2.16) if $s(d)<t$,  as $\psi\not\in \Gamma(d)+l{\widehat {\ G\  }}$ by (5.2.4), 
$$s(c_1)={\rm dim}_{{\mathds Z}/l{\mathds Z}} { \Gamma(c_1)+l{\widehat {\ G\  }}\over l{\widehat {\ G\  }}}= s(d)+1.$$

On the other hand, if $s(d)=t$ and $n(d)<r$ then $n(c_1)=n(d)+1$ as $z_1\in \Lambda^1(m)$
by (5.2.12)   and $z_0\not\in \Lambda^1(m)$ by the selection of $z_0\in $ Supp$(d)$. Hence we have in both cases
that the defect of $c_1$ is given by 
$$\Delta(c_1)=\Delta(d)-1.$$

We may by this method  construct by induction a sequence of divisors  $c_1,c_2,c_3\ldots$ in $\Lambda^r(M_r+1)$ such that their defects 
are strictly decreasing
$$\Delta(d)>\Delta(c_1) >  \Delta(c_2)> \ldots$$
and where, as in  (5.2.15),
$$l^{M_r}\vert\vert P_{c_i}{\rm \ \ for \ all \ } i.\leqno{(5.2.17)}$$
This sequence must terminate in a  divisor $c$ with zero  defect 
$\Delta(c)=0$, that is to say $\Gamma(c)={\widehat {\ G\  }}$, by Nakayama's Lemma, and $c\in \Lambda^r(m)$ and where
$$l^{M_r}\vert\vert P_{c}.\leqno{(5.2.18)}$$

 The cohomology class  $\gamma_m(c)$ is therefore defined and as $\Gamma(c)={\widehat {\ G\  }}$ we have, where $z^\times$ is a place of 
 $K(E_{l^m})$ over $z$,
$$\{ g\in G\vert\  \phi_{{\rm Frob}(z^\times)} (g)=0 {\rm \ for \ all\ } z\in 
{\rm Supp}(c)\}=0.\leqno{(5.2.19)}$$
We have  from Lemma 5.1.7(ii) that 
$$\gamma_{M_{r-1}}(c)\in 
{\rm Sel}_{l^\infty}(E/K)^{-\nu(r)}.\leqno{(5.2.20)}.$$
We have  by the definition of $M_{r-1}$ that 
$\gamma_{M_{r-1}}(c-z)=0$ for any prime divisor $z$ in the support of $c$ as 
$c-z$ has $r-1$ elements in its support. But by 
Proposition 4.5.1(iv), we have for any $z\in {\rm Supp}(c)$ where
$y$ is the place of $K$ over $z$,
$${\rm ord}\ \gamma_{M_{r-1}}(c)_{y}= {\rm ord}\ \gamma_{M_{r-1}}(c-z)_{y}.$$
Hence we obtain  for any $z\in {\rm Supp}(c)$ that    
$\gamma_{M_{r-1}}(c)_{y}=0$ where $y$ is the prime of $K$ over $z$.
 We then obtain from (5.2.19), as $g_y= \phi_{{\rm Frob}(z^\times)} (g)$ where
$y$ is the place of $K$ above $z$, that 
$$G\cap {\mathds Z}\gamma_{M_{r-1}}(c)=0.$$
Since $l^{M_r}\vert\vert P_c$, by (5.2.18),  we have that $\gamma_m(c)$ has order $l^{m-M_r}$  for all $m\geq M_r$ by Proposition 3.4.14(i).  Hence $\gamma_{M_{r-1}}(c)$ has order $l^{M_{r-1}-M_r}$ and belongs to the Selmer group  Sel$_{l^\infty}(E/K)^{-\nu(r)}$,
by (5.2.20), which proves the proposition.
${\sqcap \!\!\!\!\sqcup}$


\vskip1.6in
\noindent {\bf 5.3. Proof of  Theorem 4.1.15}

\vskip0.2in\noindent (5.3.1) The notation (5.1.1) of section \S5.1 holds also for this section.

\vskip0.2in\noindent {\bf 5.3.2. Lemma.}
{\sl  Let $A$ be a finite abelian $p$-group where $p$ is a prime number and 
with invariants $I_1\geq I_2\geq \ldots \geq I_r$. Let $B$ be a subgroup of 
$A$ with invariants $I_1\geq I_2\geq \ldots\geq I_s$ where $s\leq r$. Then
there is a subgroup $C$ of $A$ such that 
$$A=B\oplus C$$
that is to say, $A$ is the direct sum of the subgroups $B,C$. The invariants of $C$ are 
$I_{s+1}\geq I_{s+1}\geq  \ldots\geq I_r$.   }

\vskip0.2in
[The proof of this result  follows  from the 
structure theorem of finite abelian groups and is omitted.] ${\sqcap \!\!\!\!\sqcup}$

\vskip0.2in
 For the proof of Theorem 4.1.15,   as $P_0$ has infinite order 
in $E(K)$, the Tate-Shafarevich group  $\coprod\!\!\!\coprod$$
(E/K)$ is finite and $P_0$ generates a subgroup of finite index
in $E(K)$ by [1, Chapter 7, Theorem 7.6.5] which is one of the 
principal results of this monograph [1].

The image of $P_0$ 
in ${}_{l^m}E(K)$ belongs
to the $-\epsilon$-eigenspace of ${}_{l^m}E(K)$ for all $m$ (by [1, Chapter 7, Lemma 7.14.11] or Lemma 4.2.2 above). It follows 
from the decomposition into eigenspaces
$$  {}_{l^m}E(K)\cong (
{}_{l^m}E(K)  )^\epsilon\oplus ({}_{l^m}E(K))^{-\epsilon}$$
 that the $\epsilon$-eigenspace $({}_{l^m}E(K))^\epsilon$ is a finite abelian
group of order bounded  independently of $m$ and hence $E(K)^\epsilon$ is a
finite abelian group. As $E(K)$ has no $l$-torsion (as 
$l\in {\cal P}$ see Proposition
1.10.1 and Definition 3.1.3(f)) it follows that 
$$({}_{l^m}E(K))^\epsilon=0{\rm \ \  \ for\ all \ } m\leqno{(5.3.3)}$$
 and that
$$({}_{l^m}E(K))^{-\epsilon}\cong {\mathds Z}/l^m{\mathds Z} {\rm \ \ 
for \ all \ }m\leqno{(5.3.4)}$$
which proves the isomorphisms (4.1.18).

The $l^m$-Selmer group  ${\rm Sel}_{l^m}(E/K)^\pm$  belongs to an exact sequence
of eigenspaces
$$0\longrightarrow \big(  {}_{l^m}{E(K)}\big)^{\pm}\longrightarrow 
{\rm Sel}_{l^m}(E/K)^\pm\longrightarrow 
{\coprod\! \!\!\!\coprod}(E/K) 
_{l^m}^\pm\longrightarrow 0.\leqno{(5.3.5)}$$
  Hence this exact sequence
 induces an isomorphism of $\epsilon$-eigenspaces
$${\rm Sel}_{l^\infty}(E/K)^\epsilon\cong 
{\coprod\! \!\!\!\coprod}(E/K)_{l^\infty}^\epsilon.\leqno{(5.3.6)}$$

The largest invariant of the
$-\epsilon$-eigenspace
${\rm Sel}_{l^m}(E/K)^{-\epsilon}$
is at most equal to $m$. But this group ${\rm Sel}_{l^m}(E/K)^{-\epsilon}$
contains the subgroup  $({}_{l^m}E(K))^{-\epsilon}\cong {\mathds Z}/l^m{\mathds Z}$
with invariant $m$.
By Lemma 5.3.2  applied to the eigenspace  ${\rm Sel}_{l^m}(E/K)^{-\epsilon}$ and the 
subgroup $({}_{l^m}E(K))^{-\epsilon}$, we have that the 
$-\epsilon$-eigencomponent of the exact sequence
 (5.3.5) splits. Hence we obtain an isomorphism
of $-\epsilon$-eigenspaces
$${\rm Sel}_{l^m}(E/K)^{-\epsilon} \cong
{}_{l^m}E(K)\oplus
 {\coprod\! \!\!\!\coprod}(E/K) 
_{l^m}^{-\epsilon} {\rm \ \ 
for \ all \ } m.\leqno{(5.3.7)}$$
The isomorphisms (5.3.6) and (5.3.7) prove the theorem.
 ${\sqcap \!\!\!\!\sqcup}$

\vskip0.4in

\noindent {\bf 5.4. Proof of Theorems 4.1.9 and  4.1.13}

\vskip0.2in

\noindent (5.4.1)  The principle of the proof is to construct inductively 
divisors $c_{k}\in \Lambda^{k}$, $k=1,2\ldots$, such that the 
cohomology
classes $\delta_{M_{k-1}}(c_{k})$ in the Tate-Shafarevich group $ {\coprod\! \!\!\coprod}(E/K) 
_{l^\infty}$  form a basis of a maximal isotropic subgroup with respect to 
the Cassels pairing and where
$\delta_{M_{k-1}}(c_{k})$ has order $l^{N_{k}}$ for all $k$. The inductive step
is provided by the next Lemma 5.4.3.

\vskip0.2in\noindent (5.4.2) The notation and hypotheses (5.1.1) of the previous section holds also for this section. We further
denote by, where $l\in \cal P$,
\vskip0.2in
$\Sigma_K$  the set of all places of the global field $K$;

$[,]_v: {}_{l^m}E(K_v)\times H^1(K_v,E)_{l^m}\to {\mathds Z}/l^m{\mathds Z}$
the non-degenerate local pairing at

\qquad the  place $v$ of $K$ induced by the cup 
product, as in 
Theorem 2.1.6.
\vskip0.2in

\noindent If $z_i\in \Lambda^1(1)$ is a prime divisor of $F$ indexed by an integer $i$, denote by $y_i$ the prime divisor of $K$ lying above the place $z_i$ of $F$ where this place $z_i$ is inert in the field extension $K/F$.

\vskip0.2in\noindent {\bf 5.4.3. Lemma.}
{\sl Assume that $l\in {\cal P}$ is coprime to the order of {\rm Pic}$(A)$. Let $s\geq 1$ be  a positive integer and let $r,t\geq0$ be  non-negative integers. Let $S$ be a  subgroup of {\rm Sel}$_{l^s}(E/K)$. Let $e\in {\rm Sel}_{l^s}(E/K)^{-\nu(r+1)}$, $\gamma_s(c)\in {\rm Sel}_{l^s}(E/K)^{-\nu(r)}$, where $c\in \Lambda^r(s+t)$,
 be elements of the Selmer group where
$$S\cap (e,\gamma_s(c))=0$$
and where $(e,\gamma_s(c))$ is the subgroup of the Selmer group generated 
by $e,\gamma_s(c)$. Suppose also that $e, \gamma_s(c)$ both have order
$l^n$ where $n\leq s$. Then there are infinitely many prime divisors $z\in \Lambda^1(s+t)$ coprime to ${\rm Supp}(c)$ such that if $y$ is the place of $K$ over $z$ then we have}

\noindent (1) $S_y=0$;

\noindent (2)
  {\sl the value in  ${\mathds Z}/l^n{\mathds Z}$  of the local pairing at $y$
  with the class $\delta_s(c+z)$
$$[\delta_s(c+z)_y,w_y ]_y$$
has order $l^n$ where $w_y\in {}_{l^n}E(K_y)$ is a point such that $e_y=\partial_{l^n }(w_y)$.}

\vskip0.2in
\noindent {\it Proof of Lemma 5.4.3.} Note that from (3.1.10), we have that Sel$_{l^n}(E/K)$ is the subgroup
of Sel$_{l^s}(E/K)$ annihilated by $l^n$ and in particular Sel$_{l^n}(E/K)$ contains $e$ and $\gamma_s(c)$.

 Let $T$ be the subgroup of the Selmer group
${\rm Sel}_{l^s}(E/K)$ generated by $S, e$ and  $\gamma_s(c)$.
Then we have an isomorphism
$$T\cong S\oplus (e,\gamma_s(c)).$$

For a fixed non-zero element $x\in T$, the set of characters $$\chi:T\to {\mathds Z}/l^s{\mathds Z}$$ such that 
$${\rm ord}(\chi(x))<{\rm ord}(x)$$
is a proper subgroup of ${\widehat T}$. This follows as the subgroups of 
${\mathds Z}/l^s{\mathds Z}$ are linearly ordered 
${\mathds Z}/l^s{\mathds Z}\supseteq l{\mathds Z}/l^s{\mathds Z}\supseteq
l^2{\mathds Z}/l^s{\mathds Z}\ldots$.
As a  group cannot be the union of two proper subgroups,  there is then a character 
$\chi:T\to {\mathds Z}/l^s{\mathds Z}$
such that
$${\rm ord}(\chi(e))={\rm ord}(e)$$
$${\rm ord}(\chi(\gamma_s(c))={\rm ord}(\gamma_s(c))$$
$$\chi(S)=0.$$

By Proposition 3.3.8 applied to the subgroup $T$ and the character $\chi$, we may select a
prime divisor $z\in \Lambda^1(s+t)$ satisfying, where $y$ is the place of 
$K$ lying over $z$,
$${\rm ord}(e_y)={\rm ord}(e)$$
$${\rm ord}(\gamma_s(c)_y)={\rm ord}(\gamma_s(c))$$
$$S_y=0$$
where the subscript $y$ denotes the restriction at $y$ of elements of the 
Selmer group Sel$_{l^s}(E/K)$. In particular, condition (1) of the lemma is 
satisfied for this $z$.

The class $\delta_s(c+z)$, associated to the divisor $c+z\in \Lambda^{r+1}(s+t)$, then belongs to   $H^1(K,E)_{l^s}^{-\nu(r+1)}$
and $e$ belongs to $H^1(K,E_{l^s})^{-\nu(r+1)}$. 

We may select a point $w_y\in ({}_{l^n}E(K_y))^{-\nu(r+1)}$  such that $e_y=\partial_{l^n}(w_y)$ where
$e\in$ Sel$_{l^n}(E/K)$ as already noted.
Then $w_y$ has order $l^n$ in ${}_{l^n}E(K_y)$ as $e_y$ has order 
$l^n$.

Furthermore
because
$${\rm ord}(\gamma_s(c)_y)={\rm ord}(\gamma_s(c))$$
and by Proposition 4.5.1(iv)
we must have that 
$${\rm ord}(\delta_s(c+z)_y)={\rm ord}(\gamma_s(c))=l^n.$$

 The class $\delta_s(c+z)_y$, of order $l^{n}$, belongs to  the subgroup  $H^1(K_y,E)_{l^n}^{-\nu(r+1)}$
 of $H^1(K_y,E)_{l^s}^{-\nu(r+1)}$
and 
$w_y$, of order $l^{n}$, belongs to $({}_{l^n}E(K_y))^{-\nu(r+1)}$. In particular, 
$\delta_s(c+z)$ and $w_y$ both belong to the $-\nu(r+1)$-eigenspaces of their respective
spaces.
Hence the local term
$$[\delta_s(c+z)_y,w_y ]_y$$
has order $l^n$.
This follows from the non-degeneracy of 
the local pairing: 
 by Lemma 5.1.8(i), if $z\in \Lambda^1(n)$ and $y\in \Sigma_K$ is the place of $K$ over $z$ then the two elements
$$f\in {}_{l^n}E(K_{y} ), \
d\in H^1(K_{y} ,E)_{l^n}$$
give via the Tate pairing an element
$$[f,d]_{y}$$
which is non-zero if they are in the same $\tau$-eigenspace and the product of their orders is $\geq l^n$. 
Therefore condition (2) is satisfied by $z$ which proves the lemma.${\sqcap \!\!\!\!\sqcup}$

\vskip0.2in
We now prove simultaneously the two Theorems
4.1.9 and 4.1.13.  

As $P_0$ has infinite order 
in $E(K)$, the Tate-Shafarevich group  $\coprod\!\!\!\coprod$$
(E/K)$ is finite and $P_0$ generates a subgroup of finite index
in $E(K)$ by [1, Chapter 7, Theorem 7.6.5]. 

By definition of the invariants $N_i$ of the finite abelian group $\coprod\!\!\!\coprod$$
(E/K)_{l^\infty}$
in \S4.1.5, there is a maximal isotropic subgroup $D$  of $\coprod\!\!\!\coprod$$(E/K)_{l^\infty}$, with respect to the non-degenerate 
anti-symmetric Cassels pairing, where
$$D=\prod_i D_i$$
 and each $D_i$ is a cyclic group of order $l^{N_i}$ and where
$$D^{\epsilon }= \prod_{ i{\rm \ odd}} D_i
$$
and 
$$D^{-\epsilon }=\prod_{i{\rm \ even}} D_i
.$$

From Theorem 4.1.15, we have the decomposition of eigenspaces
$${\rm Sel}_{l^m}(E/K)^{\pm} \cong
\big({}_{l^m}E(K) \big)^{\pm}\oplus
 {\coprod\! \!\!\!\coprod}(E/K) 
_{l^m}^\pm {\rm \ \ 
for \ all \ } m\geq 0\leqno{(5.4.4)}$$
where the projection onto the second factor  is given by the natural surjection
$$\pi_m: {\rm Sel}_{l^m}(E/K)\longrightarrow {\coprod\! \!\!\!\coprod}(E/K) 
_{l^m}.$$

Let $m$ be an integer such that $l^m$ is greater than or equal to the exponent of the finite group ${\coprod\! \!\!\coprod}(E/K) 
_{l^\infty}$, that is to say $m\geq \max_iN_i$. For each integer $i$, let $d_i\in 
\coprod\!\!\!\coprod$$(E/K)_{l^\infty}^{-\nu(i)}$ be a generator of $D_i$ and 
let $e_i$ such that 
$$ e_i\in {\rm Sel}_{l^m}(E/K)^{-\nu(i)}, \ {\rm where\ }    \pi_m(e_i)=d_i, 
\ {\rm ord}(e_i)=l^{N_i},\leqno{(5.4.5)}$$
 be the lifting of $d_i$ to the Selmer eigenspace 
${\rm Sel}_{l^m}(E/K)^{-\nu(i)}$ via the decomposition (5.4.4) such that $e_i$ has zero component in the 
subgroup ${}_{l^m} E(K)$; in particular, 
we take  $e_i$ 
to have order equal to  $l^{N_i}$ for all $i$ and to belong to the $-\nu(i)$-eigenspace 
as $d_i$ has order $l^{N_i}$. For each valuation $v$ of $K$ and each $i$, select 
$w_{i,v}\in {}_{l^{N_i}}E(K_v)$ such that the localization $e_{i,v}$ of $e_i$ at 
$v$ satisfies
$$e_{i,v}=\partial_{l^{N_i}}(w_{i,v}).\leqno{(5.4.6)}$$
Here 
$$\partial_{l^{N_i}}: {}_{l^{N_i}}E(K_v)\to H^1(K_v, E_{l^{N_i}})\leqno{(5.4.7)}$$
denotes the connecting homomorphism associated to the morphism
$l^{N_i}:E\to E$ of multiplication by $l^{N_i}$.

The cohomology class $\gamma_{M_0+N_1}(0)$ belongs to the Selmer group
Sel$_{l^m}(E/K)^{-\epsilon }$ and has order $l^{N_1}$ as $l^{M_0}\vert\vert P_0$, by Lemma 5.1.2 or Proposition 3.4.14(i),
$${\rm ord}(\gamma_{M_0+N_1}(0))=l^{N_1}.$$
The element $e_1\in$ Sel$_{l^m}(E/K)^{\epsilon }$
has the same order
$${\rm  ord}(e_1)=  l^{N_1}.$$
 Let $S$ be the subgroup of the Selmer group
${\rm Sel}_{l^m}(E/K)$ generated by $e_i$ for all $i\geq 2$.  
The element $\gamma_{M_0+N_1}(0)$ belongs to the component ${}_{l^m}E(K)$ in the decomposition
(5.4.4) of ${\rm Sel}_{l^m}(E/K) $.
Hence we have that the subgroup of ${\rm Sel}_{l^m}(E/K)$  generated by 
$S,e_1, \gamma_{M_0+N_1}(0)$ is direct sum
$$S\oplus  {\mathds Z}e_1\oplus {\mathds Z}\gamma_{M_0+N_1}(0).$$

We may now apply Lemma 5.4.3 to $S$ and the elements $e_1, \gamma_{M_0+N_1}(0)$
where we take the parameters of the lemma to be
$$c=0, r=0, s=M_0+N_1, t=0, n=N_1.\leqno{(5.4.8)}$$
There is according to the lemma a prime divisor
 $z_1\in \Lambda^1(M_0+N_1)$ 
which satisfies the following 2 conditions, where $y_1\in \Sigma_K$ is the prime 
divisor of 
$K$ lying over $z_1$ and where the subscript $y_1$ denotes localization
at $y_1$,
and where the point 
$w_{1,y_1}$  in $(E(K_{y_1})/l^{N_1}E(K_{y_1}))^{\epsilon }$ is
such that 
$\partial_{l^{N_1}} (w_{1,y_1})=e_{1,y_1}$:-
$${\rm ord} [\delta_{M_0+N_1}(z_1)_{y_1},w_{1,y_1}]_{y_1} =l^{N_1}\leqno{(5.4.9)}$$
and
$$S_{y_1}=0.\leqno{(5.4.10)}$$
Here in (5.4.9), $\delta_{M_0+N_1}(z_1)$ is the cohomology class in 
$H^1(K,E)_{l^{M_0+N_1}}^\epsilon$ associated to $z_1$.

Let
$$\delta_{M_0}(z_1)\in {\coprod\! \!\!\!\coprod}(E/K)^{\epsilon 
}_{l^\infty}$$
be the cohomology class associated to this prime divisor $z_1\in \Lambda^1(M_0+N_1)$;
that $\delta_{M_0}(z_1)$ belongs to the Tate-Shafarevich group $\coprod\!\!\!\coprod$$(E/K)^{\epsilon 
}_{l^\infty}$ follows from  Lemma 5.1.7(i).

On the one hand, $l^{N_1}$ is the maximum order of an element of $\coprod\!\!\!\coprod$$(E/K)^{\epsilon }_{l^\infty}$. From Proposition 5.2.1 
and the isomorphism of  (5.4.4),   there is an element in the 

\noindent Selmer group ${\rm Sel}_{l^\infty}(E/K)^{\epsilon }$, and hence in the Tate-Shafarevich group $\coprod\!\!\!\coprod$$(E/K)_{l^\infty}^{\epsilon }$,
of order $l^{M_0-M_1}$. It follows that  we have the inequality
$$M_0-M_1\leq N_1.\leqno{(5.4.11)}$$

 On the other hand,  the Cassels pairing gives, as $l^n
\delta_{M_0}(z_1)=\delta_{M_0-n}(z_1)$ if $n\leq M_0$,
$$<\delta_{M_0}(z_1), l^nd>_{\rm Cassels}=
<l^n\delta_{M_0}(z_1), d>_{\rm Cassels}\leqno{(5.4.12)}$$
$$=<\delta_{M_0-n}(z_1), d>_{\rm Cassels} $$
and where we have
$$l^{N_i}\delta_{M_0+N_i-n}(z_1)=\delta_{M_0-n}(z_1).$$

By the construction of the Cassels pairing as a sum of local pairings, more precisely from
Proposition 4.6.1 and equation (4.6.2), we obtain from this 
last equation (5.4.12)  that for all $0\leq n\leq N_i-1$  where 
$l^{N_i}$ is the order of 
$d_i$
$$<\delta_{M_0}(z_1), l^nd_i>_{\rm Cassels}=\sum_{v\in \Sigma_K}
\ [\delta_{M_0+N_i-n}(z_1)_v,w_{i,v}]_{v}\leqno{(5.4.13)}$$
$$=[\delta_{M_0+N_i-n}(z_1)_{y_1} ,w_{i,y_1} ]_{y_1} $$
 where $e_{i,v} =\partial_{l^{N_i}}(w_{i,v} )$ as in (5.4.6) for all $v\in \Sigma_K$, as we have that $\delta_{M_0+N_i-n}(z_1)_v=0$ for all places $v$ of $K$ not dividing $z_1$ (by Proposition 4.5.1(i)).

The term
$$[\delta_{M_0+N_i-n}(z_1)_{y_1} ,w_{i,y_1} ]_{y_1} $$
is zero if $i\geq 2$ by (5.4.10). Hence we have from (5.4.13)
$$<\delta_{M_0}(z_1), d_i>_{\rm Cassels}=0{\rm \ \ for \ } i\geq 2.\leqno{(5.4.14)}$$

Let $i=1$. From the sum formula (5.4.13) we have 
$$<\delta_{M_0}(z_1), l^nd_1>_{\rm Cassels}=[\delta_{M_0+N_1-n}(z_1)_{y_1} ,w_{1,y_1} ]_{y_1}.\leqno{(5.4.15)} $$

By (5.4.9), the term $[\delta_{M_0+N_i}(z_1)_{y_1} ,w_{1,y_1} ]_{y_1}$ 
has order $l^{N_1}$. Hence by (5.4.15), the element $<\delta_{M_0}(z_1), d_1>_{\rm Cassels}$ 
of ${\mathds Q}/{\mathds Z}$ has order $l^{N_1}$. 
Hence the character
$$d\mapsto <\delta_{M_0}(z_1), l^nd>_{\rm Cassels},\ \ \ 
\coprod\!\!\!\!\coprod(E/K)_{l^\infty}^{\epsilon }\to {\mathds Q}/{\mathds Z},$$
is non-zero for all $n$ such that $0\leq n\leq N_1-1$ and more precisely
$$<\delta_{M_0}(z_1), l^nd_1>_{\rm Cassels}$$
is non-zero for all $n$ such that $0\leq n\leq N_1-1$.

Therefore  the character
$$d\mapsto <\delta_{M_0}(z_1),d>_{\rm Cassels},\ \ \ 
\coprod\!\!\!\!\coprod(E/K)_{l^\infty}^{\epsilon }\to {\mathds Q}/{\mathds Z},$$
vanishes on $\prod_{i\geq 2}D_i$ (by (5.4.14)) and 
 its restriction to $D_1$ generates the dual ${\widehat D}_1$. Hence the element $\delta_{M_0}(z_1)$ of
$\coprod\!\!\!\coprod$$(E/K)^{\epsilon 
}_{l^\infty}$  has order at least $l^{N_1}$
as this is the order of the cyclic group $D_1$. Since $\delta_{M_0}(z_1)$ has order at most $l^{M_0-M_1}$, by definition of the cohomology class $\delta_{M_0}(z_1)$ (see Lemma 5.1.7(iii)), we obtain that 
$$N_1\leq M_0-M_1.$$

Hence we must have from the inequality (5.4.11) 
$$N_1=M_0-M_1.\leqno{(5.4.16)}$$
It follows from this equality that $\delta_{M_0}(z_1)$ has order
$l^{M_0-M_1}$  and therefore  $l^{M_1+1}$ does not divide $P_{z_1}$ and hence we  have 
$$l^{M_1}\vert\vert P_{z_1}.$$

In summary, we have shown that $\delta_{M_0}(z_1)$ has order $l^{N_1}=l^{M_0-M_1}$, 
$l^{M_1}\vert\vert P_{z_1}$, $S_{y_1}=0$, and the character 
$d\mapsto <\delta_{M_0}(z_1),d>_{\rm Cassels}$
vanishes on $\prod_{i\geq 2}D_i$ and 
 its restriction to $D_1$ generates the dual ${\widehat D}_1$.

We now proceed by induction.
Suppose there are  prime divisors $z_1,\ldots,z_k\in \Lambda^1(M_0+N_1)$ such that
$$e_{i,y_j}=0, \ {\rm \ for \ all \ \ }i\not=j,{\rm \ and \ for \ all \ } j{\rm 
\ such \ that \ } 1\leq j\leq k\leqno{(5.4.17)}$$
and 
if 
$$c_j=z_1+\ldots+ z_j$$
then
$$l^{M_j}\vert\vert P_{c_j},\ \delta_{M_{j-1}}(c_j)\in 
{\coprod\! \!\!\!\coprod}(E/K)^{-\nu(j)} 
_{l^\infty}, \ {\rm \ for \ all \ } 1\leq j\leq k\leqno{(5.4.18)}$$
and the characters
$$\chi_j: d\mapsto <\delta_{M_{j-1}}(c_j)
,d>_{\rm Cassels}, \ 1\leq j\leq k$$
 vanish  on  $\prod_{i\geq k+1}D_{i}
$ and form a triangular basis of 
the dual of $\prod_{i\leq k} D_i$  such that 
the restriction of $\chi_j$ to the cyclic subgroup 
$D_j$ is a basis for the dual ${\widehat D}_j$ for all 
$j=1,\ldots, k$.  Suppose further we have shown that
$$M_{j-1}-M_j=N_j, {\rm \ for \ } 1\leq j\leq k.\leqno{(5.4.19)}$$
and
$${\rm ord}\ \delta_{M_{j-1}}(c_j)= l^{N_j}{\rm \ for \ } 1\leq j\leq k.\leqno{(5.4.20)}$$

We have already proved  the existence of the divisor $c_1=z_1$ and these 
properties of the previous paragraph of $e_{i,y_1}$, $P_{c_1}$, $\delta_{M_0}(c_1)$, $\chi_1$, $N_1=M_0-M_1$, including
(5.4.17),  (5.4.18), (5.4.19), and  (5.4.20) for $k=1$. 
Let  $y_i\in \Sigma_K$ be the place of $K$
above $z_i$ for all $i=1,\ldots,k$.

Let $m$ be the integer already selected such that $l^m\geq $ exp$({\coprod\! \!\!\coprod}(E/K) 
_{l^\infty})$. The  order of $\delta_{M_{k-1}}(c_k)$ in $\coprod\!\!\!\coprod$$(E/K)_{l^\infty}$ is the same as its order as a character on $D$ via the non-degenerate Cassels pairing. Since $D$ is an isotropic subgroup of 
${\coprod\! \!\!\coprod}(E/K) 
_{l^\infty}$, it follows that
$${\mathds Z}\delta_{M_{k-1}}(c_k)\cap D=0.\leqno{(5.4.21)}$$
We have
$$\gamma_{M_k+N_{k+1}}(c_k)=l^{N_k-N_{k+1}}\gamma_{M_{k-1}}(c_k)$$
as $N_k=M_{k-1}-M_k$ from (5.4.19) and where $N_k-N_{k+1}\geq 0$ by the 
definition of the integers $N_i$ as the invariants of 
$\coprod\!\!\!\coprod$$(E/K)^{
}_{l^\infty}$ in decreasing order. It follows from the induction hypothesis that the cohomology class $\delta_{M_k+N_{k+1}}(c_k)$ has order 
$l^{N_{k+1}}$, by (5.4.20),  and belongs to $\coprod\!\!\!\coprod$$(E/K)^{-\nu(k) 
}_{l^\infty}$.

We have $l^{M_k}\vert\vert P_{c_k}$ by the induction hypothesis (5.4.18). Hence by Lemma 3.4.14(i),  the cohomology class $\gamma_{M_k+N_{k+1}}(c_k)$ then has the same
order as its homomorphic image $\delta_{M_k+N_{k+1}}(c_k)$ namely $l^{N_{k+1}}$.

Let $S$ be the subgroup of the Selmer group Sel$_{l^m}(E/K)$ generated
by the elements $e_i$ for all $i\not=k+1$.
Let $T$ be the subgroup of the Selmer group Sel$_{l^m}(E/K)$ generated
by 
$\gamma_{M_k+N_{k+1}}(c_k), e_{k+1}, S$.
 From the isomorphism (5.4.4) and that ${\mathds Z}\delta_{M_k+N_{k+1}}(c_k)$
has trivial intersection with $D$ by (5.4.21),    there is an equality of subgroups of the Selmer
group ${\rm Sel}_{l^m}(E/K)$, where the sums on the right hand side are
direct,
$$T= {\mathds Z} \gamma_{M_k+N_{k+1}}(c_k)\oplus {\mathds Z}e_{k+1}\oplus S.$$

We may now apply Lemma 5.4.3 to $S$ and the elements $e_{k+1}, \gamma_{M_k+N_{k+1}}(c_k)$
where we take
$$c=c_k, r=k, s=M_k+N_{k+1}, t=M_0+N_1- s, n=N_{k+1}.\leqno{(5.4.22)}$$
There is then according to the lemma a prime divisor
 $z_{k+1}\in \Lambda^1(M_0+N_1)$ 
which satisfies the following 2 conditions
$${\rm ord} [\delta_{M_k+N_{k+1}}(c_k+z_{k+1}),w_{{k+1},y_{k+1}}]_{y_{k+1}} =l^{N_{k+1}}\leqno{(5.4.23)}$$
and
$$S_{y_{k+1}}=0\leqno{(5.4.24)}$$
 where $y_{k+1}\in \Sigma_K$ is the prime 
divisor of 
$K$ lying over $z_{k+1}$ and where the subscript $y_{k+1}$ denotes localization
at $y_{k+1}$,
and where the point 
$w_{{k+1},y_{k+1}}$  in $(E(K_{y_{k+1}})/l^{N_{k+1}}E(K_{y_{k+1}}))^{\epsilon }$ is
such that 
$\partial_{l^{N_{k+1}}} (w_{{k+1},y_{k+1}})=e_{{k+1},y_{k+1}}$.

For this selection of $z_{k+1}\in \Lambda^1(M_0+N_{1})$, note that $M_k+N_{k+1}\leq M_0+N_1$  and so $t\geq 0$,
where $t$ is the parameter of (5.4.22),  because $M_r,N_r$ are both decreasing
sequences of integers and that $H^1(K, E_{l^{M_k+N_{k+1}}})$ 
is a subgroup of 
$H^1(K, E_{l^{M_0+N_{1}}})$ by Lemma 3.1.5.

Let $c_{k+1}$ be the divisor which is the sum of the $z_i$, for $i=1,\ldots, k+1$,
$$c_{k+1}=\sum_{j=1}^{k+1} z_{j}.\leqno{(5.4.25)}$$ 
Let
$$\delta_{M_k}(c_{k+1})\in {\coprod\! \!\!\!\coprod}(E/K)^{-\nu(k+1)
}_{l^\infty}$$
be the cohomology class associated to this divisor $c_{k+1}\in \Lambda^{k+1}(M_0+N_1)$;
that $\delta_{M_k}(c_{k+1})$ belongs to the Tate-Shafarevich group $\coprod\!\!\!\coprod$$(E/K)^{-\nu(k+1)
}_{l^\infty}$ follows from  Lemma 5.1.7(i).

Then  for $0\leq n\leq N_{k+1}-1$  by the construction of the Cassels pairing as a sum of local terms, more
precisely from Proposition 4.6.1 and equation (4.6.2), we have the following sum formulae as $d_i$ has order $l^{N_i}$
$$<\delta_{M_k}(c_{k+1}), l^nd_i>_{\rm Cassels}=<\delta_{M_k-n}(c_{k+1}), d_i>_{\rm Cassels}\leqno{(5.4.26)}$$
$$=\sum_{v\in \Sigma_K}[\delta_{M_k-n+N_i}(c_{k+1})_{v}, w_{i,v}]_{v}
 =\sum_{j=1}^{k+1}
[\delta_{M_k-n+N_i}(c_{k+1})_{y_j}, w_{i,y_j}]_{y_j}.$$
because $\delta_{M_k-n+N_i}(c_{k+1})_{v}=0$ for all $v\in \Sigma_K$ not dividing an element of  Supp$(c_{k+1})$
by Proposition 4.5.1(i).
Here $ w_{i,y_j}\in {}_{l^{N_i}}E(K_{y_j})$ is an element already chosen (see (5.4.6)) such that 
$$\partial_{l^{N_i}}( w_{i,y_j})=e_{i,y_j}.$$

We have the following  sum for the Cassels pairing, obtained from those of (5.4.26),
$$ <\delta_{M_k}(c_{k+1}), l^nd_i>_{\rm Cassels}=\sum_{j=1}^{k+1}
[\delta_{M_k-n+N_i}(c_{k+1})_{y_j}, w_{i,y_j}]_{y_j}.\leqno{(5.4.27)}$$
 We have $w_{i,y_j}=0$ for all $i\not=j$ and all $j$ such that $1\leq j\leq k+1$ by (5.4.17), for $j\leq k$, and by (5.4.24), for $j=k+1$.
It follows that  for $i\geq k+1$ all terms $[\delta_{M_k-n+N_i}(c_{k+1})_{y_j}, w_{i,y_j}]_{y_j}$ 
of this sum (5.4.27)  are zero except the last $[\delta_{M_k-n+N_i}(c_{k+1})_{y_{k+1}}, w_{i,y_{k+1}}]_{y_{k+1}}$ 
 and the entire sum  
 is zero for $i>k+1$.

By (5.4.23) the local term
$$[\delta_{M_k+N_{k+1}-n}(c_{k+1})_{y_{k+1}}, w_{k+1,y_{k+1}}]_{y_{k+1}}$$ is non-zero for all integers $n$ such that $0\leq n\leq N_{k+1}-1.$

Therefore the character, by (5.4.27),
$$\chi_{k+1}: \ d\mapsto <\delta_{M_k}(c_{k+1}),d>_{\rm Cassels}$$
vanishes on $\prod_{i> k+1}D_{i}$ and 
 its restriction to $D_{k+1}$ generates ${\widehat D}_{k+1}$, as $D_{k+1}$ has order $l^{N_{k+1}}$ by definition.   Hence $\chi_{k+1}$ extends the 
triangular 
basis $\chi_1,\ldots, \chi_k$ to generate $\prod_{i\leq k+1}{\widehat D}_i$
and
$\delta_{M_k}(c_{k+1})$ has order at least $l^{N_{k+1}}$. Since it has order at most $l^{M_k-M_{k+1}}$, by definition of the cohomology class
$\delta_{M_k}(c_{k+1})$ (see Lemma 5.1.7(iii)), we conclude that
$$N_{k+1}\leq M_k-M_{k+1}\leqno{(5.4.28)}$$
and also
$$l^{N_{k+1}}\leq \ {\rm ord}\  \delta_{M_k}(c_{k+1})\ \leq l^{M_k-M_{k+1}}.\leqno{(5.4.29)}$$

Let $C_k$ be the subgroup of the Selmer eigencomponent Sel$_{l^m}(E/K)^{-\nu({k+1})}$ given by
$$C_k=\bigg(\gamma_{M_0+N_1}(0),e_1,\ldots, e_k, \gamma_{M_0}(c_1),\ldots, \gamma_{M_{k-1}}(c_k)\bigg)^{-\nu({k+1})}.$$
If $k$ is even then $C_k$ is generated by $e_1, \gamma_{M_0}(c_1), e_3, \gamma_{M_2}(c_3), \ldots , e_{k-1}, \gamma_{M_{k-2}}(c_{k-1})$
of which there are $k$ in number. If $k$ is odd then $C_k$ is generated by $\gamma_{M_0+N_1}(0),e_2, \gamma_{M_1}(c_2), e_4, \gamma_{M_3}(c_4), \ldots , e_{k-1}, \gamma_{M_{k-2}}(c_{k-1})$
of which there are $k$ in number. We then have for all integers $k$
$${\rm rank}(C_k)\leq k.\leqno{(5.4.30)}$$

The elements $e_1,\ldots, e_k$ generate a subgroup of Sel$_{l^m}(E/K)$ isomorphic to $\prod_{i\leq k} D_i$ by the decomposition 
(5.4.4) of the Selmer group. Furthermore,  the elements $\gamma_{M_0}(c_1),\ldots, \gamma_{M_{k-1}}(c_k)$ generate a subgroup 
of Sel$_{l^m}(E/K)$ isomorphic to the dual of  $\prod_{i\leq k} D_i$ as $\gamma_{M_{i-1}}(c_i)$ has the same 
order as $\delta_{M_{i-1}}(c_i)$ for all $i=1,\ldots, k$ (see  (5.4.18), (5.4.19), (5.4.20),  and Proposition 3.4.14(i)).
In the decomposition (5.4.4) of Sel$_{l^m}(E/K)$, take $S$ to be the subgroup 
${}_{l^m}E(K) \oplus_{1\leq i\leq k} \prod D_i\oplus \bigoplus_{1\leq i\leq k}  {\mathds Z} \gamma_{M_{i-1}}(c_i) $.
Then we have $C_k=S^{-\nu({k+1})}$. 

If a finite abelian group $G$ is a direct product $G_1\times
G_2$ of 2 subgroups and $g\in G$ is such that ${\mathds Z}g\cap G_1=0$ then the order of 
the element $g$ is at most the exponent  of $G_2$. Then by the previous remark where we take 
$G={\rm Sel}_{l^m}(E/K)^{-\nu({k+1})} $ and $G_1=C_k$, we have that 
$l^{N_{k+1}}$ is the maximum order of an element $c\in {\rm Sel}_{l^m}(E/K)^{-\nu({k+1})} $  if $$
{\mathds Z}c\cap C_k=0$$
by the decomposition (5.4.4) of the 
Selmer group ${\rm Sel}_{l^m}(E/K)$.

On the other hand, by Proposition 5.2.1 and (5.4.30) applied to subgroup $C_k$ of the  Selmer group there is an element in ${\rm Sel}_{l^m}(E/K)^{-\nu({k+1})}$ of order $l^{M_k-M_{k+1}}$ satisfying  $
{\mathds Z}c\cap C_k=0$. Hence we have 
$$M_k-M_{k+1}\leq N_{k+1}$$
and so by (5.4.28) we have 
$$N_{k+1}=M_k-M_{k+1}.\leqno{(5.4.31)}$$
It follows from this equality  and (5.4.29) that 
$${\rm ord}\ \delta_{M_k}(c_{k+1})= l^{M_k-M_{k+1}}.\leqno{(5.4.32)}$$Therefore 
$l^{M_{k+1}+1}$ does not divide $P_{c_{k+1}}$ and $l^{M_{k+1}}\vert P_{c_{k+1}}$ and hence we have
$$l^{M_{k+1}}\vert\vert P_{c_{k+1}}.\leqno{(5.4.33)}$$
The properties (5.4.31), (5.4.32), and (5.4.33) complete the proof of the 
 induction step and this proves the Theorems 4.1.9 and 4.1.13. ${\sqcap \!\!\!\!\sqcup}$

\vskip0.4in

\noindent {\bf 5.5. Proof of Theorems 1.1.1. and 4.1.14 }

\vskip0.2in

\noindent {\it Proof of Theorem 1.1.1.} From Proposition 2.2.5, we have as $l\in 
{\cal P}$ is an odd prime number that $ 
{\coprod\!\!\!\coprod}(E/F)_{l^\infty}\cong{\coprod\!\!\!\coprod}(E/K)_{l^\infty}^{+1}
$.
The Theorem now follows immediately from Theorem 4.1.9. ${\sqcap \!\!\!\!\sqcup}$

\vskip0.2in
\noindent {\it Proof of Theorem 4.1.14.} 
As $P_0$ has
infinite 
order,  by Lemma 5.1.2 and [1, Chap. 7, Theorem 7.6.5] we have 
that  $M_0$ is finite, the group ${\mathds Z}P_0$ has finite index in $E(K)$, and the highest power of $l$ dividing the index $[E(K):{\mathds Z}P_0]$  is equal to
$$l^{M_0}={\big\vert} (E(K)/{\mathds Z}P_0)_{l^\infty}{\big\vert}.$$

By Theorem 4.1.9, we have that  
$${\bigg\vert } {\coprod\!\!\!\!\coprod}(E/K)_{l^\infty}{\bigg\vert}  =\prod_{i\geq 0}l^{2(M_i-M_{i+1})}= l^{2(M_0-M_\infty)}$$
where 
$$M_\infty=\min_{i\in {\bf N}} M_i
$$
and where this minimum exists as the $M_i$ form a decreasing
sequence of non-negative integers by Lemma 5.1.4. ${\sqcap \!\!\!\!\sqcup}$

\vskip0.4in 
\noindent {\bf 5.6. Generators of  Tate-Shafarevich groups}
\vskip0.2in

\vskip0.2in\noindent (5.6.1)  The notation (5.1.1) of the section \S5.1 holds also 
for this section. We further
denote by
\vskip0.2in
$\Sigma_K$  the set of all places of the global field $K$;

$[,]_v: {}_{l^m}E(K_v)\times H^1(K_v,E)_{l^m}\to {\mathds Z}/l^m{\mathds Z}$
the non-degenerate local pairing at

\qquad the  place $v$ of $K$ induced by the cup 
product, as in 
Theorem 2.1.6.

\vskip0.2in
\noindent{\bf 5.6.2. Theorem.} {\sl
Let $l$ be a prime number belonging to $\cal P$; assume that $l$ is coprime to the order of ${\rm Pic}(A)$.  Suppose that $P_0$ has  infinite 
order in $E(K)$ and let $a$ be an integer such that $a\geq 2M_0$. Then we have}

\noindent (1) {\sl the 
classes $ \delta_{M_0}(c)$, for all $c\in \Lambda^1(a) $, generate $\coprod\!\!\!\coprod$$(E/K)_{l^\infty}^{\epsilon }$;}

\noindent (2)  {\sl the classes 
$ \delta_{M_1}(c)$, for all $c\in \Lambda^2(a) $, generate $\coprod\!\!\!\coprod$$(E/K)_{l^\infty}^{-\epsilon }$.}

\vskip0.2in
\noindent{\bf 5.6.3. Theorem.} {\sl Under the hypotheses of  Theorem 5.6.2, 
 the 
classes $ \delta_{M_r }(c)$, for all $c\in \Lambda^r(a) $, generate $\coprod\!\!\!\coprod$$(E/F)_{l^\infty}$ where $r=(3-\epsilon)/2$.}

\vskip0.2in

\noindent {\it Proof of Theorem 5.6.3.}   This evidently      follows from Theorem 5.6.2 and  Proposition 2.2.5. ${\sqcap \!\!\!\!\sqcup}$

\vskip0.2in
\noindent {\it Proof of Theorem 5.6.2.} By Theorem 4.1.9 above or [1, Chap. 7, Theorem 7.6.5], the group $\coprod\!\!\!\coprod$$(E/K)_{l^\infty}$ is finite.  By  Lemmas 5.1.2 and  5.1.4, the quantities $M_0,M_1,M_2,\ldots $  are all finite and form a 
decreasing sequence of non-negative integers.  We fix an integer $a\geq 2M_0$. The classes $ \delta_{M_0}(z)$,
for $z\in \Lambda^1(a) $, and the classes 
$ \delta_{M_1}(z_1+z_2)$, for $z_1+z_2\in \Lambda^2(a) $, belong to 
$\coprod\!\!\!\coprod$$(E/K)_{l^\infty}$ by Lemma 5.1.7(i). We now prove separately the two parts of the 
theorem.

\vskip0.2in
\noindent (1) Suppose that $d\in $ $\coprod\!\!\!\coprod$$(E/K)^{\epsilon }_{l^\infty}$ has order exactly $l^M$ for some $M>0$ and is in the $\epsilon$-eigenspace of
$\coprod\!\!\!\coprod$$(E/K)_{l^\infty}$. 
By Theorem 4.1.15 there is an isomorphism of  $\epsilon$-components
$${\rm Sel}_{l^\infty}(E/K)^\epsilon\cong 
{\coprod\! \!\!\!\coprod}(E/K)_{l^\infty}^\epsilon\leqno{(5.6.4)}$$
induced from the natural surjection of the Selmer group onto the 
Tate-Shafarevich group.

By  Theorem 4.1.9, we have $M\leq M_0$. By the isomorphism (5.6.4), we 
may  lift $d$ to an element of the Selmer group $e\in {\rm 
Sel}_{l^\infty}(E/K)^\epsilon $ of order $l^M$. 
The cohomology class $\gamma_{M_0+M}(0)$ belongs to the $-\epsilon$-eigenspace 
$ {\rm 
Sel}_{l^\infty}(E/K)^{-\epsilon} $ of the 
Selmer group and has order $l^M$ by Lemma 5.1.2. 

We may apply Lemma 5.4.3 to the elements $e$, $\gamma_{M_0+M}(0)$ and the 
subgroup $S=0$ where we take the parameters of the lemma to be
$$ c=0, r=0, s=M_0+M, t= a-(M_0+M), n=M .\leqno{(5.6.5)}$$
Note that $t\geq 0$ as $a\geq 2M_0\geq M_0+M$.
There is according to the lemma a prime divisor
 $z\in \Lambda^1(a)$ 
which satisfies the following  condition
$${\rm ord}[ \delta_{M_0+M}(z)_y, w_y]_y=l^M\leqno{(5.6.6)}$$
 where $y\in \Sigma_K$ is the prime 
divisor of 
$K$ lying over $z$ and where the subscript $y$ denotes localization
at $y$,
and where the point   
$w_{y}$  in $(E(K_{y})/l^{M}E(K_{y}))^{\epsilon }$ is
such that 
$\partial_{l^{M}} (w_{y})=e_{y}$.
Here in (5.6.6), $\delta_{M_0+M}(z)$ is the cohomology class in $H^1(K,E)_{l^{M_0+M}}$
associated to $z$.

Let $\delta_{M_0}(z)$ be the cohomology class of 
$\coprod\!\!\!\coprod$$(E/K)^{\epsilon }_{l^\infty}$
associated to $z$, where this class belongs to the Tate-Shafarevich group
by Lemma 5.1.7.

The   Cassels pairing gives, as $l^u\delta_{M_0}(z)=\delta_{M_0-u}(z)$ if $u\leq M_0$,
$$<\delta_{M_0}(z), l^ud>_{\rm Cassels}=<\delta_{M_0-u}(z), d>_{\rm Cassels}. \leqno{(5.6.7)}$$
As $d$ has order $l^{M}$ and that we have 
$$l^{M}\delta_{M_0+M-u}(z)=\delta_{M_0-u}(z),$$
by the construction of the Cassels pairing in terms of local Tate pairings (Proposition 4.6.1 and equation (4.6.2)) we obtain from this last equation (5.6.7)
$$<\delta_{M_0}(z), l^ud>_{\rm Cassels}=[\delta_{M_0+M-u}(z)_{y} ,w_{y} ]_{y} $$
where as above  $$e_{y} =\partial_{l^M}(w_{y} ),$$ 
and where $w_y\in {}_{l^M}E(K_y)$ as we have that $\delta_{M_0+M-u}(z)_v=0$ for all places $v$ of $K$ not dividing $z$ by Proposition 4.5.1(i).

 By (5.6.6)   the element 
given by 
the Tate pairing
 $$[\delta_{M_0+M-u}(z)_{y} ,w_{y} ]_{y} $$ of ${\mathds Q}/{\mathds Z}$
is non-zero for all integers $u$ such that  $1\leq u\leq M-1$ and hence
$$<\delta_{M_0}(z), l^{u}d>_{\rm Cassels}$$
is non-zero for all integers $u$ such that  $1\leq u\leq M-1$.
We obtain  that the character, where $ z\in \Lambda^1(a)$, 
$$f\mapsto < \delta_{M_0}(z), f>_{\rm Cassels}, 
\ \ \coprod\!\!\!\!\coprod(E/K)_{l^\infty}\to {\mathds Q}/{\mathds Z},
$$
of $\coprod\!\!\!\coprod$$(E/K)_{l^\infty}$ 
when  restricted to ${\mathds Z}d$ generates the dual 
${\widehat{ {\mathds Z}d}}$ of this subgroup  ${ {\mathds Z}d}$. 

 As $d$ is any element of the abelian group $\coprod\!\!\!\coprod$$(E/K)_{l^\infty}^{\epsilon }$, the non-degeneracy of Cassels pairing implies that  the  classes $\{ \delta_{M_0}(z),z\in \Lambda^1(a)\} $ generate
the $\epsilon$-eigenspace $\coprod\!\!\!\coprod$$(E/K)_{l^\infty}^{\epsilon }$ which proves the part of the theorem for $\coprod\!\!\!\coprod$$(E/K)_{l^\infty}^{\epsilon }$.

\vskip0.2in
\noindent (2) Suppose that $f\in $ 
  $\coprod\!\!\!\coprod$$(E/K)_{l^\infty}^{-\epsilon }$ has order exactly $l^{M'}$ and is 
in the $-\epsilon$-eigenspace of $\coprod\!\!\!\coprod$$(E/K)_{l^\infty}$. 
We have $M'\leq M_0$ by Theorem 4.1.9.

We have from Theorem 4.1.15,    an isomorphism
compatible with the $\tau$-eigenspaces
$${\rm Sel}_{l^m}(E/K) \cong
{}_{l^m}E(K)\oplus
 {\coprod\! \!\!\!\coprod}(E/K) 
_{l^\infty} {\rm \ \ 
for \ all \ } m\geq N_1$$
from which as stated in this theorem we 
obtain the isomorphism, taking $m=a\geq 2M_0\geq N_1$,
$$\Delta_a :{\rm Sel}_{l^a}(E/K)^{-\epsilon} \cong
{\mathds Z}/l^a{\mathds Z} \oplus
 {\coprod\! \!\!\!\coprod}(E/K) 
_{l^\infty}^{-\epsilon} .$$
The projection
of ${\rm Sel}_{l^a}(E/K)$ onto the second factor ${\coprod\! \!\!\coprod}(E/K) 
_{l^\infty}$ of this direct sum  is the natural surjection of the Selmer group onto the 
Tate-Shafarevich group.

Lift $f$ to $
g\in {\rm Sel}_{l^a}(E/K)^{-\epsilon}$ via this isomorphism $\Delta_a$ 
for the integer $a$ and where
$g$ has order $l^{M'}$ and has zero component in the first term ${\mathds Z}/l^a{\mathds Z}$
of this decomposition of the Selmer group.

By Theorem 4.1.9 or alternatively Proposition 5.2.1, there 
is an element  $d\in $ $\coprod\!\!\!\coprod$$(E/K)_{l^\infty}^{\epsilon }$  of order exactly $
l^{M_0-M_1}$. 
Lift, via the isomorphism (5.6.4), $d$ to an element of the Selmer group
$e\in {\rm Sel}_{l^{a}}(E/K)^\epsilon$ which is also of order
$l^{M_0-M_1}$.

The cohomology class $\gamma_{2M_0-M_1}(0)$ belongs to the $-\epsilon$-eigenspace 
$ {\rm 
Sel}_{l^a}(E/K)^{-\epsilon} $ of the 
Selmer group and has order $l^{M_0-M_1}$ by Lemma 5.1.2.

Let $T$
be the subgroup of 
$ {\rm 
Sel}_{l^a}(E/K) $ generated by the three elements  $\gamma_{2M_0-M_1}(0), 
e$ and $g$. 
As the two elements $\gamma_{2M_0-M_1}(0), g$ belong to the different 
components of $ {\rm 
Sel}_{l^a}(E/K)^{-\epsilon} $ under the isomorphism $\Delta_a$ 
and as $e$ belongs to a different eigenspace 
$ {\rm 
Sel}_{l^a}(E/K)^{\epsilon} $,  the group $T$ is the direct sum
of the subgroups generated by these three elements, that is to say we have
$$T\cong {\mathds Z}\gamma_{2M_0-M_1}(0)\oplus {\mathds Z}e\oplus
{\mathds Z}g.$$

We may apply Lemma 5.4.3 to the subgroup $S={\mathds Z}g$ and the elements 
$\gamma_{2M_0-M_1}(0)$  and $e$ of the Selmer group; we take the parameters of the lemma to be
$$  c=0,r=0, s=2M_0-M_1, t=a-(2M_0-M_1), n=M_0-M_1.\leqno{(5.6.8)}$$
Note that $t\geq 0$ as $a\geq 2M_0$. 
 There is according to the lemma a prime divisor
 $z_1\in \Lambda^1(a)$ 
which satisfies the following  2 conditions
$${\rm ord }[\delta_{2M_0-M_1}(z_1)_{y_1} ,w_{y_1} ]_{y_1} =l^{M_0-M_1}\leqno{(5.6.9)}$$
$$({\mathds Z} g)_{y_1}=0\leqno{(5.6.10)}$$
 where $y_1\in \Sigma_K$ is the prime 
divisor of 
$K$ lying over $z_1$,
and where the point   
$w_{y_1}$  in $(E(K_{y_1})/l^{M_0-M_1}E(K_{y_1}))^{\epsilon }$ is
such that  $\partial_{l^{M_0-M_1}} (w_{y_1})=e_{y_1}$ .
Here in (5.6.9), $\delta_{2M_0-M_1}(z_1)$ is the cohomology class in $H^1(K,E)_{l^{2M_0-M_1}}$
associated to $z_1$.

Let
$$\delta_{M_0}(z_1)\in {\coprod\! \!\!\!\coprod}(E/K)^{\epsilon 
}_{l^\infty}$$
be the cohomology class associated to this prime divisor $z_1\in \Lambda^1(a)$;
that $\delta_{M_0}(z_1)$ belongs to $\coprod\!\!\!\coprod$$(E/K)^{\epsilon 
}_{l^\infty}$ follows from  Lemma 5.1.7(i).

 The Cassels pairing gives, as $l^u
\delta_{M_0}(z_1)=\delta_{M_0-u}(z_1)$ if $u\leq M_0$,
$$<\delta_{M_0}(z_1), l^ud>_{\rm Cassels}=
<l^u\delta_{M_0}(z_1), d>_{\rm Cassels}\leqno{(5.6.11)}$$
$$=<\delta_{M_0-u}(z_1), d>_{\rm Cassels} $$
where
$$l^{M_0-M_1}\delta_{2M_0-M_1-u}(z_1)=\delta_{M_0-u}(z_1).$$

By the construction of the Cassels pairing as a sum of local terms, more precisely from
Proposition 4.6.1 and equation (4.6.2), we obtain from this 
last equation (5.6.11)  that for all $0\leq u\leq M_0-M_1-1$  where 
$M_0-M_1$ is the order of 
$d$
$$<\delta_{M_0}(z_1), l^ud>_{\rm Cassels}=\sum_{v\in \Sigma_K}
\ [\delta_{2M_0-M_1-u}(z_1)_v,w_{v}]_{v}\leqno{(5.6.12)}$$
$$=[\delta_{2M_0-M_1-u}(z_1)_{y_1} ,w_{y_1} ]_{y_1} $$
 where $e_{v} =\partial_{l^{M_0-M_1}}(w_{v} )$ for all $v\in \Sigma_K$, and where we have that $\delta_{2M_0-M_1-u}(z_1)_v=0$ for all places $v$ of $K$ not dividing $z_1$ (by Proposition 4.5.1(i)).

From the sum formula (5.6.12) we have 
$$<\delta_{M_0}(z_1), l^ud>_{\rm Cassels}=[\delta_{2M_0-M_1-u}(z_1)_{y_1} ,w_{y_1} ]_{y_1}. $$
By (5.6.9), $[\delta_{2M_0-M_1}(z_1)_{y_1} ,w_{y_1} ]_{y_1} $ has order $l^{M_0-M_1}$.

 Hence the map
$$s\mapsto <\delta_{M_0}(z_1), l^us>_{\rm Cassels},\ \ \ 
\coprod\!\!\!\!\coprod(E/K)_{l^\infty}^{\epsilon }\to {\mathds Q}/{\mathds Z},$$
is non-zero for all $u$ such that $0\leq u\leq M_0-M_1-1$ and more precisely
$$<\delta_{M_0}(z_1), l^ud>_{\rm Cassels}$$
is non-zero for all $u$ such that $0\leq u\leq M_0-M_1-1$.

Therefore  the character
$$s\mapsto <\delta_{M_0}(z_1),s>_{\rm Cassels},\ \ \ 
\coprod\!\!\!\!\coprod(E/K)_{l^\infty}^{\epsilon }\to {\mathds Q}/{\mathds Z},$$
is such that its restriction to ${\mathds Z}d$ generates the dual ${\widehat {{\mathds Z}d}}$. Hence 
$\delta_{M_0}(z_1)$ has order at least $l^{M_0-M_1}$
as this is the order of the cyclic group ${\mathds Z}d$. Since $\delta_{M_0}(z_1)$ has order at most $l^{M_0-M_1}$
, by definition of the cohomology class $\delta_{M_0}(z_1)$ (see Lemma 5.1.7(iii)), we obtain that 
$\delta_{M_0}(z_1)$ has order given by 
$${\rm ord}\ \delta_{M_0}(z_1)= l^{M_0-M_1}.$$

We have $l^{M_1}\vert P_{z_1}$ by  definition of $M_1$. It follows that $\gamma_{M_0}(z_1)$ has 
order $\leq l^{M_0-M_1}$ by Proposition 3.4.14(i). As $\delta_{M_0}(z_1)$ has order $l^{M_0-M_1}$
by the previous paragraph and $\delta_{M_0}(z_1)$ is a homomorphic image of 
$\gamma_{M_0}(z_1)$ we must have that $\gamma_{M_0}(z_1)$ has 
order exactly  given by
$${\rm ord}\ \gamma_{M_0}(z_1)=l^{M_0-M_1}.$$ 
Therefore we have by Proposition 3.4.14(i)
 $$l^{M_1}\vert\vert P_{z_1}.$$

We may apply Lemma 5.4.3 to the subgroup $S=0$ and the two elements 
$\gamma_{M_1+M'}(z_1)$ and $g$  of the Selmer group which both have order 
$l^{M'}$; we take the parameters of the lemma to be
$$  c=z_1,r=1, s=M_1+M', t=a-(M_1+M'), n=M'.\leqno{(5.6.13)}$$
Note that $t\geq 0$ as $a\geq 2M_0$ and $M', M_1\leq M_0$. 
 There is according to the lemma a prime divisor
 $z_2\in \Lambda^1(a)$ 
which satisfies the following   condition
$${\rm ord }[\delta_{M_1+M'}(z_1+z_2)_{y_2} ,x_{y_2} ]_{y_2} =l^{M'}\leqno{(5.6.14)}$$
 where $y_2\in \Sigma_K$ is the prime 
divisor of 
$K$ lying over $z_2$,
and where the point   
$x_{y_2}$  in $(E(K_{y_2})/l^{M'}E(K_{y_2}))^{\epsilon }$ is
such that  $\partial_{l^{M'}} (x_{y_2})=e_{y_2}$.
Here in (5.6.14), $\delta_{M_1+M'}(z_1+z_2)$ is the cohomology class in $H^1(K,E)_{l^{M_1+M'}}$
associated to $z_1+z_2$.

Let 
$$\delta_{M_1}(z_1+z_2)$$ be the cohomology class of $\coprod\!\!\!\coprod$$(E/K)_{l^\infty}^{-\epsilon}$ associated  to the divisor
$z_1+z_2\in \Lambda^2(a)$, where this class belongs to the Tate-Shafarevich
group by Lemma 5.1.7(i).

 Then  the Cassels pairing gives, as $l^u\delta_{M_1}(z_1+z_2)=\delta_{M_1-u}(z_1+z_2)$ if $u\leq M_1$,
$$<\delta_{M_1}(z_1+z_2), l^uf>_{\rm Cassels}=
<\delta_{M_1-u}(z_1+z_2), f>_{\rm Cassels}. \leqno{(5.6.15)}$$
As $f$ has order $l^{M'}$ and that we have 
$$l^{M'}\delta_{M_1+M'-u}(z_1+z_2)=\delta_{M_1-u}(z_1+z_2),$$
by the construction of the Cassels pairing (Proposition 4.6.1 and equation (4.6.2)) we obtain from this last equation (5.6.15)
$$<\delta_{M_1}(z_1+z_2), l^uf>_{\rm Cassels}=\sum_{z\in {\rm \ Supp}(z_1+z_2)} [\delta_{M_1+M'-u}(z_1+z_2)_{y} ,x_{y} ]_{y} \leqno{(5.6.16)}$$
where $$g_{y} =\partial_{l^{M'}}(x_{y} ),\leqno{(5.6.17)}$$
and where $$x_y\in ({}_{l^{M'}}E(K_y))^{-\epsilon}$$ and because we have that $\delta_{M_1+M-u}(z_1+z_2)_v=0$ for all places $v$ of $K$ not dividing $z_1+z_2$ by Proposition 4.5.1(i).
Here in (5.6.16) and (5.6.17), $z$ runs over the places $z_1,z_2$ and 
$y\in \Sigma_K$ runs over the place of $K$ above $z$, that is to say above
$z_1,z_2$.

The first term $[\delta_{M_1+M'-u}(z_1+z_2)_{y_1} ,x_{y_1} ]_{y_1}$
in the sum (5.6.16) is zero by (5.6.10) and (5.6.17).

By (5.6.14), 
the second term, an element  of ${\mathds Q}/{\mathds Z}$,
$$[\delta_{M_1+M'-u}(z_1+z_2)_{y_2} ,x_{y_2} ]_{y_2} $$ 
is non-zero for $0\leq u\leq M'-1$ and hence $\delta_{M_1}(z_1+z_2)$   has  order at  least $l^{M'}$. 
We obtain from this and (5.6.16) that the Cassels pairing 
$$<\delta_{M_1}(z_1+z_2), l^{u}f>_{\rm Cassels}=[\delta_{M_1+M'-u}(z_1+z_2)_{y_1} ,x_{y_1} ]_{y_1}$$
is non-zero for all $u$ such that $0\leq u\leq M'-1$. Hence this character 
$$s\mapsto <\delta_{M_1}(z_1+z_2), s>_{\rm Cassels},
\ \ \coprod\!\!\!\!\coprod(E/K)_{l^\infty}^{-\epsilon} \to {\mathds 
Q}/{\mathds Z}$$
has restriction to ${\mathds Z}f$ which generates the dual of this 
subgroup ${\mathds Z}f$ of order $l^{M'}$. 

As $f$ is any element of the abelian group 
$\coprod\!\!\!\coprod$$(E/K)_{l^\infty}^{-\epsilon} $, the non-degeneracy of the Cassels pairing on $\coprod\!\!\!\coprod$$(E/K)_{l^\infty}^{-\epsilon} $  implies that the classes $\delta_{M_1}(z_1+z_2)$, for divisors $z_1+z_2\in \Lambda^2(a)$, generate  the $-\epsilon$ eigenspace $\coprod\!\!\!\coprod$$(E/K)_{l^\infty}^{-\epsilon} $ which proves the  part of the theorem
for $\coprod\!\!\!\coprod$$(E/K)_{l^\infty}^{-\epsilon} $.  ${\sqcap \!\!\!\!\sqcup}$

\vskip0.4in\vfil\eject
\noindent {\bf References}

\noindent [1] Brown M.L., Heegner Modules and Elliptic Curves. Springer Lecture Notes in Mathematics No. 1849. Springer Verlag, Berlin-Heidelberg-New York, 2004.

\noindent [2] Gross  B.H., Kolyvagin's work on modular elliptic curves. In: $L$-functions and Arithmetic (eds. J.H. Coates and M.J. Taylor) Cambridge University Press 1990, pp. 235-256.

\noindent [3]  Gross B.H.,  Zagier D., Heegner points and derivatives of $L$-series. Invent. Math. {\bf 84}, 225-320 (1986)

\noindent  [4]  Kolyvagin V.A., Euler systems. In: The Grothendieck Festschrift (Vol. II). P. Cartier et al., eds., pp. 435-483. Progress in Mathematics 8.   Birkh\"auser, Boston 1990

\noindent [5]  Kolyvagin V.A., Finiteness of $E({\mathds Q})$ and $\coprod\!\!\!\coprod$$(E/{\mathds Q})$ for a class of Weil curves. Math. USSR Vol.  {\bf 32} (1989) No. 3.

\noindent [6]  Kolyvagin V.A., On the structure of Shafarevich-Tate groups, Proceedings of USA-USSR Symposium on Algebraic Geometry, Chicago 1989.  Lecture Notes in Mathematics No. 1479, Springer Verlag, Berlin-Heidelberg-New York  1991

\noindent [7]  Kolyvagin V.A., On the structure of Selmer groups, Math. Annalen {\bf 291}, 253-259 (1991)

\noindent [8]  Mccallum W.G., Kolyvagin's work on Shafarevich-Tate groups. 
In: $L$-functions and Arithmetic (eds. J.H. Coates and M.J. Taylor) Cambridge University Press 1990, pp. 295-316.

\noindent [9]  Milne J.S., Arithmetic Duality Theorems, Perspectives in Mathematics, Academic Press, Boston 1986

\noindent [10] Raynaud, M., Caract\'eristique d'Euler-Poincar\'e d'un faisceau et
cohomologie des vari\'et\'es ab\'eliennes. In: Dix Expos\'es sur la Cohomologie des Sch\'emas. Ed. A. Grothendieck, pp. 12-30. North Holland Publ. Comp. Amsterdam 1967.

\end{document}